\documentclass[aps,pre,preprint,groupedaddress,showpacs]{revtex4-2}
\usepackage{graphicx,amsmath,amssymb}
\fontfamily{qhv}\selectfont
\usepackage{capt-of}
\usepackage[font=sf]{caption}

\begin{document}

\title{A GHOST PERTURBATION SCHEME TO SOLVE ORDINARY DIFFERENTIAL EQUATIONS }

\author{P.L. Garrido}
\email[]{garrido@onsager.ugr.es}
\affiliation{Instituto Carlos I de F{\'\i}sica Te{\'o}rica y Computacional. Universidad de Granada. E-18071 Granada. Spain }

\date{\today}
\begin{abstract}
We propose an algebraic method that finds a sequence of functions that exponentially approach the solution of any second-order ordinary differential equation (ODE) with any boundary conditions. We define an extended ODE (eODE) composed of a linear generic differential operator that depends on free parameters, $p$, plus an $\epsilon$ perturbation formed by the original ODE minus the same linear term. After the eODE's formal $\epsilon$ expansion of the solution, we can solve order by order a hierarchy of linear ODEs and we get a sequence of functions $y_n(x;\epsilon,p)$ where $n$ indicates the number of terms that we keep in the $\epsilon$-expansion. We fix the parameters to the optimal values $p^*(n)$ by minimizing a distance function of $y_n$ to the ODE's solution, $y$,  over a given  $x$-interval.
We see that the eODE's perturbative solution converges exponentially fast in $n$ to the ODE solution when $\epsilon=1$:  $\vert y_n(x;\epsilon=1,p^*(n))-y(x)\vert<C\delta^{n+1}$ with $\delta<1$. The method permits knowing the number of solutions for Boundary Value Problems just by looking at the number of minima of the distance function at each order in $n$, $p^{*,\alpha}(n)$, where each $\alpha$ defines a sequence of functions $y_n$ that converges to one of the ODE's solutions. We present the method by its application to several cases where we discuss its properties, benefits and shortcomings and some practical algorithmic improvements on it.
\end{abstract}
\pacs{18-3e}
\maketitle
\section{Introduction}
Physics is a thriving part of science because it manages that theories and experiments concur together to understand Nature. In the last centuries, we have built several fundamental theories such as Thermodynamics, Mechanics (Classical and Quantum), and Relativity that describe and predict with precision some of the variate behaviour of most of the systems we chose to scrutiny.  Typically, we apply the appropriate theory to some simplified ideal models of the real system to solve it exactly or approach them with some extra assumptions. In this way, we extract many valuable generic properties of such systems. Moreover, linear perturbation schemes to our initial ideal model permit us to go systematically a little beyond by introducing some realistic features to it. This scheme has been very successful for many years. However, it has been harder and harder to get meaningful theoretical results associated with phenomena far beyond linear approximations that are the ones that captive our primary attention. To get some insight into such cases, we use numerical approximated methods to directly solve the corresponding equations associated with the problem and/or look for some strong theoretical extra-assumptions to characterise the phenomena.

The author's primary interest is in systems at non-equilibrium states, an exciting field where the situation described above is ubiquitous \cite{Garrido0}. There are well known set of differential equations that describe, for instance, a fluid, a chemical reaction, or the dynamics of a set of neurons  \cite{Batchelor}. Those equations are complex, and we can get some straightforward solutions only for some simple academic cases or by linearizing the equations near a known solution that typically corresponds to the equilibrium case.  However, the far from equilibrium phenomena are out of our theoretical reach.  One typical example is the stationary Fourier's Law for heat transfer . It describes the temperature behaviour of a system characterised by a thermal conductivity coefficient $\kappa$ (that depends on the intrinsic properties of the material and, typically on the local temperature) that it extends in a spatial domain $\Lambda$  with a set of boundary conditions. The stationary Fourier's law is written by the partial differential equation:
\begin{equation}
\nabla\cdot\left(\kappa(T)\nabla T\right)=0\quad,\quad T(x_0)=\tilde T(x_0)\quad \forall x_0\in\partial \Lambda
\end{equation}
with a given $\tilde T(x)$-function.
This equation has a simple ``ideal'' solution when the temperature at the boundaries are all equal: $\tilde T(x)=T_0$ that implies $T(x)=T_0$ $\forall x\in\Lambda$. In some cases we can go beyond this ideal solution. For instance, let us assume that our system is a three dimensional cubic box of side $L$ with temperatures $T_0$ and $T_1$ at the faces at $z=0$ and $z=L$ respectively and periodic boundary conditions on the rest of the faces. The Fourier's equation is reduced to a one-dimensional ordinary differential equation (ODE):
\begin{equation}
\frac{d}{dz}\left[\kappa(T)\frac{dT(z)}{dz}\right]=0\quad T(0)=T_0\quad,\quad T(L)=T_1
\end{equation}
with the implicit formal solution:
\begin{equation}
\int_{T_0}^{T(z)}dT \kappa(T)=\frac{Jx}{L}\quad,\quad J=\int_{T_0}^{T_1}dT \kappa(T)
\end{equation}
When the temperature difference is very small, say $\Delta T\equiv T_1-T_0\simeq 0$ we can linearize Fourier's equation and the solution is just the linear profile: $T(z)=T_0+\Delta T z/L$ and the heat current across the system is $J=\kappa(T_0)\Delta T/L$. Moreover, we can systematically go beyond this linear solution order by order in powers of $\Delta T$. Let us remark a couple of things for this well-known case: (1) In general, one cannot obtain formal solutions of Fourier's law with non-symmetric domains and/or boundary distribution of temperatures. (2) In any case, we can apply a perturbative scheme because there is a natural parameter on the problem and a known reference's solution. A different situation occurs in developing interdisciplinary fields such as Biology or  Social Sciences. There,  the corresponding theories are also based on non-linear differential equations. However, frequently they do not have parameters that permit us to define (at least) a linear description around some reference state due to their intrinsic complex nature.
 And (3) except for simple cases (for instance, $\kappa(T)=cte$ in our example or see academic examples in ref.\cite{Bird}), it is necessary to use numerical tools to get some insight into the solution of the equation.

Things become theoretically more challenging if the differential equations that define our system's behaviour are just the deterministic stationary part of a fluctuating Langevin equation. Moreover, we may be interested, for instance, in the spatial correlations of our state variables at the stationary state. In the context of the Fourier's Law, such a situation is contemplated by the Fluctuating Hydrodynamics Theory \cite{Sengers}.
From a theoretical point of view, we need to know first the stationary solution. Then, we should solve a functional equation where the unknown is the correlation (see, for instance, \cite{Garrido}). As we already commented, we cannot typically obtain the stationary solution analytically. Therefore, from the beginning, we are restricted to doing, if possible, perturbation expansions around a reference state (near the equilibrium in the case of Fourier's Law). In all the other instances, the unique way to get a flavour of the correlation's behaviour is by doing direct simulations of microscopic models.  

In this context, we think it is of great interest to look for any algebraic way to approach the solutions of those classes of non-linear differential equations and look for a  perturbation scheme that connects an exactly solvable case with the real non-linear one. Once we have settled on our primary goal, we focus in this paper on solving generic non-linear second order differential equations as the first step before going to other ODEs, partial differential equations or Langevin equations.

Many powerful and successful algorithms solve ODEs that mainly depend on derivative's discretisations, and the application of an iterative method \cite{odesol} (see, for instance, the software package BVPSolve that solves Boundary Value Problems (BVP) \cite{bvp} with methods adapted to different situations). However, their results are not suitable for use in other parts of the theoretical effort to describe a problem. Nevertheless, some methods attempt to get analytical approximations to the solutions without any use of discretisation tools: from Taylor expansions, WKB type of approximations or variational methods (see, for instance, reference \cite{rev} to find some bibliography to get a general view of the different strategies). Let us focus on a few methods that initially fit our general strategy and are related to our proposal.

The first technique we want to recall is the {\it Adomian decomposition method} \cite{Adomian0} where the original non-linear differential equation $N(x)y=0$  is decomposed into two pieces: a linear part $L(x)y$ that it is easely invertible (for instance $L(x)y=y''$) and the rest $N(x)y-L(x)y$. Moreover, it is introduced a formal parameter $\epsilon$ in such a way that the differential equation becomes:
\begin{equation}
N(x)y=0\quad\Rightarrow\quad L(x)y+\epsilon\left[N(x)y-L(x)y \right]=0 \label{3}
\end{equation}
where the new extended ODE (eODE) maintains the original boundary conditions. When $\epsilon=1$ we recover trivially the original ODE. Then, it is assumed the existence of a formal decomposition of the solution:
\begin{equation}
y(x;\epsilon)=\sum_{m=0}^\infty \tilde y_m(x)\epsilon^m\label{1}
\end{equation}
that it is used to expand eODE in $\epsilon$-powers including any nonlinear term. That permits the algebraic computation of the coefficients $\tilde y_n(x)$ order by order in $\epsilon$. Finally one makes $\epsilon=1$ and it is expected that $y_n(x)=\sum_{m=1}^n \tilde y_m(x)$  converge uniformly to the solution $y(x)$ when $n\rightarrow\infty$. This method has been widely used with some success. Aside from that, there is no rigorous proof on its premises. The main problem with this method is that there is no way to control the convergence rate to the solution or even its existence. Nevertheless, it contains two helpful ideas: a formal solution expansion that permits dealing with any nonlinearity in the ODE and the observed sequence's uniform convergence in many different examples. Another criticism of this technique is about the computational cost to do the $\epsilon$-expansion of the nonlinear term. This fact can be dismissed by the use of algebraic mathematical software such as MATHEMATICA or MAPLE, where easily we can do expansions of such terms up to $\epsilon^{50}$ in a few seconds on an ordinary laptop.

The second technique we want to comment is the Homotopy method \cite{He1}. In this case the extended ODE is written:
\begin{equation}
(1-\epsilon)\left[L(x)y-L(x)y_0\right]+\epsilon N(x)y=0
\end{equation}
where $L(x)$ is a linear operator and $y_0(x)$ is a convenient chosen zeroth order guess function. From here, the method follows the same path as the Adomian decomposition. We assume an epsilon expansion of the solution and solve order by order in epsilon the corresponding ODE to find the $m$-th expansion term $\tilde y_m(x)$. The main problem with this method is that it depends strongly on the initial guess $y_0(x)$ in such a way that a wrong choice may lead to divergent solutions. An optimisation of the Homotopy method has been proposed to control and enhance the expansion's convergence rate to the solution:
\begin{equation}
(1-\epsilon)L(x)y+H(\epsilon;p) N(x)y=0\quad ,\quad H(\epsilon;p)=\epsilon\sum_{n=0}^\infty p_n\epsilon^n
\end{equation}
where the $p_n$ are constants. After expanding in $\epsilon$ the differential equation, we get order by order $\tilde y_n(x;p)$ that now depend on the parameters $p$. The parameters are found by minimizing the functional
\begin{equation}
d(p)=\int_{a}^b dx (N(x) y(x;p))^2\quad ,\quad y(x;p)=\sum_{m=0}^\infty \tilde y_m(x;p) \label{2}
\end{equation}
The minimisation procedure for an infinite number of parameters is one of the principal difficulties of this method when going to high orders in the expansion. Nevertheless, this last method includes two more exciting ideas: first, it is unnecessary to use a guess function $y_0(x)$.  And second, we can get better results by including extra parameters. They are fixed by minimising a residual function $d(p)$ that controls the distance to the exact solution. 

We present in this paper a way to generate approximations to the solution of any second order's ODE that includes, in our opinion, some of the most interesting elements of the above methods. Our scheme rests on the following items:
\begin{itemize}
\item (1) The eODE is similar to the Adomian method (\ref{3}): Linear operator plus an epsilon nonlinear deviation.
\item (2) The Adomian decomposition of the solution (\ref{1}) and the $\epsilon$ expansion of the eODE. 
\item (3) The inclusion of parameters that are fixed by minimizing a distance to the exact value (\ref{2}).
\end{itemize} 
The distinctive part of our method is that we assume a generic linear operator that may depend on four parameters (see eq.(\ref{linear})).
We'll show that in this way, the approximate solution of the eODE, $y_n(x;\epsilon,p)$ is highly sensitive to the parameter values. That is,  small changes on some $p$-values imply significant changes on $y_n$ for any given $x$ and $\epsilon$. This property can be considered unwanted, but it is beneficial to our goals. It allows the sequence $y_n$ to adapt exponentially fast to the solution as $n$ increases, and it efficiently detects the existence of several solutions in some boundary value problems. 

The paper is structured in the following form. In Section II we define the ODE, the eODE and the $\epsilon$-expansion. We also write down the recurrence to obtain the coefficients $\tilde y_n(x)$. We also define two possible ways to measure the distance between the $n$-th approximation to the exact unknown solution. Afterward, we state the conjectures we expect our method to accomplish, which we will prove in concrete examples in the paper. The first conjecture states that the distance of $y_n$ to the exact solution as a function of the parameters $p$ has several local minima whose structure is maintained with $n$. Each of them corresponds to a sequence de $y_n$ that fits one of the possible solutions of the ODE.  This property is very relevant in the case of boundary value problems where the number of solutions is a priori not known. Finally, we introduce the ghost expansion concept. We use the exponential convergence rate of the sequence $y_n$ to define a perturbative expansion for the solution that can be used as in other theoretical studies.  

The following sections of the paper are designed to show the method's properties, conjectures, and some other aspects by studying concrete examples. Section III is devoted to apply our method to the BVP $\xi y''-y=0$ with $y(0)=1$ and $y(1)=0$ where the exact solution is known. We show step by step how are computed $\tilde y_n$ up to order $n=40$. We analyse how some distances as (\ref{2}) are related to the distance to the exact solution and how such distance decays exponentially fast. We build the corresponding Ghost Expansion and comment on some convergence properties of our Adomian expansion. Section IV studies the Bratu differential equation, a BVP exactly soluble with two solutions. We show how our method detects both solutions, and we introduce a way to accelerate the convergence once we know a good approximation to any solution. We also study the use of  $\epsilon$ as a minimisation parameter. We show that the overall convergence rate increases and $\epsilon\rightarrow 1$ as we increase $n$. Section V is devoted to study the BVP $y''+\xi(y'+y^2)=0$ with $y(0)=0$ and $y(1)=1$ where no analytical solution is known. We show that only one solution exists  whenever $\xi<\xi_c=3.7681..$. Section VI introduce the application of the method to an Initial Value Problem (IVP): the Lane-Emdem equation $y''+2y'/x+y^m=0$ with $y(0)=1$ and $y'(0)=0$. Exact solutions for this ODE are known for $m=0,1$ and $5$. We observe the exponential decay rate to the exact solution for a given $x$-interval $T$. We show that the number of the expansion terms that we need to reach a given precision grows logarithmic with $T$  Finally, we propose a way to extend our method to large $T$ values without losing precision. 

\section{The Method}
Let us introduce the second-order Ordinary Differential Equation (ODE):
\begin{equation}
N(x)y\equiv g(x,y,y')y''+h(x,y,y')=0\label{ode}
\end{equation}
where $g$ and $h$ are, in principle, well-defined analytic functions on the domain where the solutions, if any, exist and $N(x)$ is the formal non-linear differential operator associated to this ODE. The boundary conditions determine the solution's existence and their properties. In this paper we study two of them:
\begin{itemize}
\item {\it (1) Boundary Value Problem (BVP):} $y(0)=\bar y_0$, $y(1)=\bar y_1$, with $x\in[0,1]$.
\item {\it (2) Initial Value Problem (IVP):} $y(0)=\bar y_0$, $y'(0)=\bar y_1$, with $x\geq 0$.
\end{itemize}
The goal of this work is to design an algebraic perturbative method capable of approximating (with arbitrary precision) the solution or solutions of the ODE. (\ref{ode}). The method is based in three main ingredients: (1) The extension of the ODE (eODE) by adding a finite number of parameters $p\equiv(p_0,p_1,\ldots)$ and a perturbative one, $\epsilon$.  This extension should be such that when $\epsilon\rightarrow 1$ or/and $p\rightarrow 0$ we recover the original ODE. Moreover, eODE should have an algebraic solution when $\epsilon\rightarrow 0$. (2) The introduction of a perturbative expansion on the eODE around $\epsilon=0$ and the obtention of and algebraic formal solution at all orders in $\epsilon$. (3) The definition of measures based on the ODE that define some kind of distance of a function to the real solution. 

\subsection{The extended ODE (eODE) and its perturbation expansion}
Let us define the {\it extended ODE (eODE)} as a quasi-linear second order differential operator $L$ plus a non-linear correction $N-L$:
\begin{equation}
L(x;p)y+\epsilon \left(N(x)y-L(x;p)y\right)=0\label{eode}
\end{equation}
where $p=(p_0,p_1,p_2,p_3)$ are arbitrary parameters and
\begin{equation}
L(x;p)y=p_0\frac{d^2y}{dx^2}+p_1\frac{dy}{dx}+p_2 y+p_3\label{linear}
\end{equation}
and $N(x)$ is the differential operator defined in eq.(\ref{ode}) that represents our original ODE.
We assume that the eODE have the same boundary conditions as the original problem (BVP or IVP). Let $y(x)$ and $y(x;p,\epsilon)$ be solutions of the ODE (\ref{ode}) and eODE (\ref{eode}) respectively. Let us assume, for simplicity, that they have the limiting properties:
\begin{eqnarray}
&&\text{(a)} \lim_{p\rightarrow 0} y(x;\epsilon,p)=y(x)\quad\forall\,\epsilon\nonumber\\
&&\text{(b)}  \lim_{\epsilon\rightarrow 1} y(x;\epsilon,p)=y(x)\quad\forall\, p\nonumber\\
&&\text{(c)} \lim_{\epsilon\rightarrow 0} y(x;\epsilon,p)=y_0(x;p)\qquad\text{with}\quad L(x;p)y_0=0\nonumber\\
&&\text{(d)} \lim_{s\rightarrow \infty} y(x;\epsilon,s p)=y_0(x;p)\quad\forall\,\epsilon\label{cond}
\end{eqnarray}
That is, we are assuming a nice, regular behavior of the solutions on their definition's domain.  Observe that we recover the solution of the ODE in two limits: $p\rightarrow 0$ and $\epsilon\rightarrow 1$.

We now define a perturbative expansion around $\epsilon=0$. Let us assume that the solution/s of the eODE can be written:
\begin{equation}
y(x;\epsilon,p)=\sum_{k=0}^\infty \tilde y_k(x;p)\epsilon^n\label{y1pert}
\end{equation} 
Let also introduce the $n$-th approximation as:
\begin{equation}
y_n(x;\epsilon,p)=\sum_{k=0}^n \tilde y_k(x;p)\epsilon^k\label{y1pert2}
\end{equation} 

Therefore we can expand the $g$ and $h$ functions on (\ref{eode}):
\begin{eqnarray}
g(x,y(x;\epsilon,p),y'(x;\epsilon,p))&=&\sum_{n=0}^\infty g_n(x;p)\epsilon^n\nonumber\\
h(x,y(x;\epsilon,p),y'(x;\epsilon,p))&=&\sum_{n=0}^\infty h_n(x;p)\epsilon^n\label{y2pert}
\end{eqnarray}
Observe that the dependence on $x$ of the coefficients $g_n$ and $h_n$ is explicitly on $x$ or through the functions $\tilde y_0(x;p)$,...,$\tilde y_n(x,p)$, $\tilde y'_0(x;p)$, $...$, $\tilde y'_n(x;p)$. We substitute eqs.(\ref{y1pert}) and (\ref{y2pert}) into (\ref{eode}) and we get a hierarchy of closed equations order by order in $\epsilon$ that allows us to determine the unknowns $\tilde y_n(x;p)$:
\begin{equation}
p_0 \tilde y''_n+p_1 \tilde y'_n+p_2 \tilde y_n=F_n(x;p)\label{basicODE}
\end{equation}
where
\begin{eqnarray}
n=0: F_0(x;p)&=&-p_3 \nonumber\\
n>0: F_n(x;p)&=&p_0\tilde y''_{n-1}-h_{n-1}-\sum_{l=0}^{n-1}g_{n-1-l}\,\tilde y''_l+p_1 \tilde y'_{n-1}+p_2\tilde y_{n-1}+p_3\delta_{n,1}
\end{eqnarray}
We do not write the function's arguments to simplify the notation. We observe that this eODE structure permits an algebraic solution order by order in $\epsilon$ for any boundary conditions. It is a matter of very simple analysis to find the general solution for eq.(\ref{basicODE}) with BVP and IVP boundary conditions: 
\begin{itemize}
\item BVP: $y(0)=\bar y_0$, $y(1)=\bar y_1$. 
\begin{eqnarray}
\tilde y_n(x;p)&=&\frac{1}{p_0w(e^w-1)}\biggl[e^{w_+x}\biggl(p_0w(\bar y_1^{(n)}e^{-w_-}-\bar y_0^{(n)})+\int_0^1 du e^{-w_-u}F_n(u;p)\nonumber\\
&+&(e^w-1)\int_0^x due^{-w_+u}F_n(u;p)-e^w\int_0^1du e^{-w_+u}F_n(u;p)\biggr)\nonumber\\
&+&e^{w_-x}\biggl(p_0w(\bar y_0^{(n)}e^{w}-\bar y_1^{(n)} e^{-w_-})-\int_0^1 du e^{-w_-u}F_n(u;p)\nonumber\\
&-&(e^w-1)\int_0^x du e^{-w_-u}F_n(u;p)+e^w\int_0^1du e^{-w_+u}F_n(u;p)\biggr)\biggr]\label{solBVP}
\end{eqnarray}

\item IVP: $y(0)=\bar y_0$, $y'(0)=\bar y_1$.
\begin{eqnarray}
\tilde y_n(x;p)&=&\frac{1}{p_0w}e^{w_+x}\left(p_0(\bar y_1^{(n)}-w_-\bar y_0^{(n)})+\int_0^x du e^{-w_+ u}F_n(u;p)\right)\nonumber\\
&+&\frac{1}{p_0w}e^{w_-x}\left(p_0(w_+\bar y_0^{(n)}-\bar y_1^{(n)})-\int_0^x du e^{-w_- u}F_n(u;p)\right)\label{solIVP}
\end{eqnarray}
\end{itemize}
Where $\bar y_0^{(0)}=\bar y_0$,  $\bar y_1^{(0)}=\bar y_1$ and $\bar y_0^{(n)}=\bar y_1^{(n)}=0$ for $n>0$. Moreover, $w_\pm=(-p_1\pm(p_1^2-4p_0p_2)^{1/2})/(2p_0)$ and $w=w_+-w_-$. 

These equations define the perturbation expansion entirely in both cases. Observe that the boundary conditions are included order by order naturally.

\subsection{Ways to measure the distance to the solutions}
We may expect that $y_n(x;\epsilon,p)$ to be an approximation of a ODE's solution $y(x)$ at least when $\epsilon\rightarrow 1$, and/or $p\rightarrow 0$ consistently with the property (\ref{cond}). Moreover, we hope for the likeness with the exact solution to increase with the perturbation order. In any case, we need some way to measure the distance between our perturbative solution and the {\it unknown} real solution. We have studied two distances:
\begin{eqnarray} 
d_1(n;\epsilon,p)&=&\biggl\{\int_0^1 dx \biggl[g(x,y_n(x;\epsilon,p),y'_n(x;\epsilon,p))y_n''(x;\epsilon,p)\nonumber\\
&+&h(x,y_n(x;\epsilon,p),y'_n(x;\epsilon,p))\biggr]^2\biggr\}^{1/2}\label{mea1}\\
d_2(n;\epsilon,p)&=&\left[\int_0^1 dx \left[y_n(x;\epsilon,p)-\bar y_n(x;\epsilon,p)\right]^2\right]^{1/2}\label{mea2}
\end{eqnarray}
where $\bar y_n(x;\epsilon,p)$ is solution of the differential equation:
\begin{equation}
\bar y''_n=-\frac{h(x,y_n(x;\epsilon,p),y'_n(x;\epsilon,p)))}{g(x,y_n(x;\epsilon,p),y'_n(x;\epsilon,p)))}\equiv G_n(x;\epsilon,p)
\end{equation}
with the corresponding original ODE's boundary conditions. 
Its solution is a particular case of eq.(\ref{basicODE}) with $p_0=1$, $p_1=p_2=p_3=0$ and $F_n(x;p)=G_n(x;\epsilon,p)$. That is, for BVP we get:
\begin{equation}
\bar y_n(x;\epsilon,p)=\bar y_0+x(\bar y_1-\bar y_0)+\int_0^x du (x-u)G_n(x;\epsilon,p)-x\int_0^1du (1-u)G_n(x;\epsilon,p)
\end{equation}
and for IVP:
\begin{equation}
\bar y_n(x;\epsilon,p)=\bar y_0+x \bar y_1+\int_0^x du (x-u)G_n(x;\epsilon,p)
\end{equation}
$d_1$ mesures the average deviation of $y_n$ to be locally the ODE's solution and $d_2$ measures the mistmach between the approximation $y_n$ and the integrated result $\bar y_n$. Let us remark that both distances are equal to zero when the approximate solution, $y_n(x;\epsilon,p)$, is equal to the exact one $y(x)$.

\subsection{The scheme to get increasingly good approximations to the ODE's solutions}

The general scheme we present here has been derived after studying several examples like those we expose later in the paper. Therefore, this section is just an effort to propose a generalized set of well-defined conjectures that resume the behaviors we have seen in particular problems. We hope they are confirmed by rigorous works or application to other cases. 

One of the most relevant ingredients in our method was to introduce the set of parameters $p$ in the quasi-linear operator $L$ (\ref{linear}). We found that, at each perturbative level, $n$, we could improve the solution $y_n(x;\epsilon,p)$ by choosing the set of parameters $(\epsilon^*,p^*)$ that minimize any of the distances $d_1$ or $d_2$. Moreover, we immediately realized that the improvement with $n$ was very fast in all the studied cases. Therefore, we claim that, in general, the following conjectures may be true:

{\bf\boldmath \underline{(Strong) Conjecture:} For any $n>n_0>0$ and given a distance $d(n;\epsilon,p)$ (for instance $d_1$ or $d_2$ defined in (\ref{mea1}) and (\ref{mea2})), there exists a {\it l(n)}-set of parameter values $\{\epsilon_k^*(n),p_k^*(n)\}_{k=1}^{l(n)}$ that are local minima for $d(n;\epsilon,p)$ such that
\begin{equation}
\lim_{n\rightarrow\infty}d(n;\epsilon_k^*(n),p_k^*(n))=0 \quad,\quad k=1,\ldots,l(n)
\end{equation}
}
We will see on the examples that $\epsilon_k^*(n)\rightarrow 1$ as $n\rightarrow\infty$ that is coherent with property (b) on eq.(\ref{cond}). Therefore, we could have fixed $\epsilon=1$ from the begining and only minimize with respect the parameters $p$. The convergence is, in this case, a bit slower, but we think that the conjecture still applies and it can be written:

{\bf\boldmath \underline{(Restricted) Conjecture:} For any $n>n_0>0$ there exists {\it l(n)}-set of parameter values $\{p_k^*(n)\}_{k=1}^{l(n)}$ that are local minima for a given measure ($d_1$ or $d_2$) with $\epsilon=1$ such that
\begin{eqnarray}
\lim_{n\rightarrow\infty}d(n;\epsilon=1,p_k^*(n))=0 \quad k=1,\ldots,l(n)
\end{eqnarray}
}
That is, the $n$th order approximation of the algebraic solution of the eODE converges to the ODE solution when we tune the value of the parameters $p(n)$ by the ones that minimize the distance defined in (\ref{mea1}) or (\ref{mea2}).  

Let us make more precise comments on the behaviors we have found in all the examples we have studied:
\begin{itemize}
\item The number of local minima of $d(n;\epsilon,p)$ may depend on $n$. 
\item The set of minima for all the $n$'s  have an overall structure. For instance, we may have two minima when $n$ is odd and one minima when $n$ is even. Therefore, the minima from the even values define a sequence, and the first and second minima from the odd $n$-values another two sequences respectively. In general, we assume that a given ODE have a minima structure with periodicity $\bar k$. Therefore, all the expansion order $n=\bar k\bar n+k$ for each $k\in[1,\bar k]$ have the same number of minima, $\bar l(k)$. We can define each sequence by two numbers: $(k,l^{(k)})$ where  we have ordered the minima in some way, $l^{(k)}\in[1,\bar l(k)]$:
$$c(k,l^{(k)})\equiv\left\{(\epsilon_{l^{(k)}}(\bar k\bar n+k),p_{l^{(k)}}^*(\bar k\bar n+k)\right\}_{\bar n=0}^\infty$$ 
Therefore, K=$\cup_{k=1}^{\bar k}\bar l(k)$ is the total number of different sequences.
\item A Boundary Value Problem may have $s$-solutions $\left\{y^{(j)}(x)\right\}_{j=1}^s$ and it is expected that $s\leq K$. Therefore more than one sequence may converge to the same solution and all solutions are assumed to be described:
\begin{equation}
\lim_{\bar n\rightarrow\infty}y_{\bar k\bar n+k}(x;c(k,l^{(k)})_{\bar n})=y^{(j)}(x)\quad j\in[1,\ldots,s]\quad \forall (k,l^{(k)})
\end{equation}
\item {\bf The Ghost Expansion:}  Let us define the distance of one element $\bar n$ of a given sequence $(k,l^{(k)})$: $d^*(\bar n;(k,l^{(k)}))=d(\bar n;c(k,l^{(k)}))$. Typically (at least on the examples we shown in this paper) $d^*(\bar n;(k,l^{(k)}))\simeq \delta(k,l^{(k)})^{\bar n}$ when $\bar n>>1$ with $\delta<1$. That is, the approximation converges uniformly on $x$ and exponentialy fast with $n$ to the ODE's solution. Therefore, we can write:
\begin{equation}
y_{\bar k\bar n+k}(x;c(k,l^{(k)})_{\bar n})=\sum_{m=0}^{\bar n} w_m(x;(k,l^{(k)})) d^*(m;(k,l^{(k)}))\label{ghost1}
\end{equation}
where, by construction, 
\begin{eqnarray}
w_m(x;(k,l^{(k)}))&=&\frac{y_{\bar k m+k}(x;c(k,l^{(k)})_{m})-y_{\bar k (m-1)+k}(x;c(k,l^{(k)})_{m-1})}{d^*(m,(k,l^{(k)}))}\quad m>0\nonumber\\
w_0(x;(k,l^{(k)}))&=&\frac{y_{k}(x;c(k,l^{(k)})_{0})}{d^*(0,(k,l^{(k)}))}\label{ghost2}
\end{eqnarray}
where $\lim_{m\rightarrow\infty}w_m(x,(k,l^{(k)}))=w(x,(k,l^{(k)}))$ is of order one. Therefore, the ODE's $j$-solution can be naturally expanded using any sequence $(k,l^{(k)})$ that converges to it:
\begin{equation}
y^{(j)}(x)=\sum_{m=0}^\infty w_m(x;(k,l^{(k)})) d^*(m;(k,l^{(k)}))
\end{equation}
We call this the {\bf Ghost Expansion} of each ODE's solution. Observe that the perturbative parameter, $d^*$, depends on the boundary conditions, the ODE's structure, and the distance used in a highly non-trivial way. Let us stress that we have found a way to create a perturbative expansion of the ODE's solution that can be helpful when dealing with theories lacking intrinsic natural perturbative parameters.
\end{itemize}

Let us study a few examples of ODEs to detail how the method can be applied. In this way we will discuss the properties of the $\epsilon$-expansion, how the distances $d_{1,2}$ behave and several other interesting issues as the existence of solutions in BVP.

\section{Example 1: \boldmath $\xi y''-y=0$ (BVP)}

Let us to illustrate the method with the BVP of one of the simplest differential equation:   
\begin{equation}
\xi y''-y=0\label{odep1}
\end{equation}
with $\xi>0$ and $y(0)=1$, $y(1)=0$. The solution of this ODE is:
 \begin{equation}
y(x)=\frac{e^{-x/\sqrt{\xi}}-e^{(x-2)/\sqrt{\xi}}}{1-e^{-2/\sqrt{\xi}}}\label{odep1e}
\end{equation}
That will help us to analyze the perturbative expansion (\ref{y1pert}) and its convergence behavior to the exact solution. Moreover, we will check how the distances $d_{1,2}$ given by eqs. (\ref{mea1},\ref{mea2}) work. 

\begin{itemize}

\item {\bf eODE's perturbative expansion:}
The ODE (\ref{odep1}) is represented in our notation by $g(x,y,y')=\xi$ and $h(x,y,y')=-y$. 
We restrict our general parametric extended ODE (\ref{eode}) to the  $p_0\neq 0$ and $p_1=p_2=p_3=0$ case. This particular choice gives already excellent results, and its simplicity allows us to show the method neatly. The eODE is given by:
\begin{equation}
\left[p_0+\epsilon\left(\xi-p_0\right)\right]y''-\epsilon y=0\label{eode2}
\end{equation}
The perturbative expresions for $g$ and $h$ in eq. (\ref{y2pert}) are in this case:
\begin{equation}
g_n(x;p_0)=\xi \delta_{n,0}\quad,\quad h_n(x;p_0)=-\tilde y_n(x;p_0)
\end{equation}
The differential equations to be solved order by order are 
\begin{equation}
p_0 \tilde y''_n=F_n(x;p_0)
\end{equation}
where
\begin{equation}
F_0(x;p_0)=0\quad ,\quad F_n(x;p_0)=(p_0-\xi)\tilde y''_{n-1}+\tilde y_{n-1}\quad n>0
\end{equation}
In particular, the first three orders are:
\begin{eqnarray}
p_0 \tilde y''_0&=&0\quad (\tilde y_0(0)=1, \tilde y_0(1)=0)\Rightarrow \tilde y_0(x;p_0)=1-x\nonumber\\
p_0 \tilde y''_1&=&1-x  \quad (\tilde y_1(0)=0, \tilde y_1(1)=0)\Rightarrow \tilde y_1(x;p_0)=-\frac{1}{6p_0} x(1-x)(2-x) \nonumber\\
p_0 \tilde y''_2&=&\frac{1-x}{6p_0}\left[6(p_0-\xi)-x(2-x)\right]\quad (\tilde y_2(0)=0, \tilde y_2(1)=0)\Rightarrow\nonumber\\
 &&\tilde y_2(x;p_0)=-\frac{1}{360p_0^2}x(1-x)(2-x)(-4+60(p_0-\xi)-6x+3x^2)\label{expan}
\end{eqnarray}
We use Mathematica software to iterate the process up to 40 perturbative terms and we fix $\xi=1/10$ for the rest of the analysis. The algebraic solutions permit precise control of our method and favor its detailed analysis.  
\begin{center}
\vskip -1cm
\includegraphics[height=9cm]{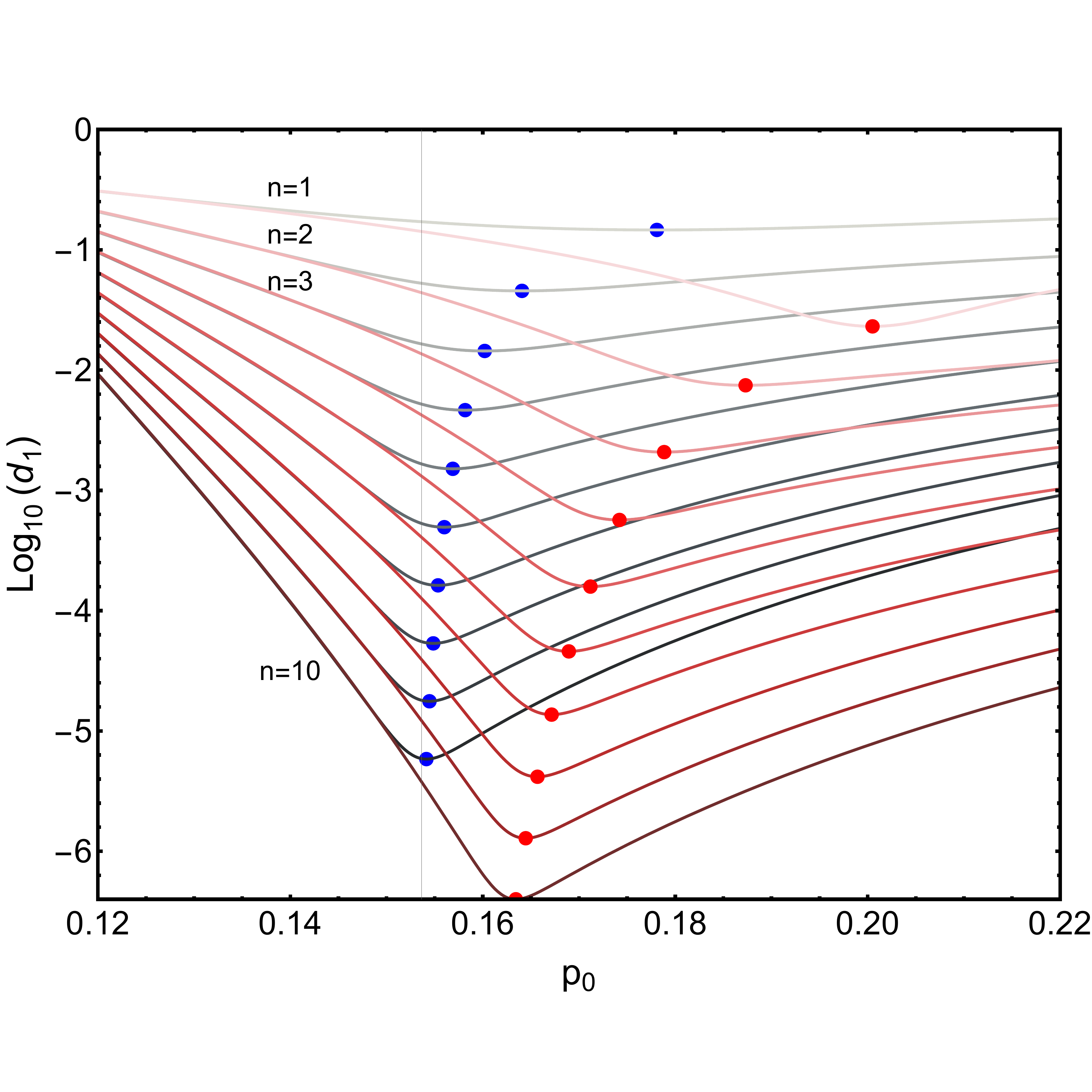}  
\captionof{figure}{Example 1, ODE eq.(\ref{odep1}). Decimal logarithm of distances $d_{1,2}(n;\epsilon=1,p_0)$ defined by eqs. (\ref{mea1}) and (\ref{mea2}) versus $p_0$ for each $n$th-perturbation approximation. The $y_n$ are obtained from the eODE's perturbative expansion (\ref{eode2}) for $\xi=1/10$. Gray-Black curves are for $d_1$ and Red-Cherry curves are for $d_2$. From top to bottom $n=1, 2,\ldots,10$. Blue and Red dots are the corresponding minima.}\label{fig1} 
\end{center}

\item {\bf \boldmath The distances $d_{1,2}(n;\epsilon=1,p_0)$ (Restricted Conjecture):} We see in figure \ref{fig1} the behavior of the $\log_{10}(d_{1,2}(n;\epsilon=1,p_0))$ as a function of $p_0$ for some values of $n$. The curves present a well defined minimum at each $n$-perturbative level. Moreover, the minimums, $p_0(n;\alpha)^{*}$ with $\alpha=1,2$ are located always in sharp and narrow basins around them. In fact, just moving the minimum by $0.02$, the distances $d_{1,2}$ increase by two orders of magnitude. Therefore, the minima are, by far, the optimal values to get the minimum distance to the exact result at each order of the perturbative expansion. Let $d_{\alpha}^*(n)\equiv d_{\alpha}(n,\epsilon=1,p_0(n;\alpha)^*)$. We can now address some open questions: Is $d_{\alpha}^*(n)$, $\alpha=1,2$, a well defined distance to the exact solution? and What are the differences between $d_1$ and $d_2$ defined above? Let us give some insight on those issues.
\begin{center}
\vskip -1cm
\includegraphics[height=9cm]{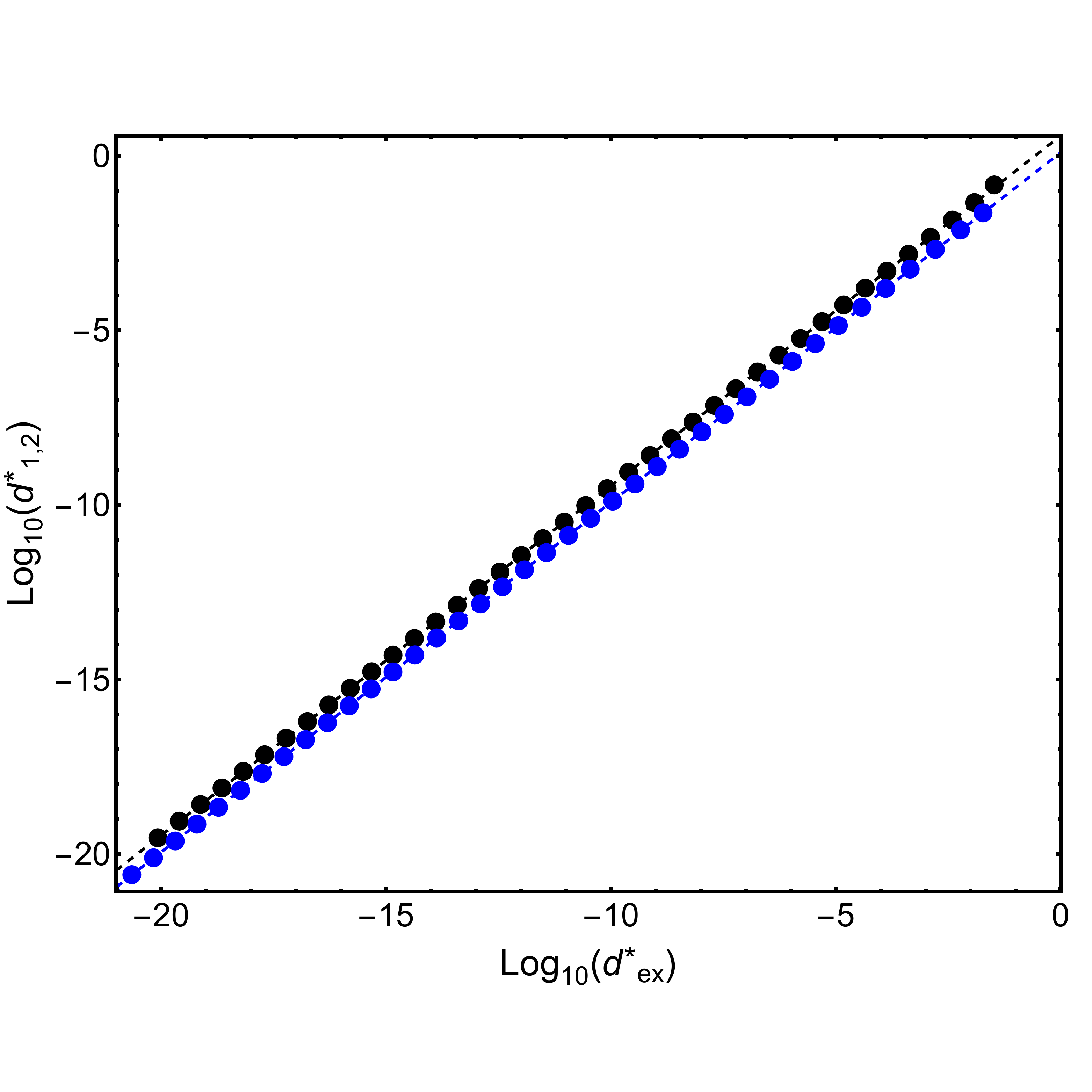}  
\vskip -1cm
\captionof{figure}{Example 1, ODE eq.(\ref{odep1}). Pairs of data: $(\log_{10}d_{ex}^*(n),\log_{10}d_{\alpha}^*(n))$ for $n=1,\ldots,40$ and $\alpha=1,2$ black and blue dots respectively. Dotted lines are linear fits to the data (see main text).}\label{fig2}  
\end{center}
 In this example we know that the exact ODE's solution is $y(x)$ given by eq.(\ref{odep1e}). We define the distance of any function $\tilde y(x)$ to it by:
\begin{equation}
d_{ex}[\tilde y]=\left[\int_0^1 dx \left[\tilde y(x)-y(x)\right]^2\right]^{1/2}\label{meae}\\
\end{equation}
In particular the distance of $y_n^*(x;\alpha)\equiv y_n(x;\epsilon=1,p_0(n;\alpha)^{*})$ to the exact solution is $d_{ex}^*(n;\alpha)\equiv d_{ex}[y_n^*(\alpha)]$ where remind that $\alpha=1,2$ stands for the two different measures we have defined. We show in Figure \ref{fig2} the set of pairs $(\log_{10}d_{\alpha}^*(n),\log_{10}d_{ex}^*)$ for $n=1$ up to $40$ and for both distances. The best fit we have found is the linear one: $\log_{10}d_\alpha^*=a_\alpha+b_\alpha\log_{10}d_{ex}^*$ with $a_1=0.57(0.01)$, $b_1=1.002(0.001)$, $a_2=0.090(0.007)$ and $b_2=1.0016(0.0006)$. That is, the distances $d_{1,2}$ are almost proportional to the distance to the exact solution $d_{ex}^*$: $d_{\alpha}^*(n)=C_\alpha (d_{ex}^*(n))^{b_\alpha}$ at least in the interval studied $n\in[1,40]$. Moreover, we see that $d_1^*(n)\simeq 3d_2^*(n)$. Observe that each distance selects different values for the optimal parameter for the same approximate algebraic solution $y_n$. Nevertheless, they converge exponentially fast to the ODE's solution.

From this example, it seems convenient to use the distance $d_2$ instead of $d_1$ because it gives a better approximation to the ODE's solution at each order in $n$. However, that is not so evident from a computational point of view. $d_1$ depends on derivatives and $d_2$ on integrals, and each may have different computational speeds. Therefore, the distance having the shorter computational time to reach a given precision  may depend on each problem. In any case, we may conclude that the distances $d_{1,2}$ (defined by Eqs. (\ref{mea1}) and (\ref{mea2})) suitable measures of the distance to the exact solution, and we can use both just by paying attention to their computational efficiency.
\begin{center}
\vskip -0.5cm
\includegraphics[height=5cm]{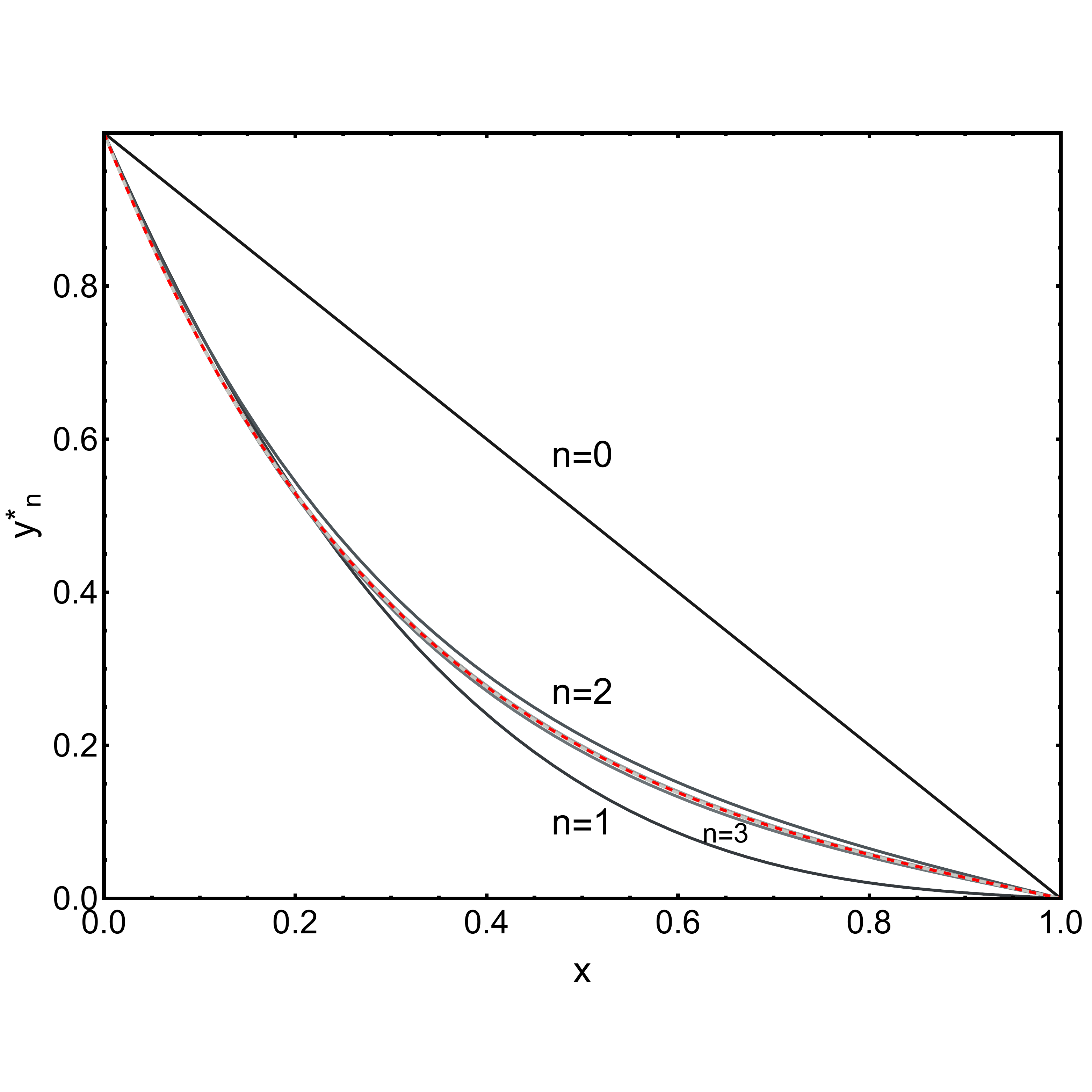}  
\includegraphics[height=5cm]{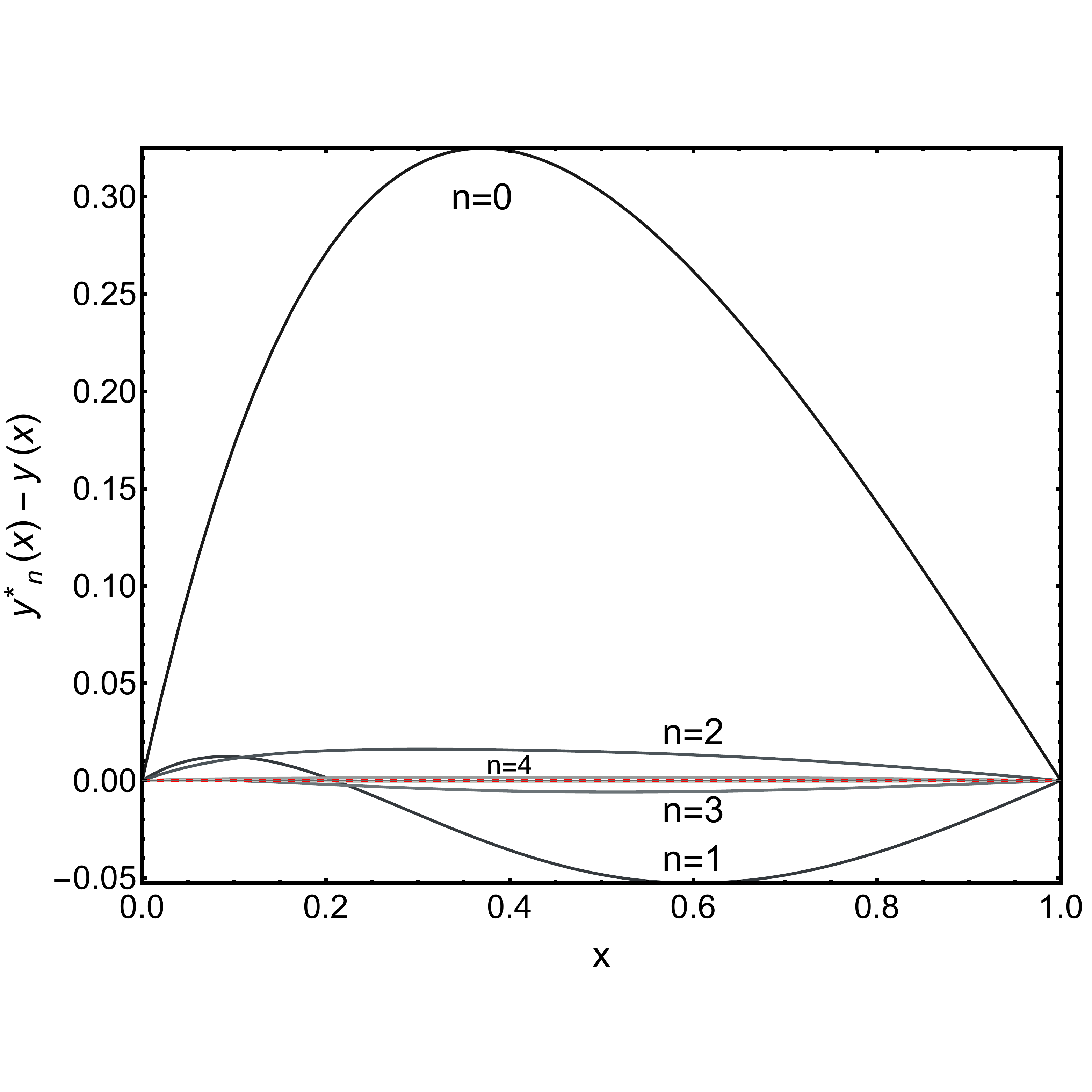}
\includegraphics[height=5cm]{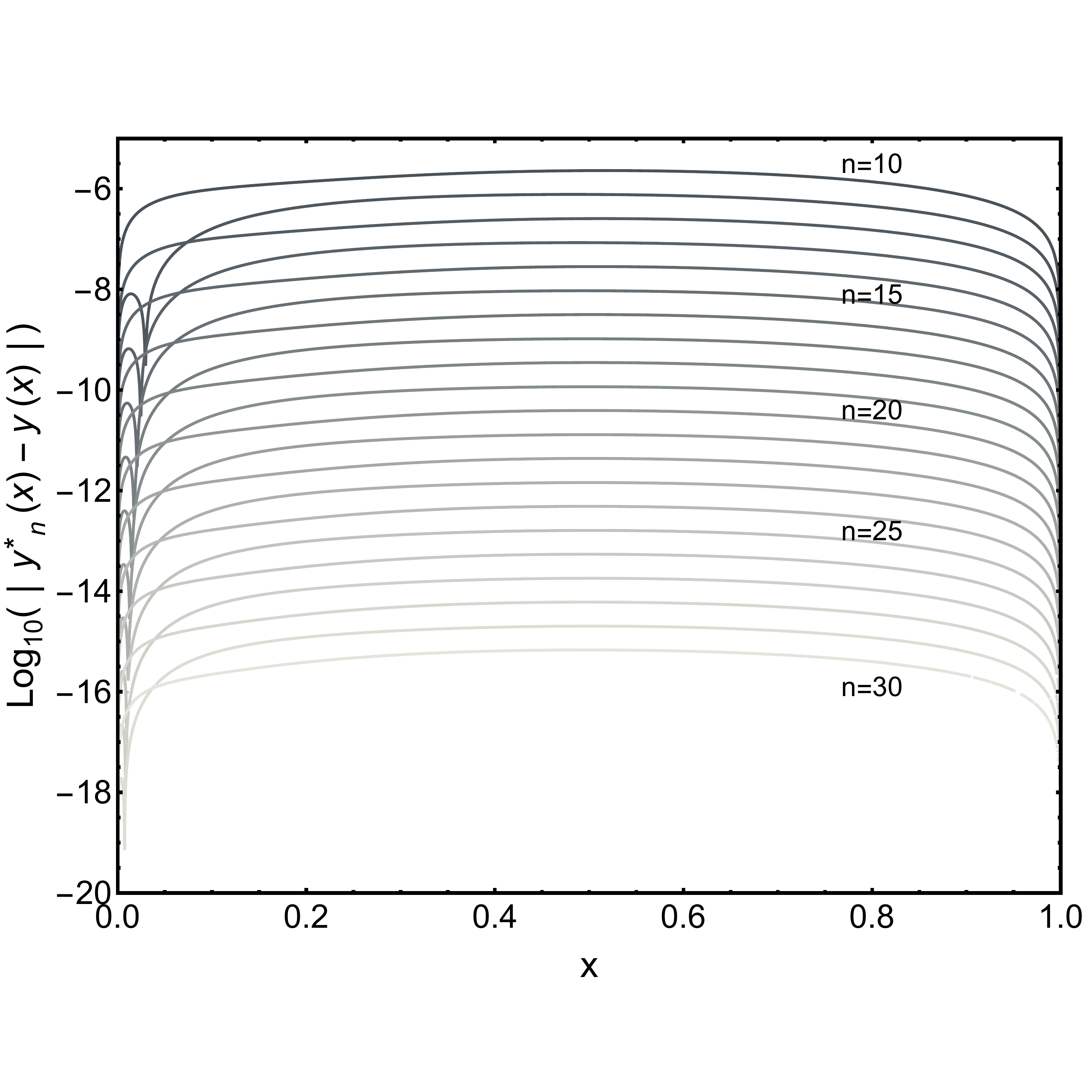}
\vskip -0.5cm
\captionof{figure}{Example 1, ODE eq.(\ref{odep1}). Behavior of the several $n$-th approximation $y_n^*(x;\alpha=1)\equiv y_n(x;\epsilon=1,p_0(n;\alpha=1)^{*})$ using the distance $d_1$. Left figure: $y_n^*(x)$ vs $x$ for $n=0, 1,\ldots 10$. Central figure: $y_n^*(x)-y(x)$ vs. $x$. Right figure: $\log_{10}\vert y_n^*(x)-y(x)\vert$ vs. $x$ for $n=10, 11,\ldots 30$. Red dotted line is the exact solution.} \label{fig3}  
\end{center}
\item {\bf \boldmath  The approximate solutions $y_n^*(x)=y_n(x;\epsilon=1,p_0^*(n))$:} We see in Figure \ref{fig3} the visual convergence of the approximations $y_n^*$ compared with the exact solution for the sequence corresponding to the minima computed with the distance $d_1$. We observe how the convergence is very fast and for $n=4$ we have already a reasonable approximation. For $n=10$ each point is on average at $10^{-6}$ (compared with $1$) distance to the exact solution. The distance diminish regularly by one order of manitude from $n$ to $n+2$ and, for instance, it of order $10^{-16}$ for $n=30$.
\begin{center}
\vskip -0.5cm
\includegraphics[height=7cm]{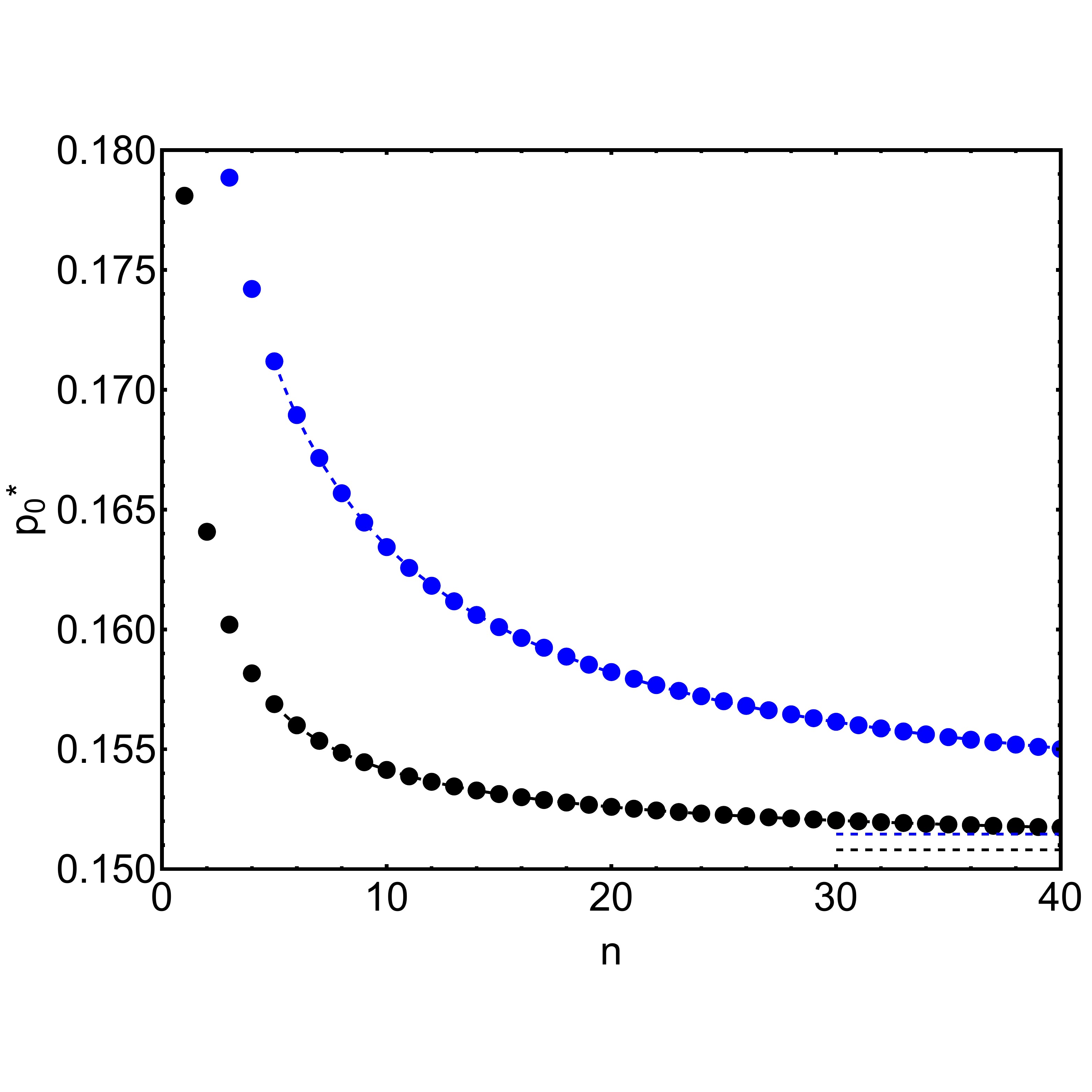}  
\includegraphics[height=7cm]{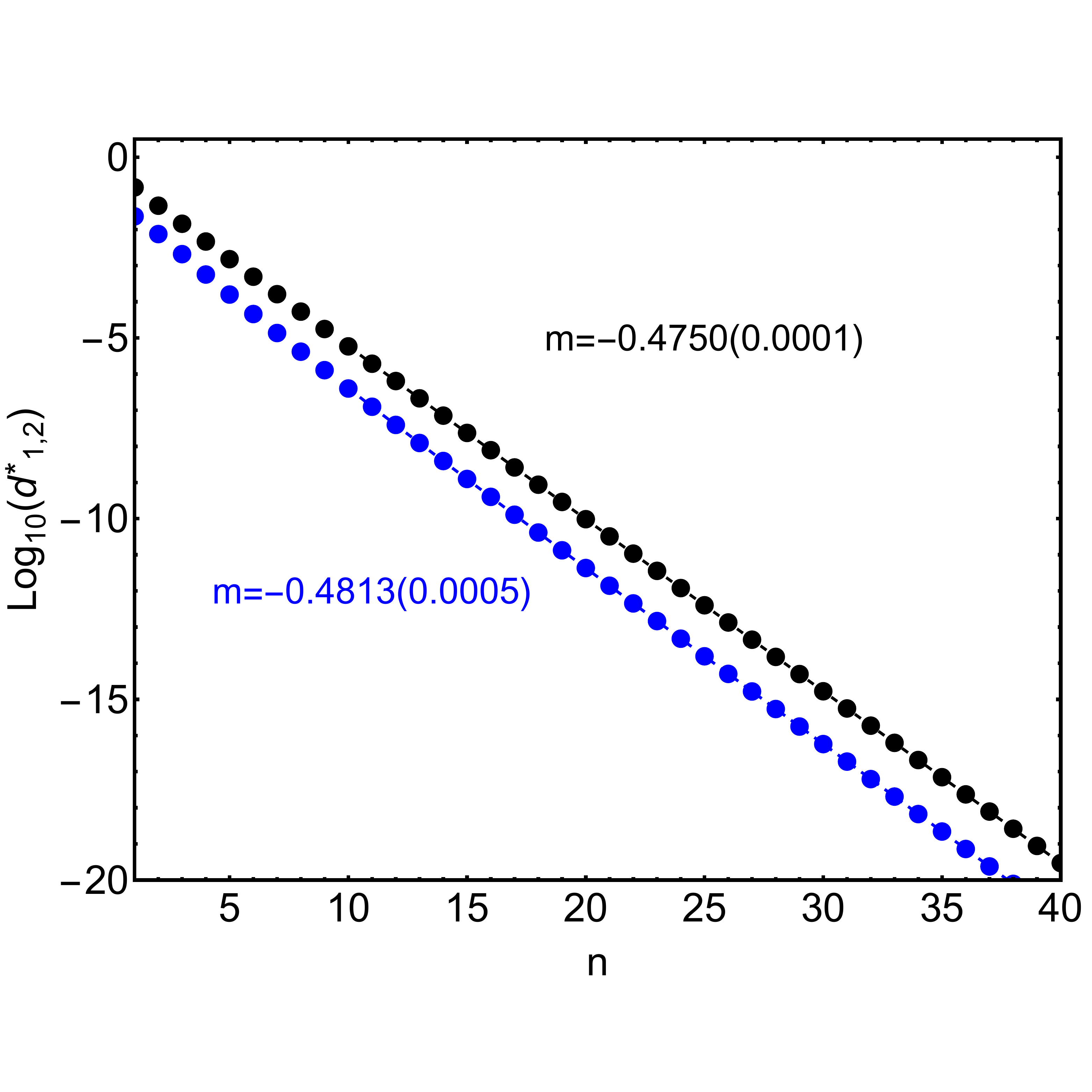}
\vskip -0.5cm
\captionof{figure}{Example 1, ODE eq.(\ref{odep1}). Asymptotic behavior with the $n$-th perturbative approximation of the values where the minima of $d_{1,2}$ are located, $p_0^*$ (left figure) and the values of such distances at the minima, $d_{1,2}^*$. Black and blue dots are for $d_1$ and $d_2$ respectively. Dashed lines are fits explained in the main text.} \label{fig4}  
\end{center}
\item {\bf \boldmath The asymptotic regime ($n\rightarrow\infty$):} We study how the the minima change with the $n$-th approximation. In figure \ref{fig4} (left) we see how the values of minima $p_0^*(n)$ decrease with  $n$ for the distance $d_1$ (black dots) and $d_2$ (blue dots). The best (simple) fit we have done to the data is a power-law decay (we tried exponential-type, but the fits were much worse): $p_0^*(n,\alpha)=a_0(\alpha)+a_1(\alpha)/n+a_2(\alpha)/n^2+a_3(\alpha)/n^3$ for all data in the interval $n\in[5,40]$. We get $a_0(1)=0.1508(0.0001)$, $a_0(2)=0.1515(0.0001)$, $a_1(1)=0.0387(0.0005)$, $a_1(2)=0.152(0.003)$, $\ldots$. Therefore, the numerical limiting values $p_0^*(\infty,\alpha)$ are $0.1508(0.0001)$ and $0.1515(0.0001)$ for $d_1$ and $d_2$ respectively. The small mismatch is probably due to the fitted function and/or to the need of data for much larger values of $n$. Let us stress that there is a finite asymptotic value for $p_0^*(n)$ that it is associated to the exact solution. We also see how the minima's depth, $d_{1,2}^*(n)$ behaves with $n$ and that we already know that it is related to the distance to the exact solution. In figure \ref{fig4} (right) we have plot the $\log_{10}d_{1,2}^*(n)$ vs. $n$. We have fit to the data a corrected linear behavior: $\log_{10}d_{\alpha}^*(n)=b_0(\alpha)n+b_1(\alpha)+b_2(\alpha)/n$ at the interval $n\in[10,40]$. We obtain $b_0(1)=-0.4750(0.0001)$ and $b_0(2)=-0.4813(0.0005)$. That is, both distances decay exponentially fast with $n$: 
\begin{equation}
d_{\alpha}^*(n)\simeq \delta(\alpha)^{n}\quad,\quad \delta(\alpha)=10^{-b_0(\alpha)}
\end{equation}
with $\delta(1)=0.3350..$ and $\delta(2)=0.3301..$ which is coherent with the small differentces found when we compared the distances with the distance to the exact solution (see above). Observe that the decay rates depend on the distance we use. 
\begin{center}
\vskip -1cm
\includegraphics[height=7cm]{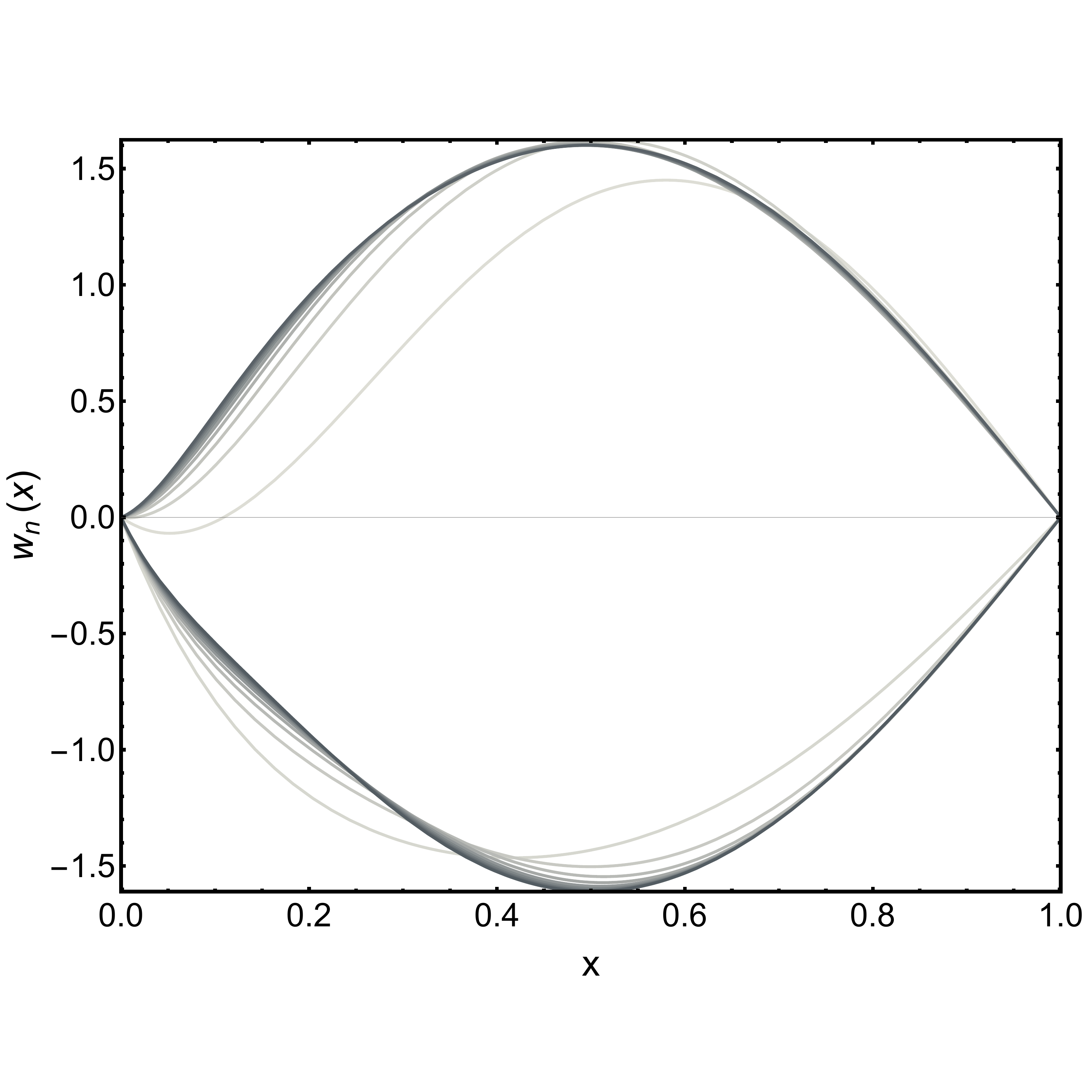}
\includegraphics[height=7cm]{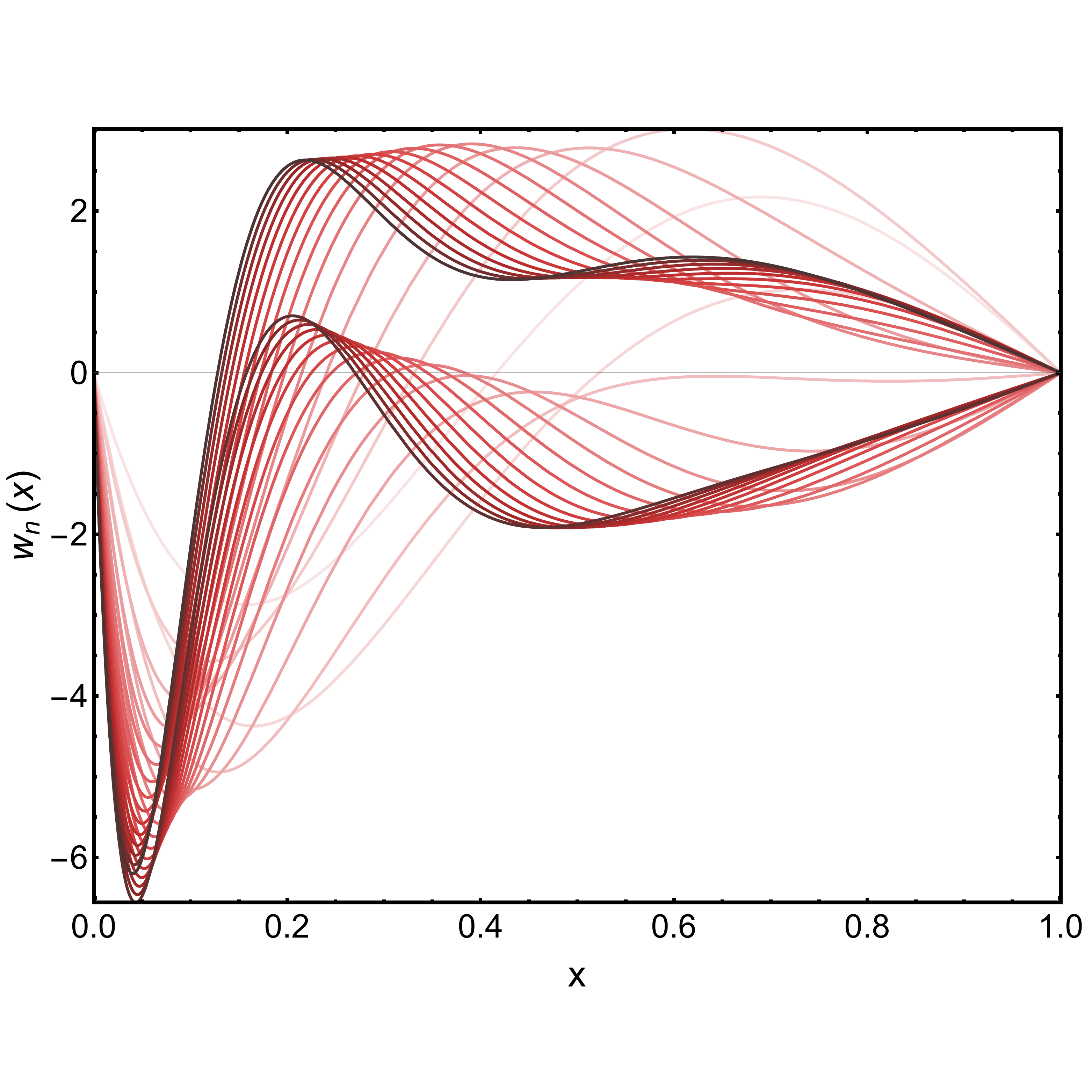}  
\vskip -1cm
\captionof{figure}{Example 1, ODE eq.(\ref{odep1}). Ghost expansion terms $w_n(x)$ defined in eq.(\ref{ghost2}) for distance $d_1$ (left) and $d_2$ (right). The darkness of gray and red colors increase with  $n\in[5,30]$}\label{fig5}  
\end{center}
\item {\bf The Ghost expansion:} The Ghost expansion is defined by eqs.(\ref{ghost1}) and (\ref{ghost2}). We show in figure \ref{fig5} how the functions $w_n(x)$ behave for the distances $d_1$ and $d_2$. We see that the such functions can be bounded by constants that do not depend on the perturbative level. Moreover, they seem to converge for large $n$-values to two limiting regular functions for $n$ even and odd values. Just to give a flavour about the form of such $w$-functions we show the first orders for $\xi=1/10$ and distance $d_1$:
\begin{eqnarray}
w_2(x)&=&(2-x) (1-x) x(-6.79498 x^2+13.59 x-1.40314)\nonumber\\
w_3(x)&=&(2-x) (1-x) x(-3.34634 x^4+13.3854 x^3-17.9727 x^2+9.17472 x-5.376)\nonumber\\
w_4(x)&=&(2-x) (1-x) x(-0.947722 x^6+5.68633 x^5-16.1041 x^4+26.5075 x^3\nonumber\\
&-&28.0748 x^2+18.2981 x-0.275838)\label{wn}
\end{eqnarray}
\begin{center}
\vskip -1cm
\includegraphics[height=7cm]{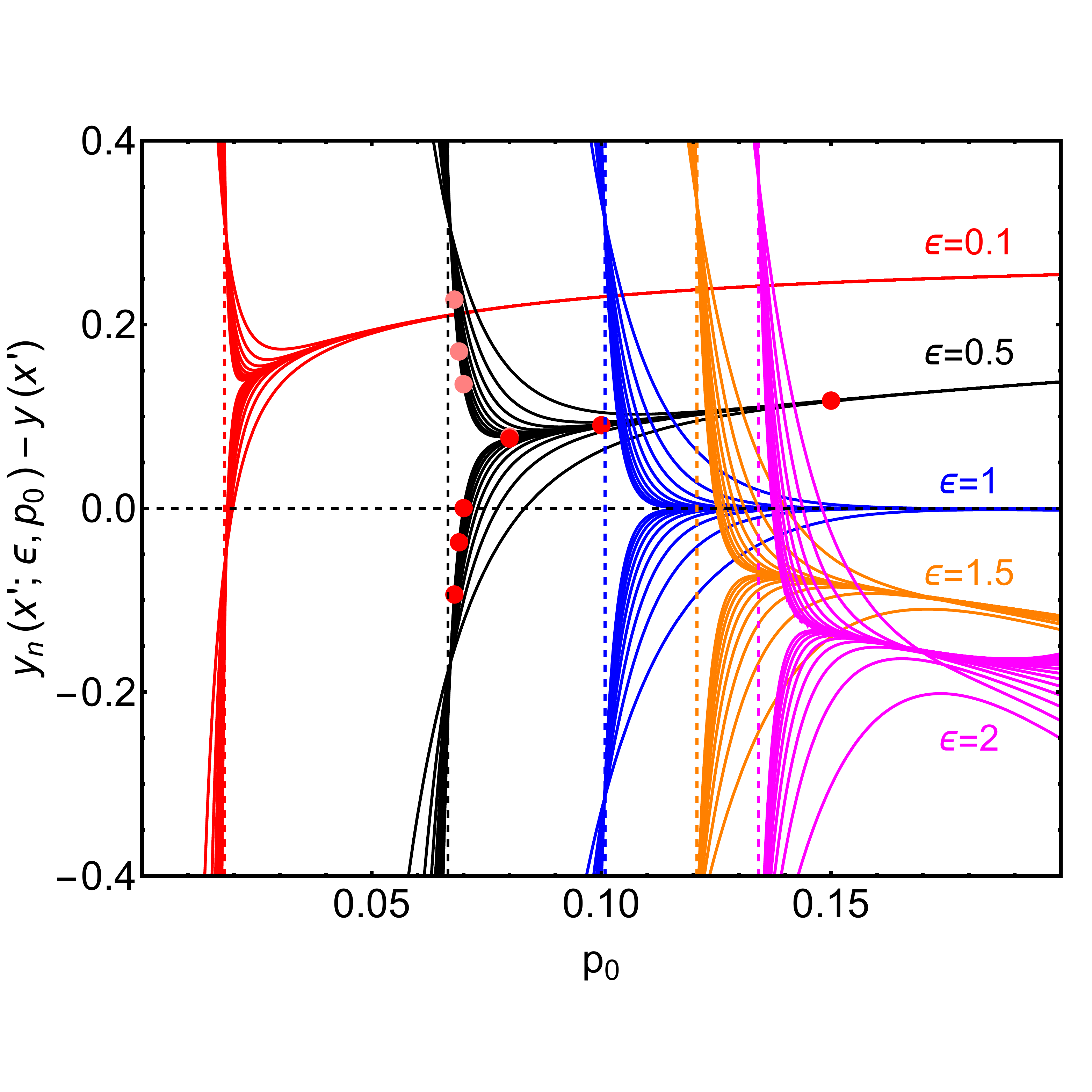}
\includegraphics[height=7cm]{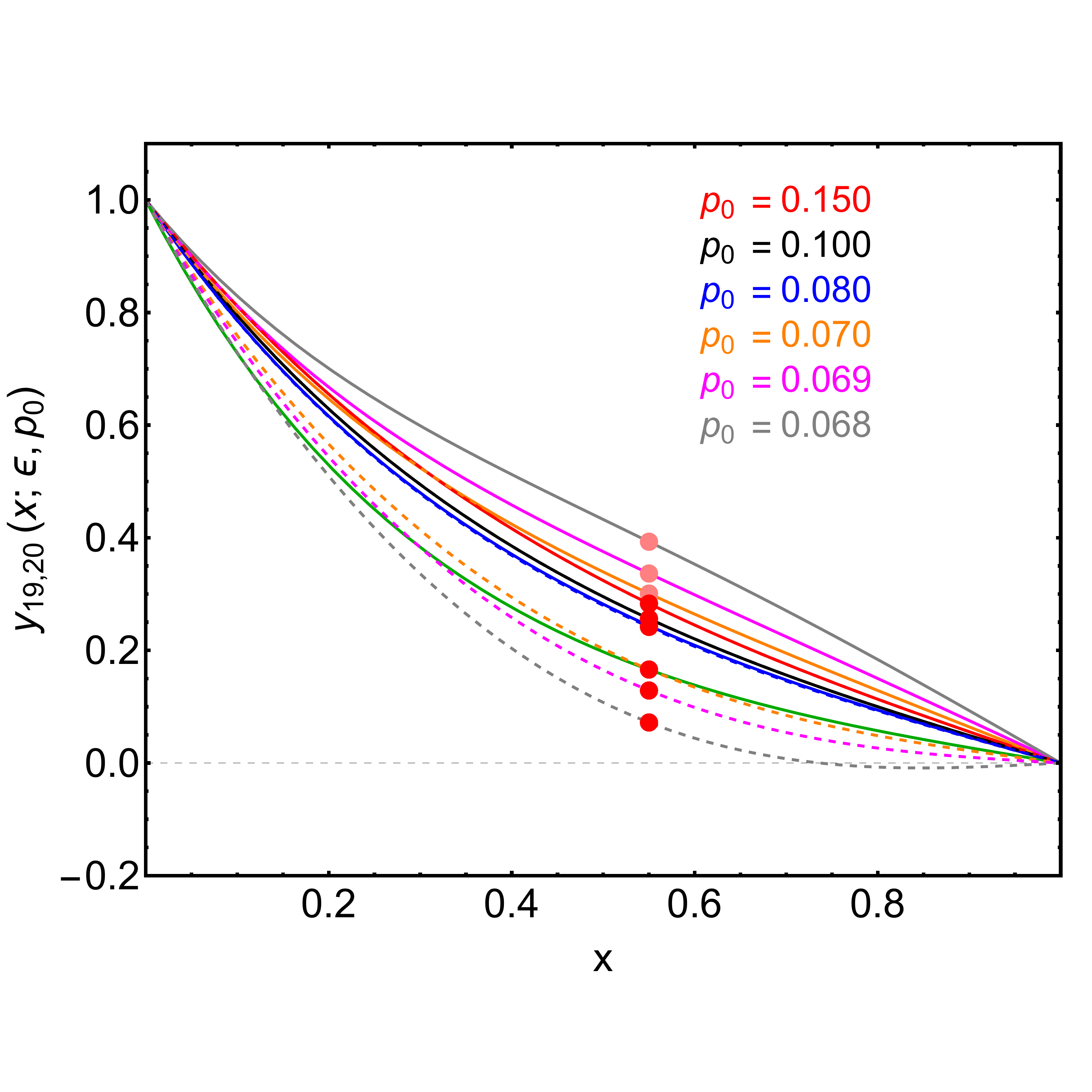}  
\vskip -1cm
\captionof{figure}{Example 1, ODE eq.(\ref{odep1}). Perturbative analysis of the BVP for the eODE given by eq. (\ref{eode2}) for $\xi=1/10$. Left: Difference between $n$'th approximate solution, $y_n$, and the exact one, $y$ given by eq. (\ref{odep1e}) at the point $x'=0.55$ as a function of $p_0$. There are shown the cases $\epsilon=0.1$, $0.5$, $1$, $1.5$, $2$ (red, black, blue, orange and magenta lines) for $n=3,\ldots,15$. Horizontal dotted line shows the zero value. Vertical dotted lines are the $p_0^c$ values such that below them the perturbative series is divergent. Red (Pink) dots represent particular values of $p_0$ (from left to right): $0.068$, $0.069$, $0.07$, $0.08$, $0.1$ and $0.15$ for $n=14(15)$. Right: Approximate perturbative solutions for  $n=19$ (solid lines)  and $n=20$ (dotted lines) for $\epsilon=0.5$ and $p_0$ taking values shown in the figure.  Red  and Pink dots are the values for $x'=0.55$ and correspond to the ones on the left figure.}\label{fig6}  
\end{center}
\item {\bf The eODE's $\epsilon$-expansion behavior:} We have used the method by fixing $\epsilon=1$ and looking for the minima of a given distance. It is interesting to show the general behavior of the original  perturbative expansion as a function of $p_0$ for a given $x$ and different $\epsilon$-values. Let us remind that, before our perturbative scheme, we expect that the solution of the eODE (\ref{eode}) will converge to the solution of the ODE (\ref{ode}) when $p_0\rightarrow 0$ $\forall\epsilon$. We show in figure \ref{fig6} (left) the difference of $y_n(x';\epsilon,p_0)$ with the known exact result given by eq.(\ref{odep1e}) at $x'=0.55$ and for $\epsilon=0.1, 0.5, 1, 1.5$ and $2$. We see the typical behavior of an asymptotic expansion:
\begin{itemize}
\item There is an apparent convergence towards the exact solution when $p_0$ moves towards zero from large initial values. For instance, we observe that for any $n$, $y_n(x'=0.55,\epsilon,p_0)-y(x'=0.55)$, decreases (when $\epsilon<1$) or increases (when $\epsilon>1$) apparently towards zero as we decrease $p_0$ from $0.2$. That occurs up to some value ($p_0\simeq 0.05$ for $\epsilon=0.1$ or $p_0\simeq 0.15$ for $\epsilon=0.5$ for example). That is the normal behavior of an analytical function having  a well defined limit at $p_0=0$.
\item When $p_0$ is small enough, the $p_0^{-1}$ terms dominate (see the expansion (\ref{expan})) for any $\epsilon$. Then, $y_n(x';\epsilon,p_0)-y(x')$  splits in two branches. The odd $n$-terms climb up to positive values and the even $n$-terms go down to negative ones quite abruptly. Therefore, it is imposible to reach the $p_0$ limit using these approximations. Observe that we could do some limit by sending $\epsilon\rightarrow 0$ and $p_0\rightarrow 0$ at the same time. We have checked such posibility by taking $\epsilon=c p_0$ and studying the series $\lim_{p_0\rightarrow 0}y_n(x;\epsilon=c p_0,p_0)=\bar y_n(x;c)$. We see that there is also an optimal $c$, but the distance to the exact solution doesn't go to zero as we increment $n$. Probably our perturbative expansion is too pathological to have a joint limit $p_0,\epsilon\rightarrow 0$ that converges to the exact solution. Nevertheless, it is something to be explored with rigorous tools.
\item The approximate solutions $y_n(x,\epsilon,p_0)$ are, typically far from the exact form for any $p_0$ and $\epsilon\neq 1$ (see for instance figure \ref{fig6} right). However, we see that the case $\epsilon=1$ is clearly singular. There is an interval of values $p_0$ where the distance to the exact solution is minimal (at least for a given $x'=0.55$ in the figure \ref{fig6} left). Moreover, such interval increases in size as we increment the approximation order. We checked that this behaviour also happens for any other $x'\in[0,1]$ value. This behaviour confirms that $\epsilon=1$ is the more convenient case to get approximations that are near the ODE's solution.
\begin{center}
\vskip -1cm
\includegraphics[height=7cm]{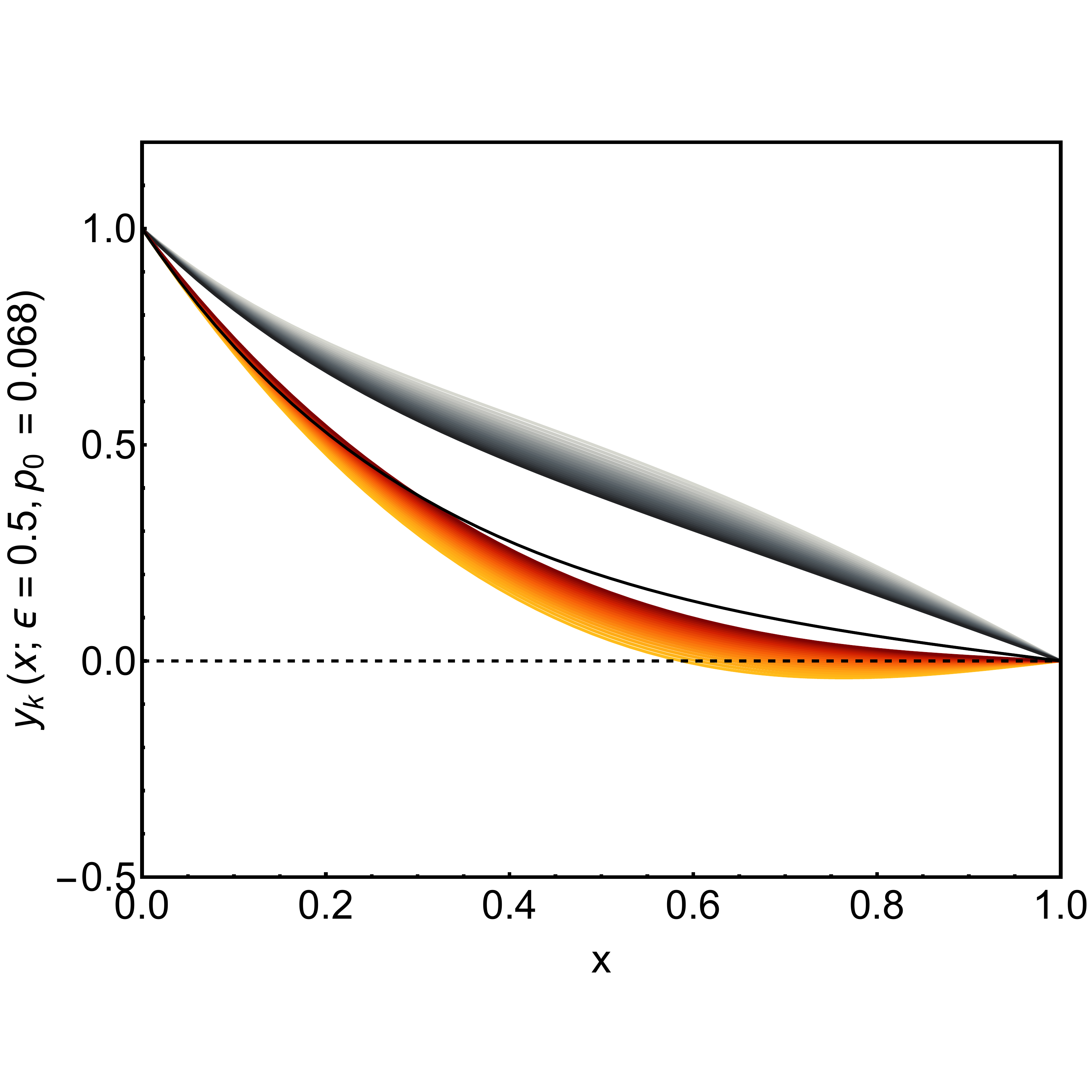}
\includegraphics[height=7cm]{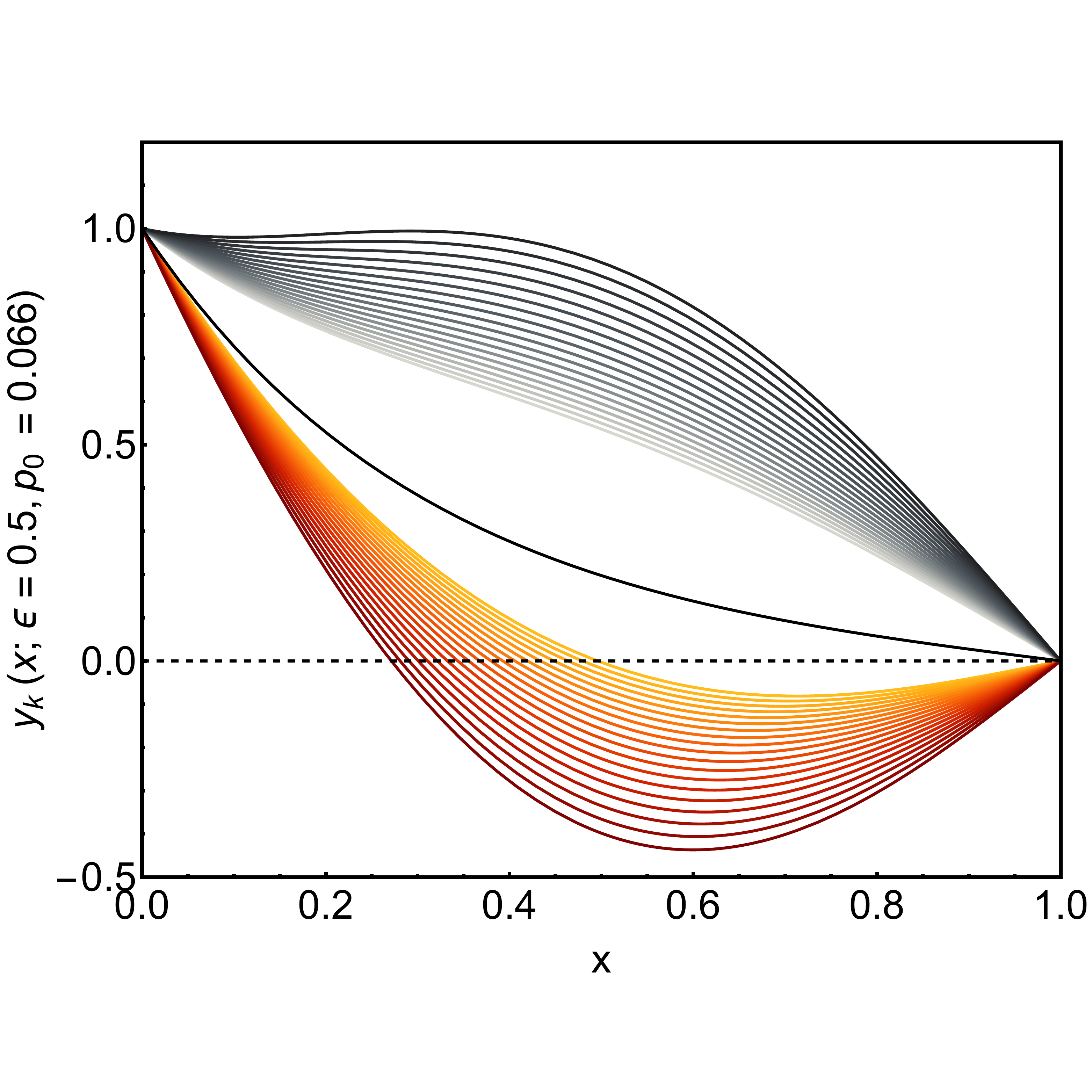}  
\vskip -1cm
\captionof{figure}{Example 1, ODE eq.(\ref{odep1}). Sequence of approximations $y_k(x;\epsilon=0.5,p_0)$ with $p_0=0.068$ (left) and $p_0=0.066$ (right) for $k=4,5,\ldots,40$. Black (orange) curves are for even (odd) $k$-values. The intensity of the color increases from small $k$-values to large ones. The solid black standalone curve corresponds to the known ODE solution.}\label{fig7}  
\end{center}
\item The sequence of functions $\left\{y_k(x;\epsilon,p_0)\right\}_{k=1}^n$  converges to a limiting function when  $n\rightarrow\infty$ whenever $p_0<p_0^c(\epsilon)$ and diverges otherwise. See for instance figure \ref{fig7} where a sequence of approximations are shown for $\epsilon=0.5$ and $p_0=0.068$ (left) and $p_0=0.066$ (right). We observe in the first case (left) how the approximation $y_k$ converges to two well-defined limits  for odd and even values of $k$ respectively. In contrast, there is a systematic distance's growth to the exact solution in the second case. In order to get a precise value of $p_0^c(\epsilon)$ we define the parameter 
\begin{equation}
s_n(\epsilon,p_0)=\left(\int_0^1dx \vert y_n(x;\epsilon,p_0)\vert^2\right)^{1/2} \label{quo}
\end{equation}
We compute $s_n(\epsilon,p_0)/s_{n-2}(\epsilon,p_0)$ for increasing values of $n$ for a given $\epsilon$ and $p_0$. When $y_n$ converges to a finite curve, such quotient tends to one. If it diverges, it approaches a constant greater than one.
\end{itemize}
\begin{center}
\vskip -1cm
\includegraphics[height=7cm]{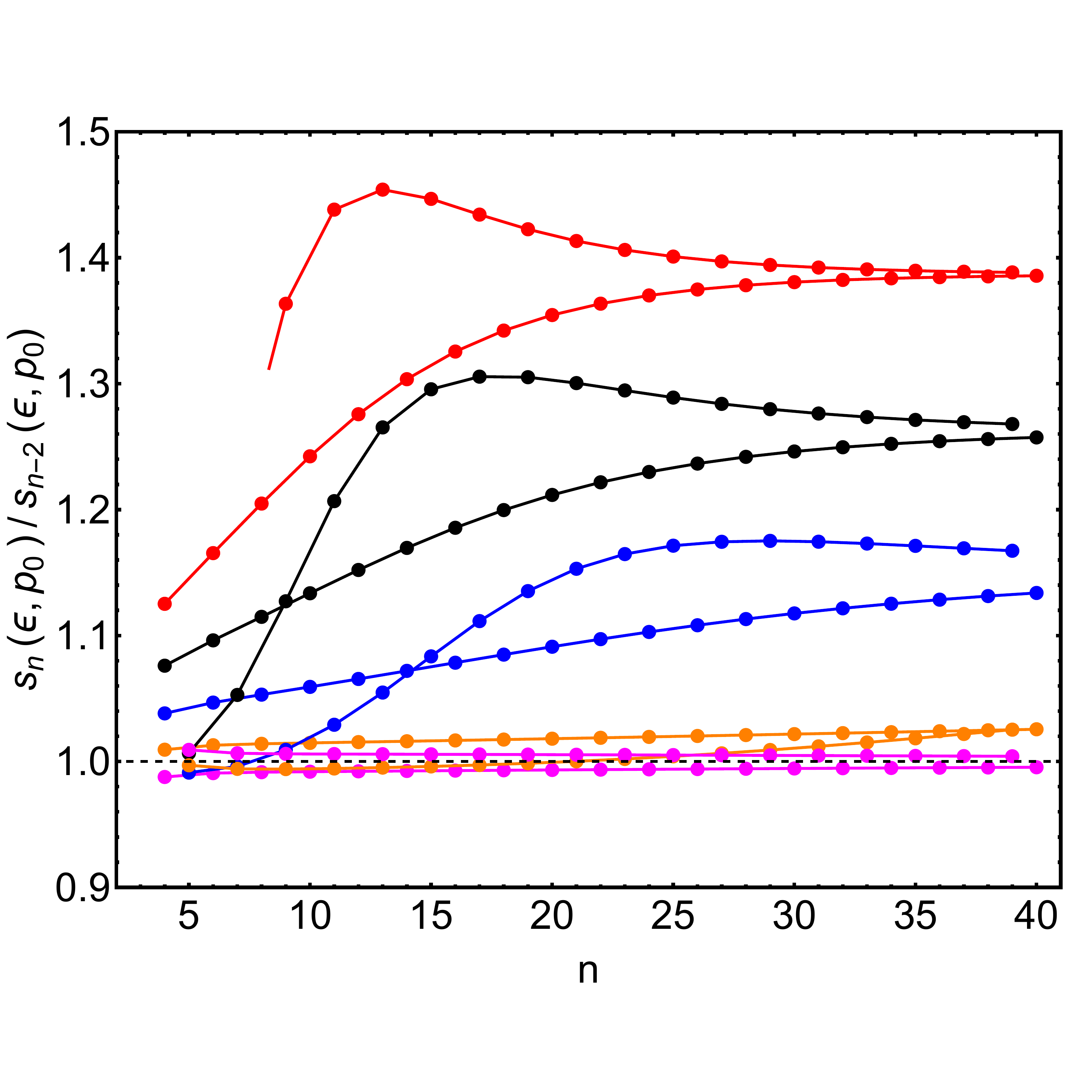}
\includegraphics[height=7cm]{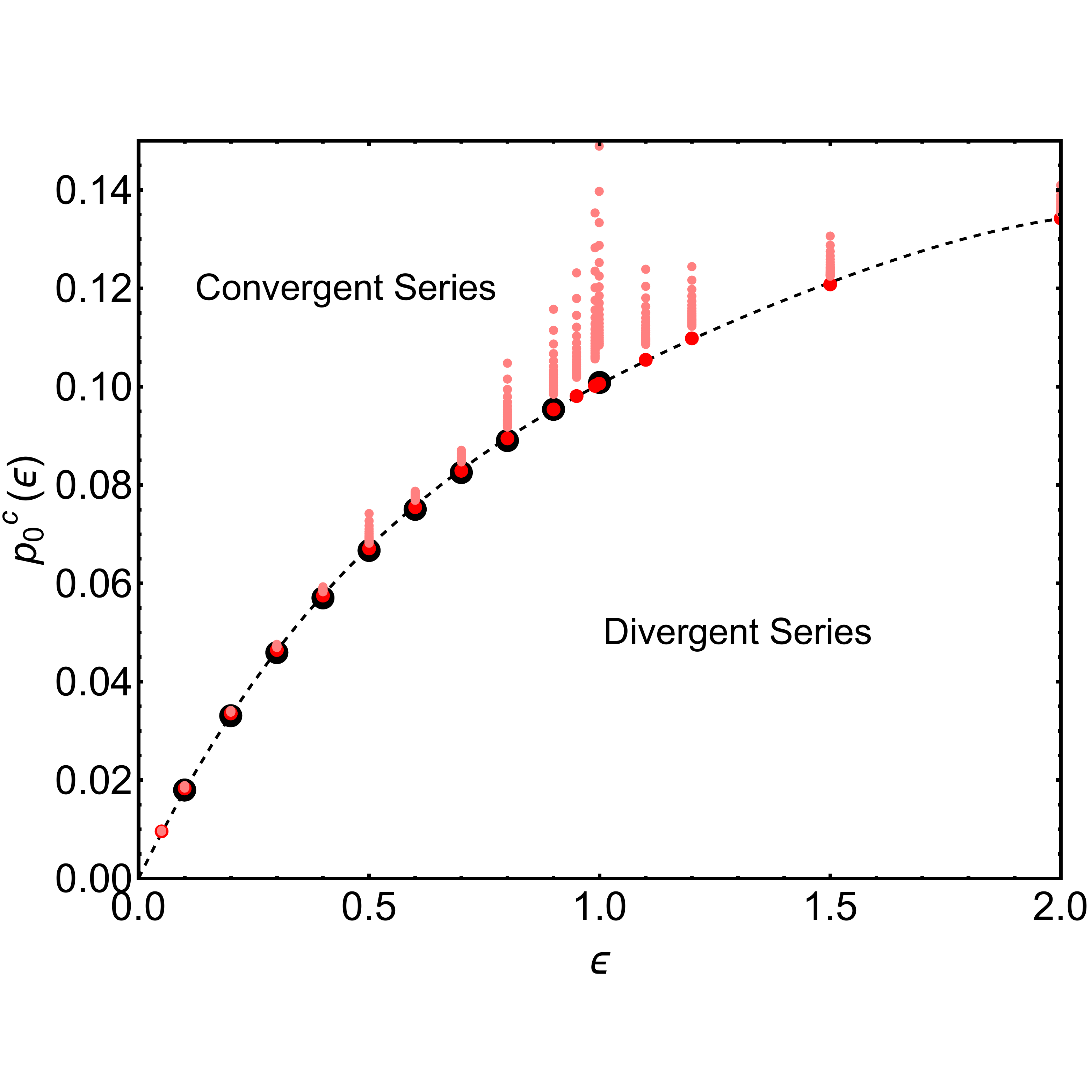}  
\vskip -1cm
\captionof{figure}{Example 1, ODE eq.(\ref{odep1}). Left: Behavior of $s_n/s_{n-2}$ vs $n$ for $\epsilon=0.5$ and $p_0=0.06, 0.062, 0.064, 0.066, 0.068$ (Red, Black, Blue, Orange, Magenta).  $s_n$ is defined by eq.(\ref{quo}). For a given colour, the top curve is for odd values of $n$ and the bottom for the even ones. Right: Black dots are the numerical computation of $p_0^c(\epsilon)$ by using the parameter $s_n/s_{n-2}$. Pink dots are the values of $p_0$ that solve the equation: $y_{2k-1}(x';\epsilon,p_0)=y(x')$ for $x'=0.5$, several $\epsilon$'s and $k\in[3,19]$. Dotted line is a help's eye.}\label{fig8}  
\end{center}
We show in figure \ref{fig8} (left) how the ratio $s_n/s_{n-2}$ is larger than one for $p_0= 0.06$, $0.062$, $0.064$, $0.066$ and is is equal to one (or it converges to it) for $p_0=0.068$. Observe that for each value of $p_0$ there are two curves corresponding to odd (even) values of $n$ (top and bottom respectively). Nevertheless, both asymptotically converge to the same value for large $n$'s. Figure \ref{fig8} (right) shows the numerical critical values (black dots) of $p_0^c(\epsilon)$ obtained by using the ratio $s_n/s_{n-2}$. The method described in this paper have sense whenever the minimum values for the distances $d_{1,2}$ are larger than $p_0^c(\epsilon)$. 

Finally, let us comment a structural behavior of $y_n$.
At some point we realized that the odd $n$-approximations $y_{n}(x;\epsilon,p_0)$ where systematically crossing the exact $y(x)$ at $p_0^z(x,n,\epsilon)$ given $x$ and $\epsilon<1$ . We explored the behavior of such zeros with respect the minimum distances $p_0^*(n,\epsilon)$ to see it was possible there were some relation. We observed that $p_0^z(x,n,\epsilon)$ tends, for increasing values of $n$, to be very near the critical value $p_0^c(\epsilon)$ (see pink dots at figure \ref{fig8} (right)). Moreover, for $\epsilon=1$ only the even approximations have crossing points and we found that $\lim_{n\rightarrow\infty}p_0^z(x,n,\epsilon=1)=\lim_{n\rightarrow\infty}p_0^*(n)$ for any $x$. This is a different way to show that the approximations converge to the ODE's solution. Nevertheless, observe that the  odd $n$ approximations do not cross $y(x)$ for any $p_0$ value and it also converges to the ODE's solution when $n$ tends to infinity. 
\end{itemize}

A general conclusion of this analysis is that the instability of our expansions permits them to explore a larger space of possible solutions efficiently just by minor variations of the parameters. Eventually, they cross the exact solution for a given $x$ and $\epsilon$. Nevertheless, only the $\epsilon=1$ makes possible a uniform approach to the solution for any $x$-value. This high sensibility to the parameter values will be crucial when looking at several solutions in BVP.

\section{Example 2: Bratu differential equation \boldmath $e^{-y}y''+1=0$ (BVP)}

We test our method with a more complex BVP. The ODE is given by
\begin{equation}
e^{-y} y''+1=0\quad,\quad y(0)=\bar y_0\quad ,\quad y(1)=\bar y_1 \label{odep2}
\end{equation}
This well known ODE has the generic solution:
\begin{equation}
y(x)=\ln\left[\frac{2B^2}{\cosh^2\left(Bx+C\right)}\right] \label{odep2e}
\end{equation}
where $B$ and $C$ are constants that are found from the equations:
\begin{eqnarray}
\sinh^2B&=&2B^2\left[4e^{-(\bar y_0+\bar y_1)/2}\sinh^2\frac{B}{2}-\left(e^{-\bar y_1/2 }-e^{-\bar y_0/2}\right)^2\right]\nonumber\\
\cosh C&=&\sqrt{2}e^{-\bar y_0/2}B\label{odep2e2}
\end{eqnarray}
$B$ can be found numerically once we fix the boundary conditions.
\begin{center}
\includegraphics[height=6.5cm]{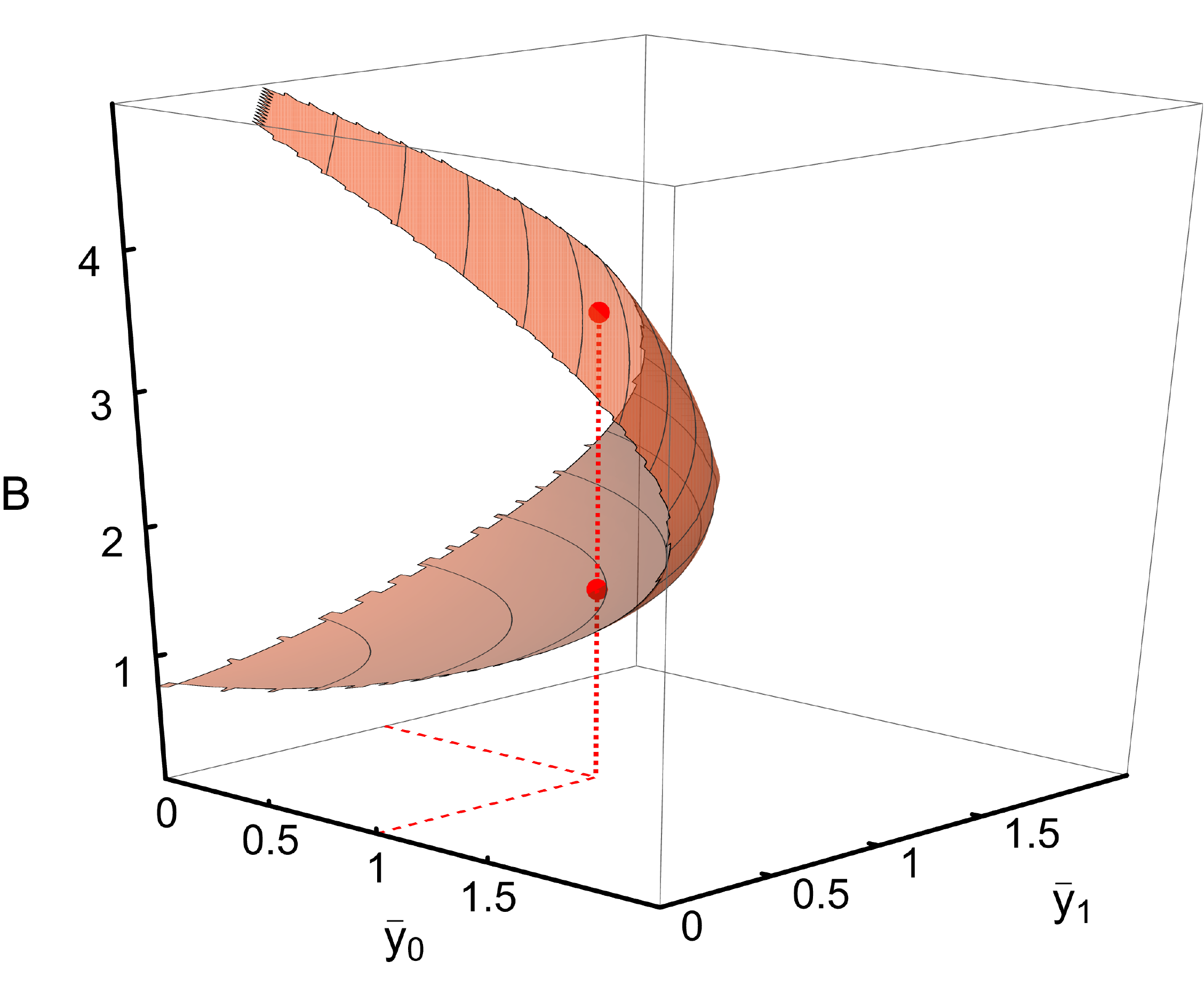}
\includegraphics[height=6.5cm]{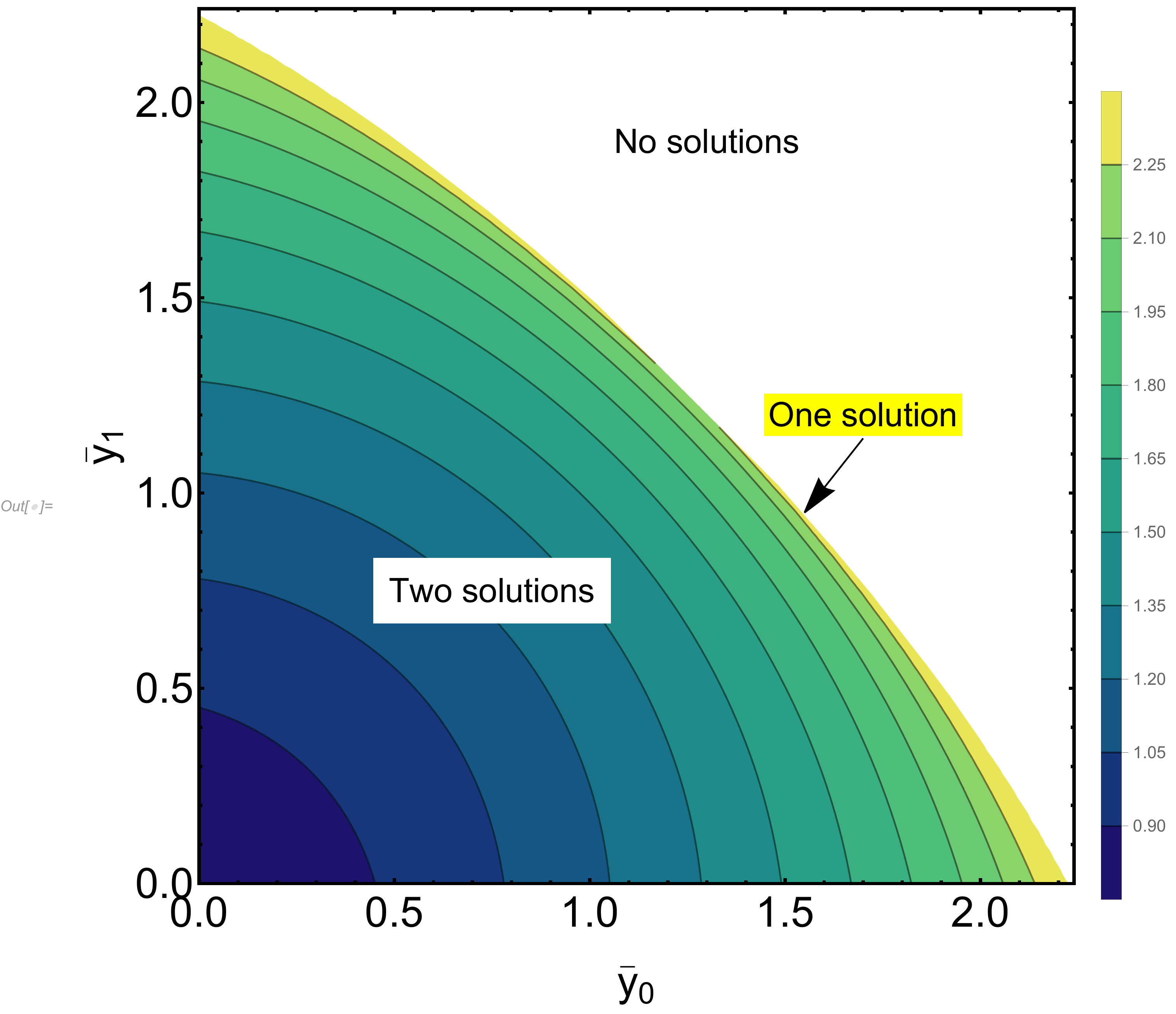}
\captionof{figure}{Numerical solution of equation (\ref{odep2e2}) for $B$-constant as a function of the initial conditions $\bar y_0$ and $\bar y_1$. Left: Dotted lines and red points localize the $B$ values for $\bar y_0=\bar y_1=1$. Right: Contour Plot}\label{fig9}  
\end{center}
We show on figure \ref{fig9} the $B$-values as a function of $\bar y_0$ and $\bar y_1$. We see that there is a region with two solutions, a curve with one solution an the rest with no solutions. In the case $\bar y_0=\bar y_1=\bar y$ the implicit equation for $B$ (\ref{odep2e2}) becomes $\cosh B/2=\sqrt{2}B\exp{(-\bar y/2)}$. Therefore, the limiting value $B^*$ with only one solution is obtained from the equation: $\tanh(B^*/2)=2/B^*$ that implies $B^*\simeq 2.39936..$. We will get two solutions when $\bar y<\bar y^*\equiv 2\log(2\sqrt{2}/\sinh(B^*/2))=1.2567...$, one solution when $\bar y =\bar y^*$ and no solution when $\bar y >\bar y^*$. We focus in this section on the case $\bar y=1$ that has two solutions: $B_1^*=1.51812..$ and $B_2^*=3.5675..$ We show in figure \ref{fig10} the two solutions corresponding to $\bar y_0=\bar y_1=1$.
\begin{center}
\includegraphics[height=7cm]{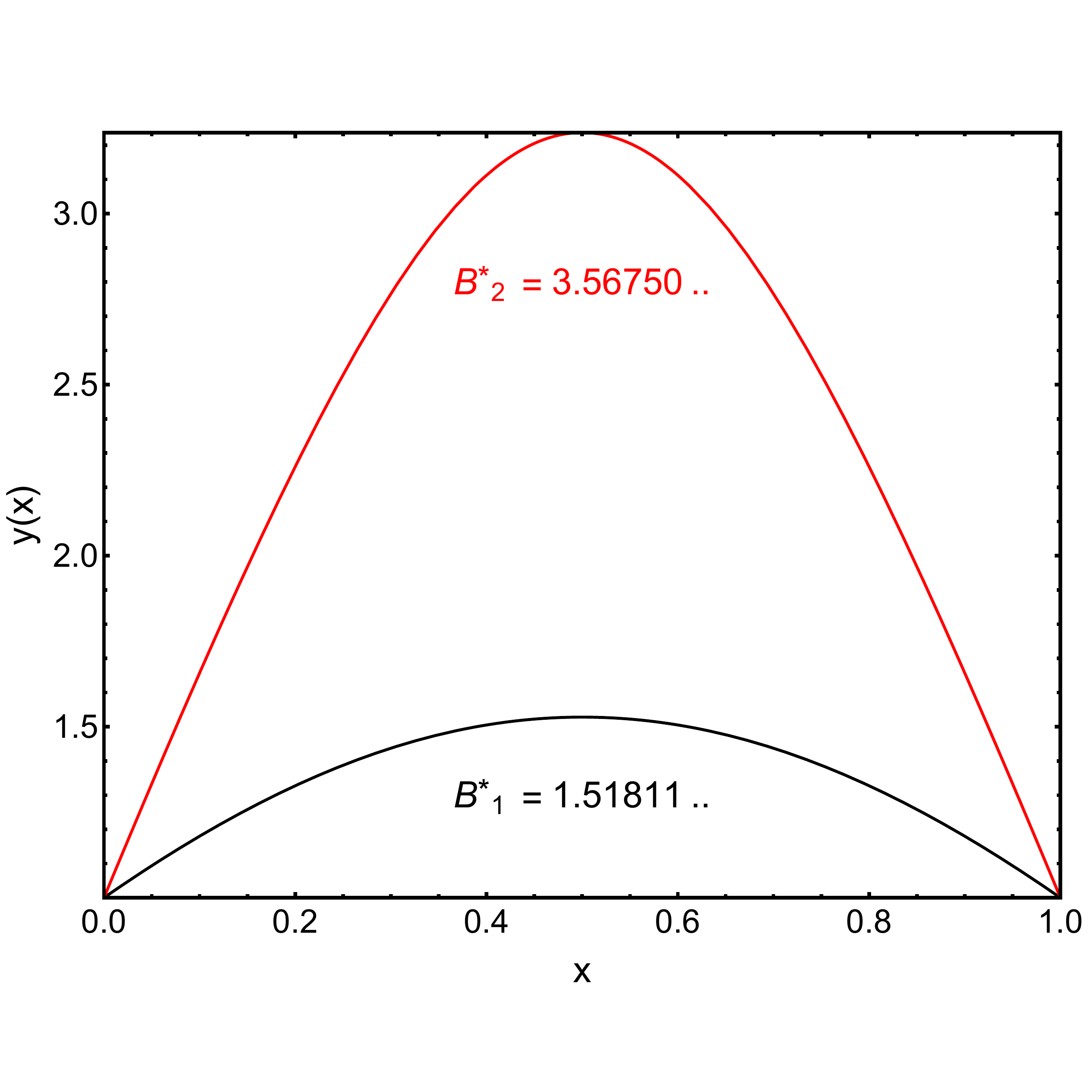}
\vskip-0.75cm
\captionof{figure}{Solutions of ODE's equation (\ref{odep2}) with $\bar y_0=\bar y_1=1$.  $B_{1,2}^*$-constants are the solutions of eq.(\ref{odep2e2}) in this case. The solutions are given by eq.(\ref{odep2e})}\label{fig10}  
\end{center}

\begin{itemize}
\item {\bf eODE's perturbative expansion:}
The ODE (\ref{odep2}) is represented in our notation by $g(x,y,y')=\exp(-y)$ and $h(x,y,y')=1$. 
We restrict our general parametric extended ODE (\ref{eode}) to the case $p_0\neq 0$ and $p_1=p_2=p_3=0$. We will discuss later other cases. The eODE is given by:
\begin{equation}
\left[p_0+\epsilon\left(e^{-y}-p_0\right)\right]y''+\epsilon =0\label{eode22}
\end{equation}
The first terms of the perturbative expansion for $g$  are:
\begin{eqnarray}
g_0(x;p_0)&=&e^{-\tilde y_0(x,p_0)}\nonumber\\
g_1(x;p_0)&=&-e^{-\tilde y_0(x,p_0)}\tilde y_1(x,p_0)\nonumber\\
g_2(x;p_0)&=&e^{-\tilde y_0(x,p_0)}\left[\frac{1}{2}y_1(x,p_0)^2-y_2(x,p_0)\right]\nonumber\\
&&\ldots
\end{eqnarray}
The differential equations to be solved order by order are
\begin{equation}
p_0 \tilde y''_n=F_n(x;p_0)
\end{equation}
where $F_n$'s are generated order by order from eq.(\ref{basicODE}). The first orders are:
\begin{eqnarray}
F_0(x;p_0)&=&0\nonumber\\
F_1(x;p_0)&=&-1\nonumber\\
F_2(x;p_0)&=&-\left(e^{-\tilde y_0(x;p_0)}-p_0\right) \tilde y_1''(x;p_0)
\end{eqnarray}
After solving the ODEs we get $\tilde y_n(x;p_0)$. Their explicit expressions for the first orders are:
\begin{eqnarray}
\tilde y_0(x;p_0)&=&1\nonumber\\
\tilde y_1(x;p_0)&=&\frac{1}{2p_0}x(1 - x)\nonumber\\
\tilde y_2(x;p_0)&=&\frac{(e p_0-1) }{2 e p_0^2}x(1-x) \nonumber\\
\tilde y_3(x;p_0)&=&\frac{1}{24 e^2 p_0^3}\left[12 e^2 p_0^2+e \left(-24 p_0-x^2+x+1\right)+12\right]x(1-x)
\label{expan_2}
\end{eqnarray}
We have computed with Mathematica up to $\tilde y_{40}$.
\begin{center}
\vskip -1cm
\includegraphics[height=7cm]{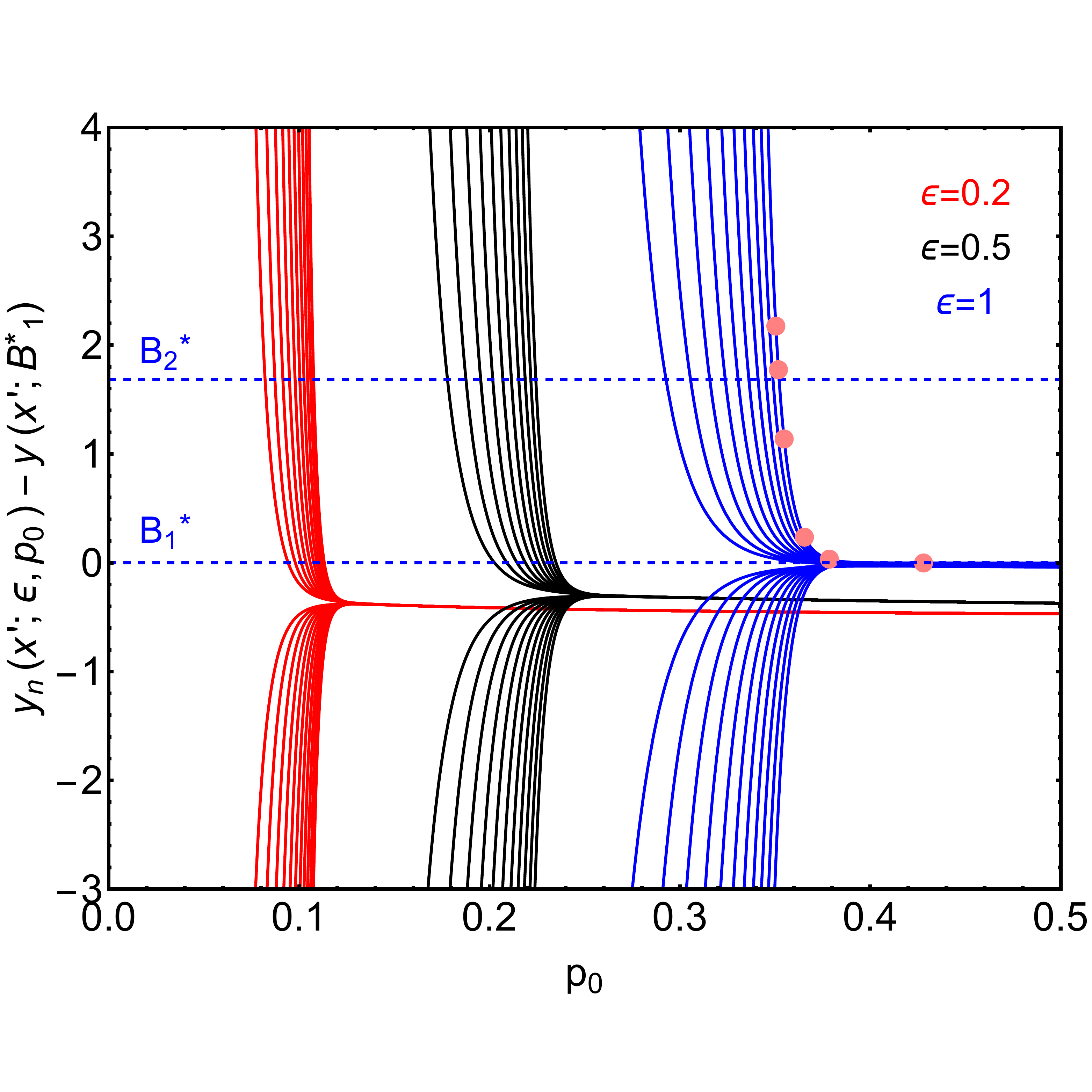}
\includegraphics[height=7cm]{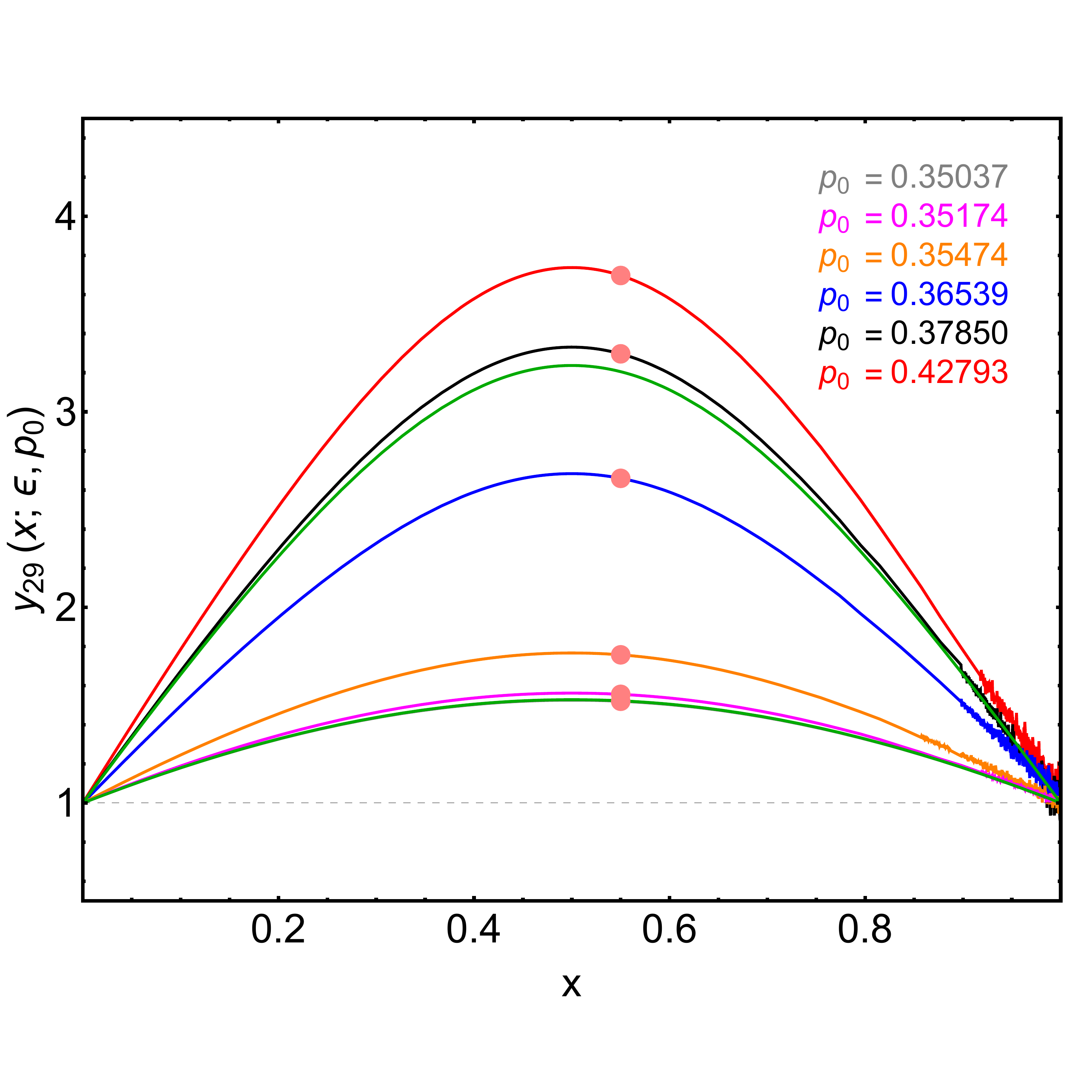}  
\vskip -1cm
\captionof{figure}{Example 2, ODE eq.(\ref{odep2}): Perturbative analysis of the BVP for the eODE given by eq. (\ref{eode22}). Left: Difference between $n$'th approximate solution, $y_n$, and the exact one, $y(x;B_1^*)$ given by eq. (\ref{odep2e}) at the point $x'=0.55$ as a function of $p_0$. There are shown the cases $\epsilon=0.2$, $0.5$ and $1$ (red, black and blue lines respectively) for $n=10,\ldots, 30$. The horizontal dotted line with labels $B_{1,2}^*$ shows the values for the known exact solutions $y(x'=0.55';B_{1,2}^*)$. Pink dots represent particular values of $p_0$ (from left to right): $0.35037$, $0.35174$, $0.35474$, $0.36539$, $0.37850$ and $0.42793$ for $n=29$. Right: $y_{29}(x;p_0)$ and $p_0$ taking values shown in the figure. Pink dots are the values for $x'=0.55$ and correspond to the ones on the left figure. Green curves are the two exact solutions $y(x;B_{1,2}^*)$.}\label{fig11}  
\end{center}
We see in figure \ref{fig11} (left) the behavior of $y_n(x'=0.55;p_0)-y(x'=0.55;B_1^*)$ for different $n$-values. Odd (even) $n$-values up (down) when $p_0$ decreases. As we discussed in Example 1 above, we see how $y_n$, $n$ even or odd, for the case $\epsilon=1$ follows one of the exact solutions, $y(x;B_1^*)$, for sufficiently large values of $p_0$ (for example in figure \ref{fig11} (right) we see that the case $p_0=0.35037$ is indistinguishable from the exact solution to the naked eye). Observe that we target the second solution $y(x;B_2^*)$ thanks to the unstable behavior of $y_n$ (n-odd) when $p_0$ goes to zero.  In some sense, such blowing up behavior with $p_0$ permits the $y_n$-expansion to explore the space of all possible solutions. However, we will see that we will pay the price for this: the convergence to the second solution will be non-trivial. Therefore, from our perturbative scheme point of view, we can say that $y(x;B_1^*)$ is a {\it regular solution} because the even or odd expansion of $y_n$ will converge smoothly to the exact result and $y(x;B_2^*)$ is a {\it singular solution} because only the odd expansion of $y_n$ will converge in a non-trivial way to it.
\begin{center}
\includegraphics[height=7cm]{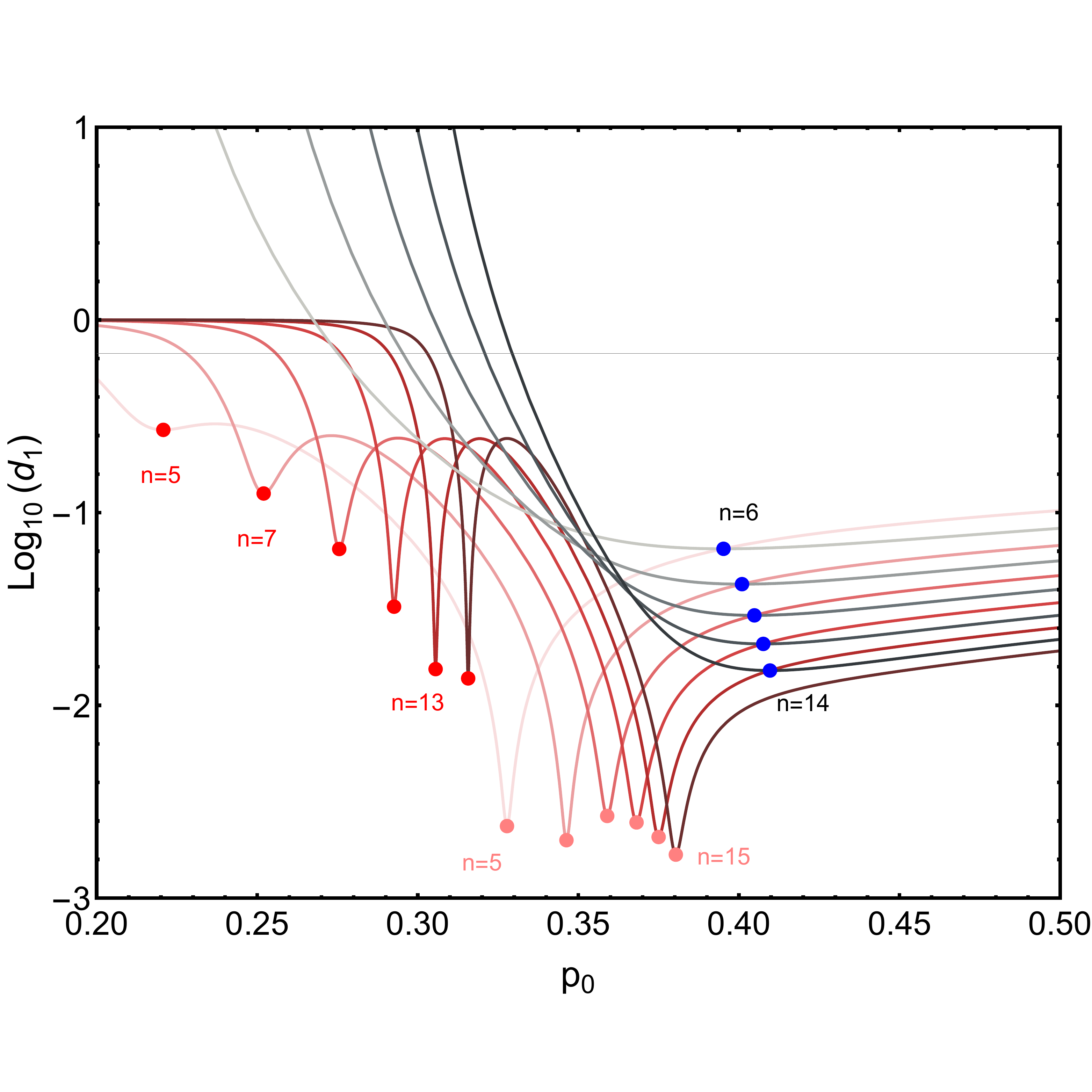}
\includegraphics[height=7cm]{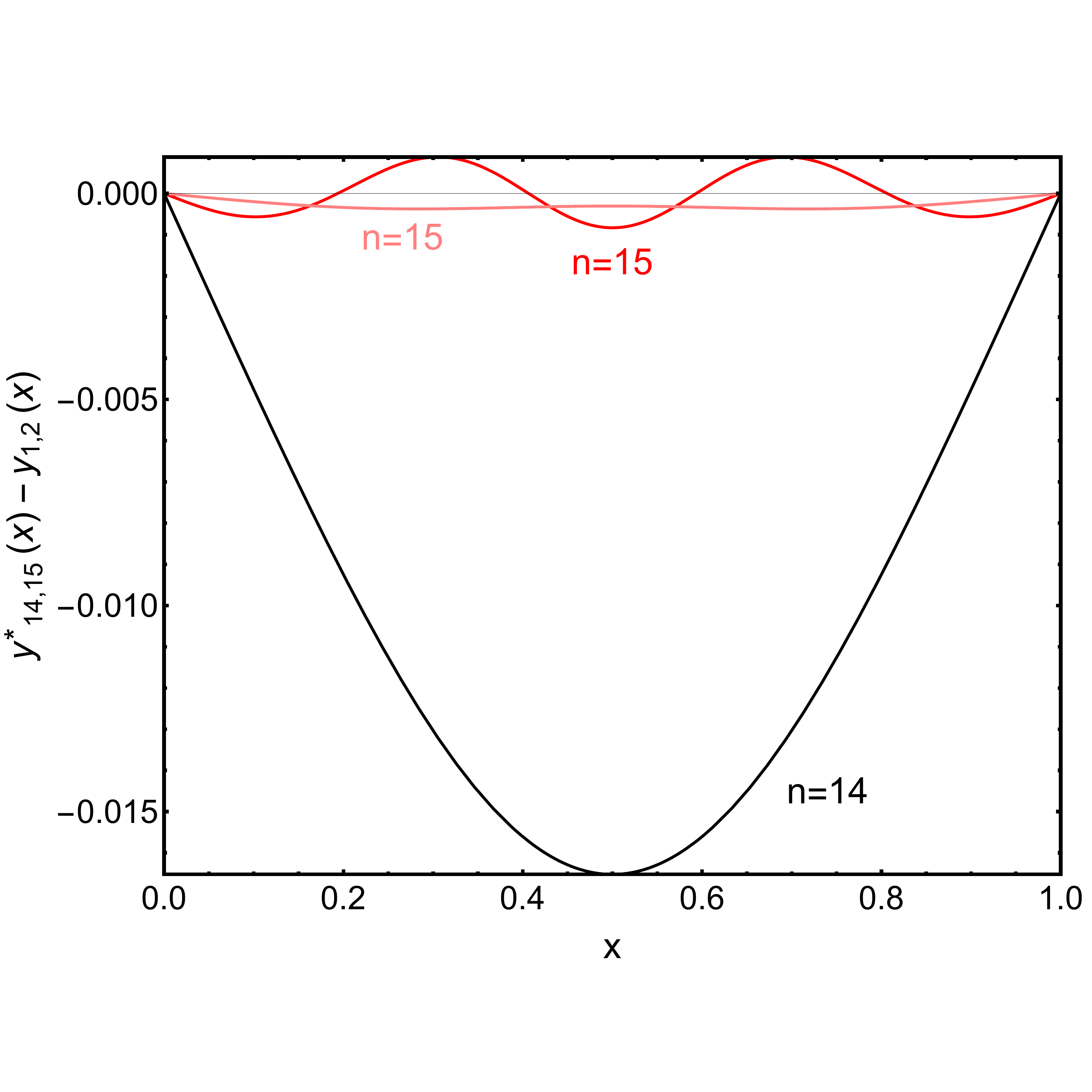}
\vskip-0.75cm
\captionof{figure}{Example 2, ODE eq.(\ref{odep2}). Left: Decimal logarithm of distance $d_{1}(n;\epsilon=1,p_0)$ defined by eq.(\ref{mea1}) versus $p_0$ for each $n$th-perturbation approximation. The $y_n$ are obtained from the eODE's perturbative expansion (\ref{eode22}). Gray-Black curves are for increasing even values of $n$ and Red-Cherry curves are for increasing odd values of $n$. We plot only $n=5,\ldots ,15$. From top to bottom $n=1, 2,\ldots,10$. Blue dots shows the unique minima  when $n$ is even and Red-Pink dots shows the existence of two  minima when $n$ is odd. Right: Differences between $y_n(x)$ and the exact known solutions $y(x;B_{1,2}^*)$. $n=14$, $p_0^*(14)$ and $B_1^*$ (Black curve), $n=15$, $p_0^{(2)*}(15)$ and $B_1^*$ (Pink curve) and $n=15$, $p_0^{(1)*}(15)$ with $B_2^*$ (Red curve). $p_0^*(n)$ are the minima of $d_1(n;p_0)$ (see the left figure)}\label{fig12}  
\end{center}

\item {\bf \boldmath The distance $d_{1}(n;\epsilon=1,p_0)$ (Restricted Conjecture):} We use $\tilde y_n(x;p_0)$ to compute $d_1(n;\epsilon=1,p_0)$ defined by eq.(\ref{mea1}) that it is shown in figure \ref{fig12}-left for $n\in[5,15]$. We observe that there is a different behavior with respect the parity of $n$. When $n$ is even $d_1(n)$ have an unique minimum that moves smoothly to larger values as we increase $n$. However, for odd $n$ values, the distance $d_1(n)$ presents a doble minimum structure. Therefore, we have three natural sets of sequences $y_n^*(x)\equiv y_n(x,p_0^*(n))$: when $n=2k$ with an unique minimum $p_0^*(2k)$ and two more when $n=2k+1$ (odd) associated with the other two local minima: $p_0^{(1)*}(2k+1)$ and $p_0^{(2)*}(2k+1)$. For instance in figure \ref{fig12}-right we show the differences: $y_{14}(x;p_0^*(14))-y(x;B_1^*)$ (black curve), $y_{15}(x;p_0^{(2)*}(15))-y(x;B_1^*)$ (pink curve) and $y_{15}(x;p_0^{(1)*}(15))-y(x;B_2^*)$ (red curve). The differences are, at most, of order $10^{-2}$ at each $x$-value.
\begin{center}
\vskip -0.5cm
\includegraphics[height=7cm]{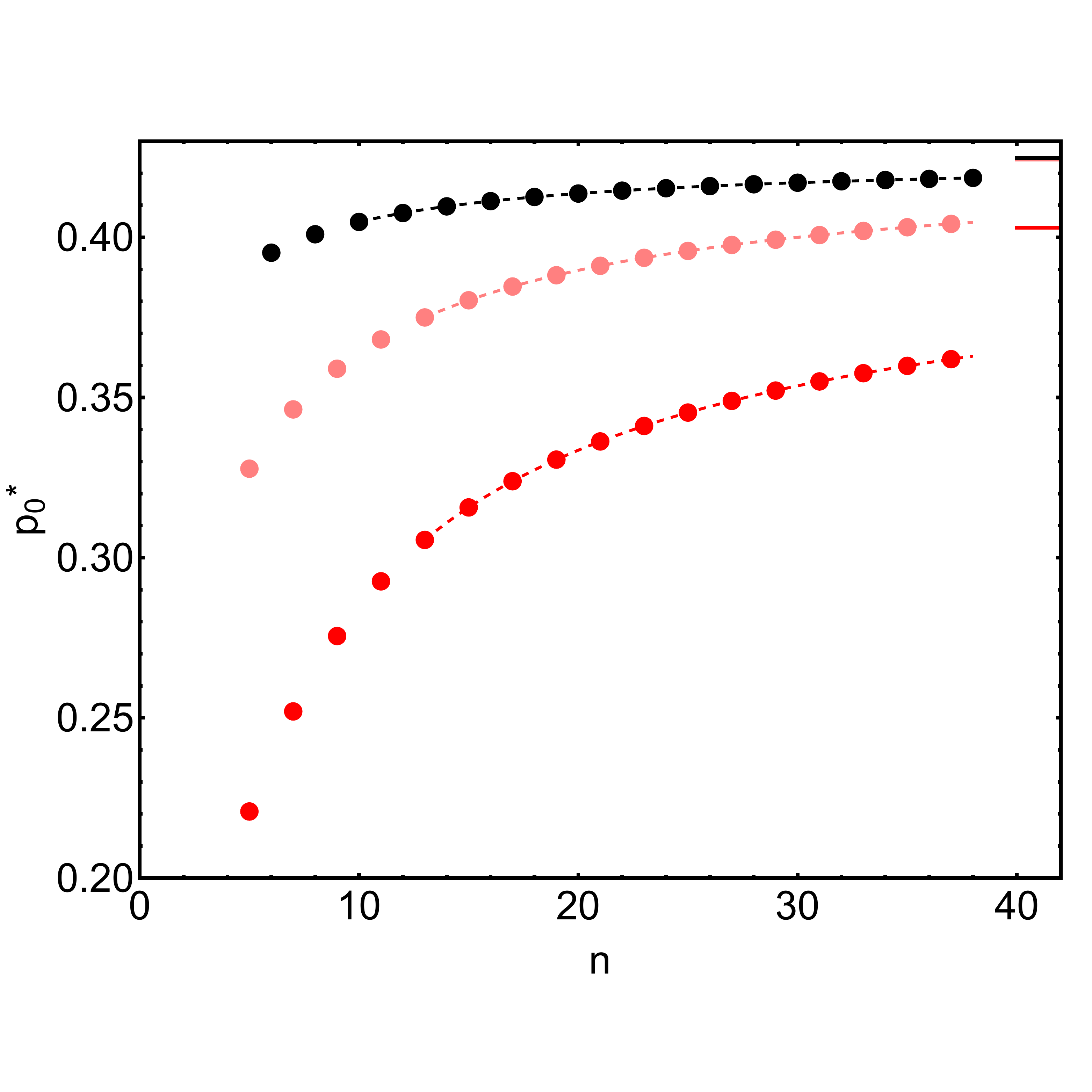}  
\includegraphics[height=7cm]{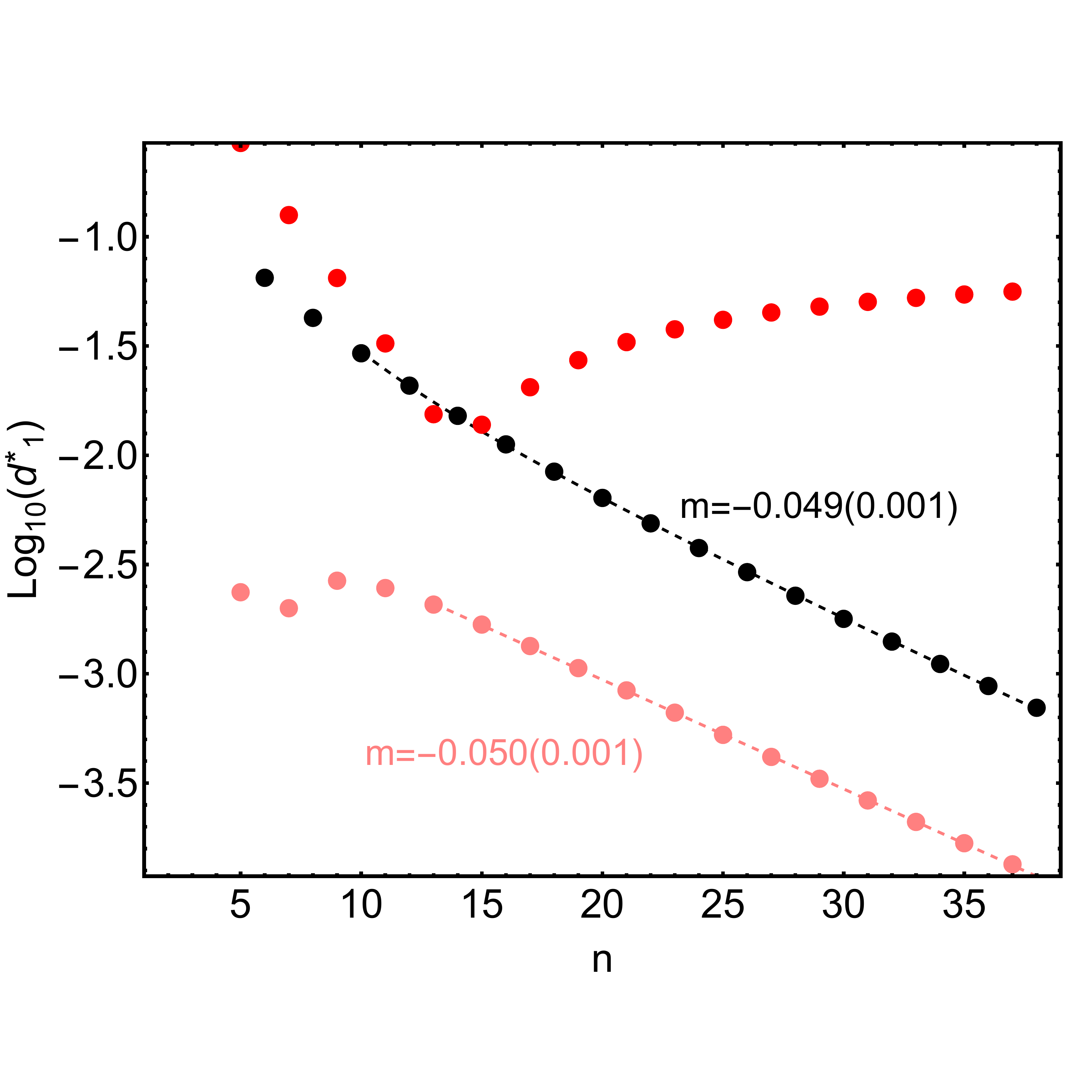}
\vskip -0.5cm
\captionof{figure}{Example 2, ODE eq.(\ref{odep2}): Asymptotic behavior with the $n$-th perturbative approximation of the values where the minima of $d_{1}$ are located, $p_0^*$ (left figure) and the values of such distances at the minima, $d_{1}^*$. Black dots are the unique minimum for each $n$-even and pink-red dots are the two local minima corresponding to $n$-odd. Dashed lines are fits explained in the main text.} \label{fig13}  
\end{center}

We show in figure \ref{fig13} how the value of such minima, $p_0^*$, and their correspondig distance, $d_1^*$ change with the $n$'th approximation order. The position of the local minima, $p_0^*$, increase with $n$ smoothly in all cases and we are able to fit the function: $p_0^*(n)=a_0+a_1/n+a_2/n^2$ (dotted lines in the figure). The assymptotic values $a_0$ are: $0.42465(0.00002)$ ($n$ even, black dots in figure \ref{fig13}), $0.42424(0.00006)$ ($n$ odd, the second minimum $p_0^{(2)*}$, pink dots) and $0.4031(0.0003)$ ($n$-odd, the first minimum $p_0^{(1)*}$, red dots). We see that the black and pink dots converge to the same asymptotic value and therefore, they are approaching to the same solution (in this case $y(x;B_1^*)$). The red dots converge to a different value (associated to the other solution $y(x;B_2^*)$). The distance $d_1^*$ decrease exponentialy fast with $n$ for $n$-even and $p_0^{(2)*}$, and for $n$-odd (see figure \ref{fig13}-right). In fact we fit $\log_{10}d_1^*(n)=m n+s$ with $m=-0.049(0.001)$ for $n$-even and $m=-0.050(0.001)$ when $p_0^{(2)*}$ and $n$-odd. These values show that the exponential decay is slow. In fact, the distance diminish by an order of magnitude each $20$ iteration steps. 

Let us mention the unexpected behavior of the sequence approaching the solution $y(x;B_2^*)$. While the associated minimum, $p_0^{(1)*}$ increases monotonously with $n$ to a limiting value (red dots in figure \ref{fig13} -left), its corresponding distance diminishes exponentially fast up to $n=13$ and then increases a bit and tends to a finite value (see red dots in figure \ref{fig13}-right). We will discuss later why we think that this is happening, how to deal when this effect occurs and, in any case, how to increase in a systematic way the precision of an approximated solution. Nevertheless, we are convinced that this different behavior is connected with the ``singular'' character of the solution (as we called it).

In figure \ref{fig14} we show the relation between $d_1^*$ to the distance to the exact solutions, $d_{ex}^*$ (\ref{meae}) in this test problem where we know them explicitly. 
We see that the sequence $y_{2k}^*$ (black dots) has a distance $d_1^*(2k)\simeq  d_{ex}^*(2k)^{a}$ with $a=1.006(0.003)$ for large values of $k$ (straight dotted line in the figure). The distance associated to the odd-$n$ values and $p_0^{(2)*}$ have a similar behavior: $d_1^{(2)*}(2k+1)\simeq  d_{ex}^{(2)*}(2k+1)^{a}$ with $a=1.13(0.04)$ for large values of $k$. We can conclude again that both $d_1^*$ distances are good measures of the distance to the exact result and it makes our scheme a self-consistent method to approach the exact solutions.
Finally,  the behavior of the distance $d_1^{(1)*}(2k+1)$ as a function of $d_{ex}$ (red points in figure \ref{fig14} reflects the singular behavior when approaching to the $y(x;B_2^*)$ solution.
\begin{center}
\vskip -1.0cm
\includegraphics[height=8cm]{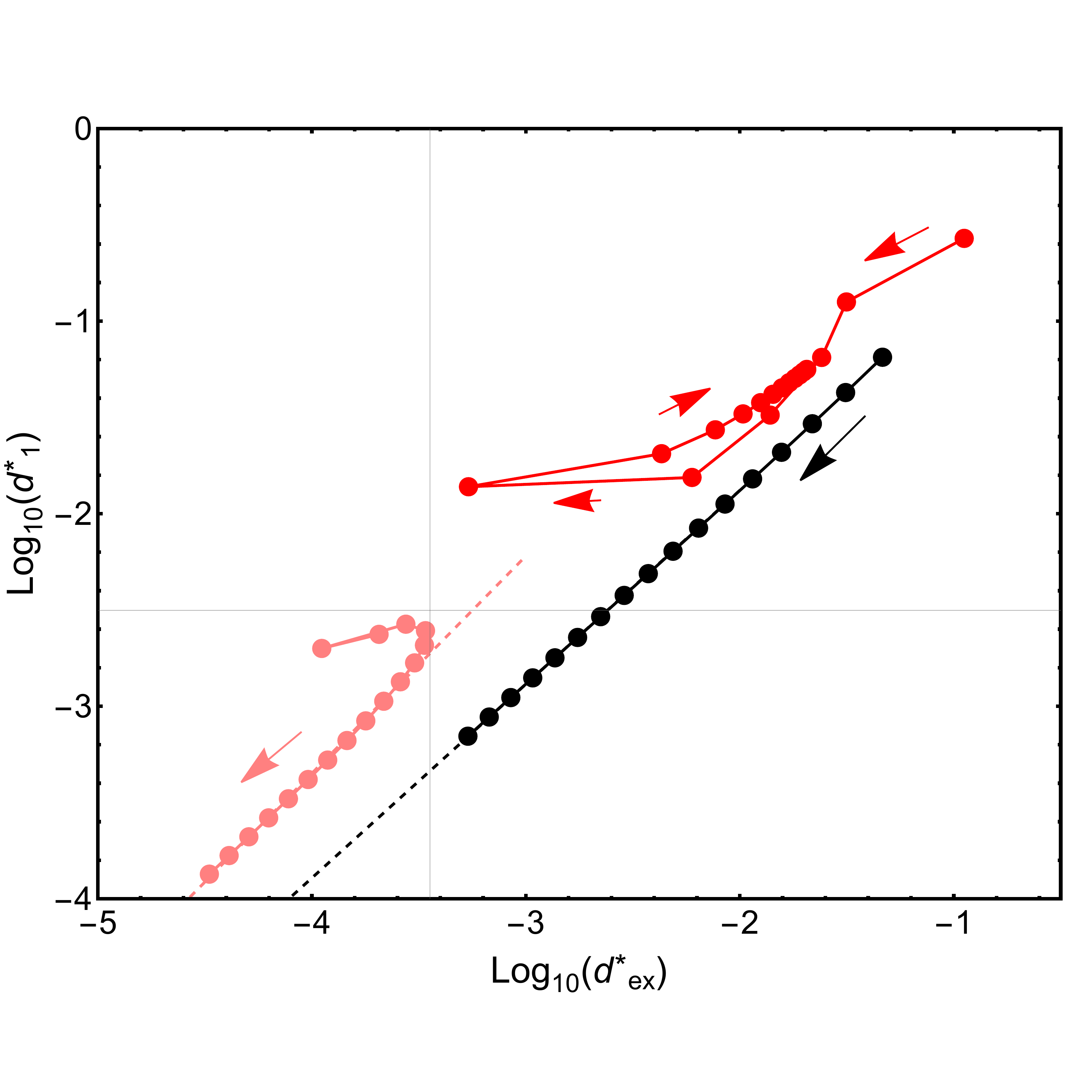}
\vskip-0.75cm
\captionof{figure}{Example 2, ODE eq.(\ref{odep2}): Decimal logarithm of the optimal distance $d_{1}^*$ defined by eq.(\ref{mea1}) versus the one for the distance to the exact solution, $d_{ex}^*$, for each $n$th-perturbation approximation. The arrows indicate increasing values of $n$ from $n=5$ up to $n=40$. Black dots correspond to even $n$-values and pink and red dots for odd $n$-values and for the two minima of $d_1$ found in this case. Dashed lines are linear fits that are explained in the main text.}\label{fig14}  
\end{center}
\begin{center}
\vskip -1.0cm
\includegraphics[height=7cm]{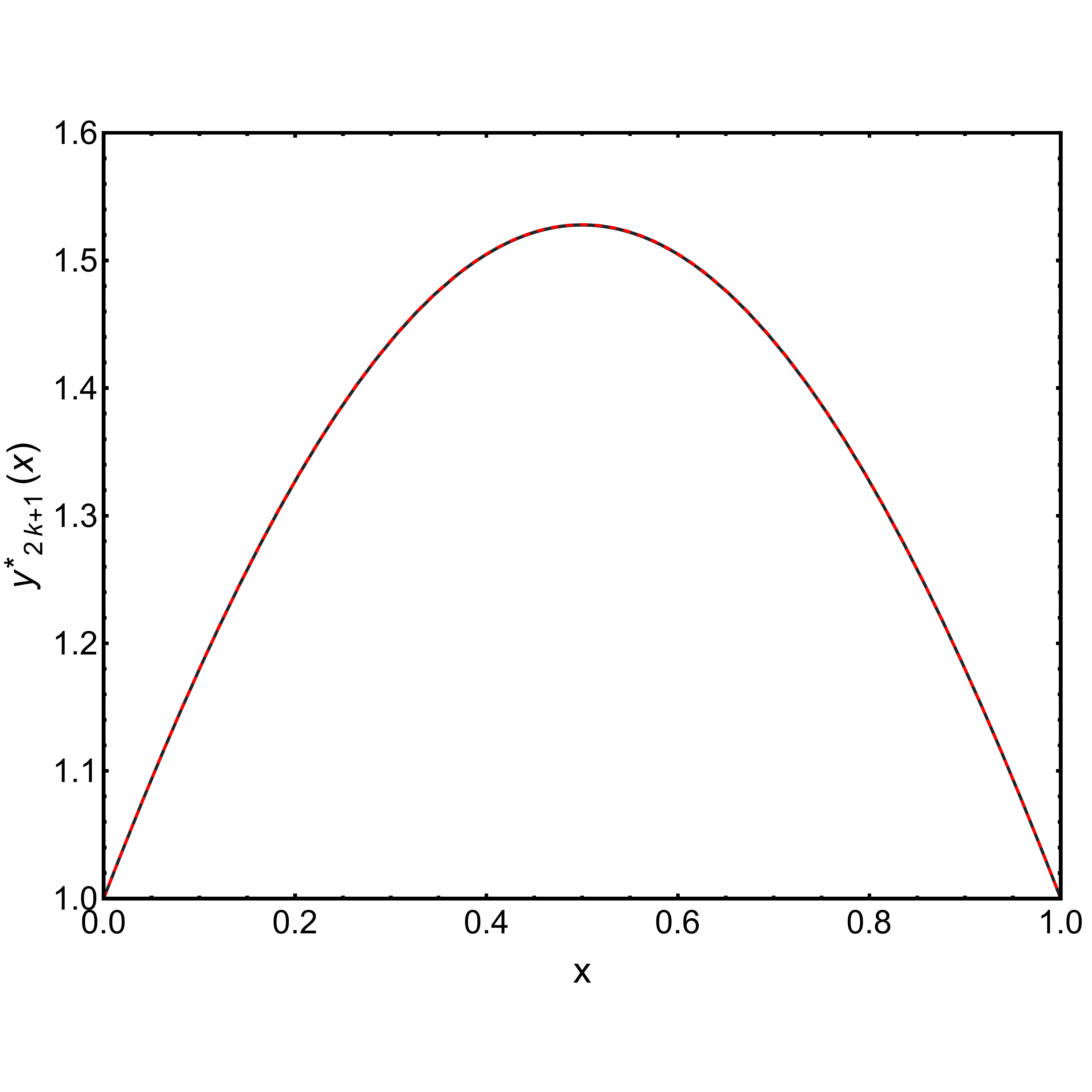}
\includegraphics[height=7cm]{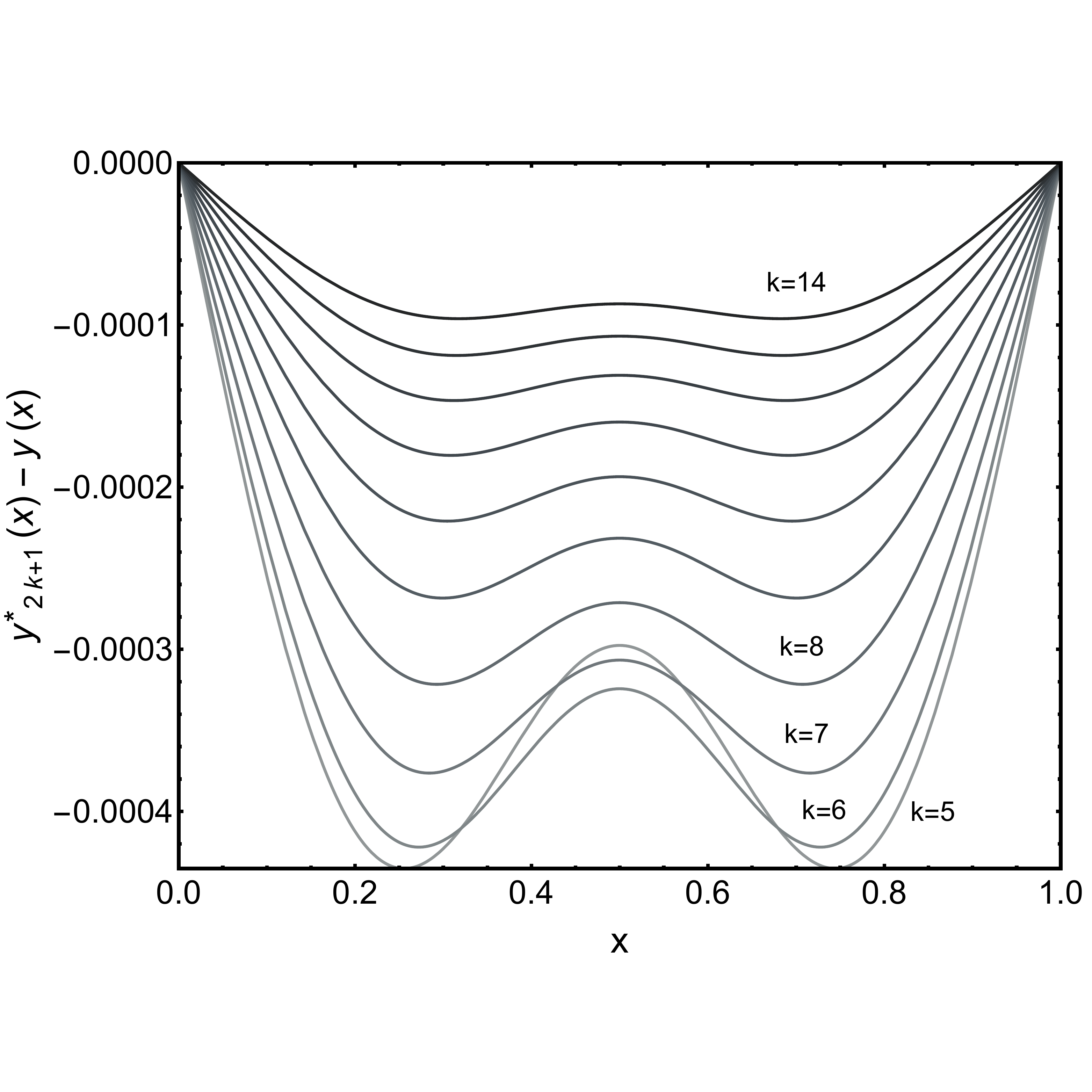}
\includegraphics[height=7cm]{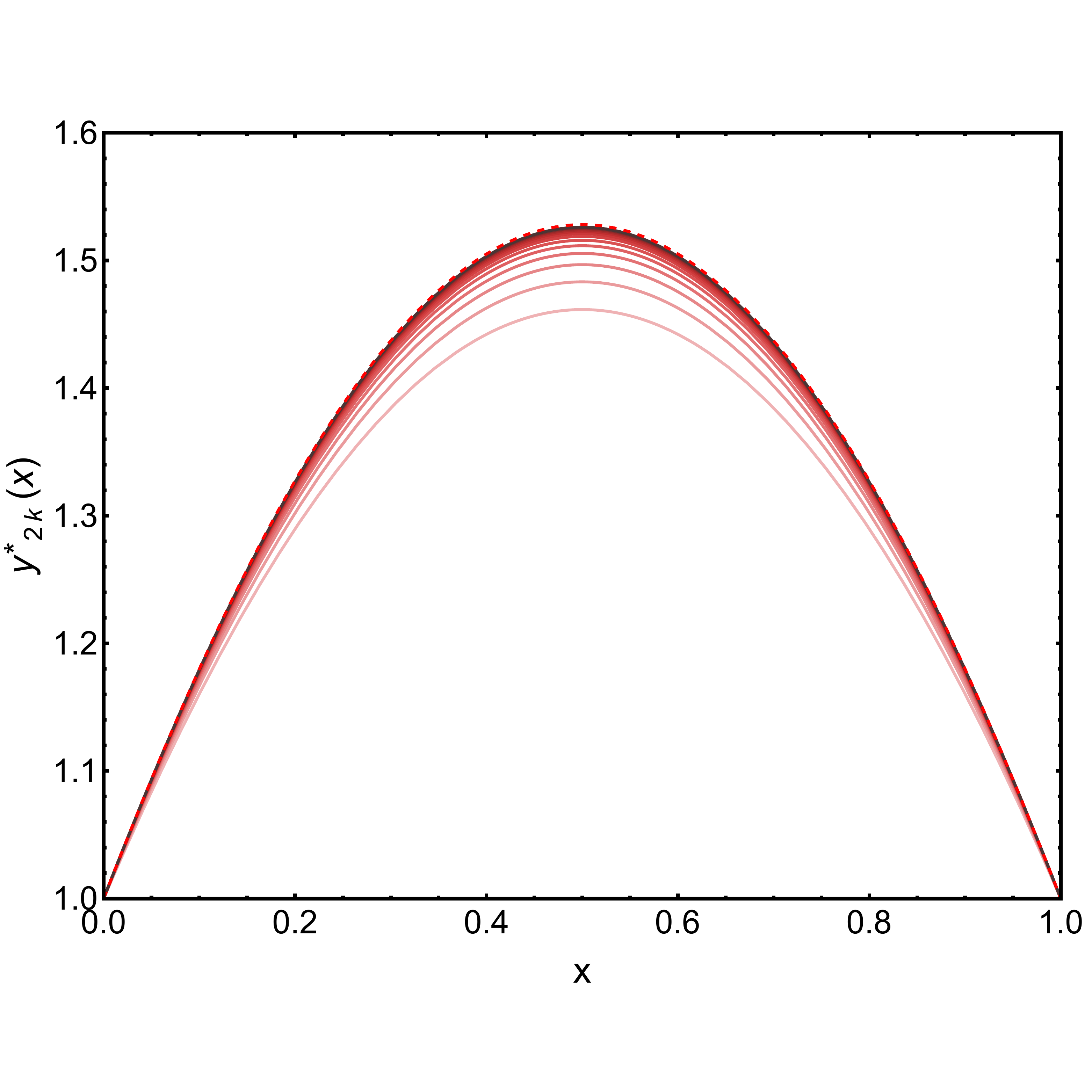}
\includegraphics[height=7cm]{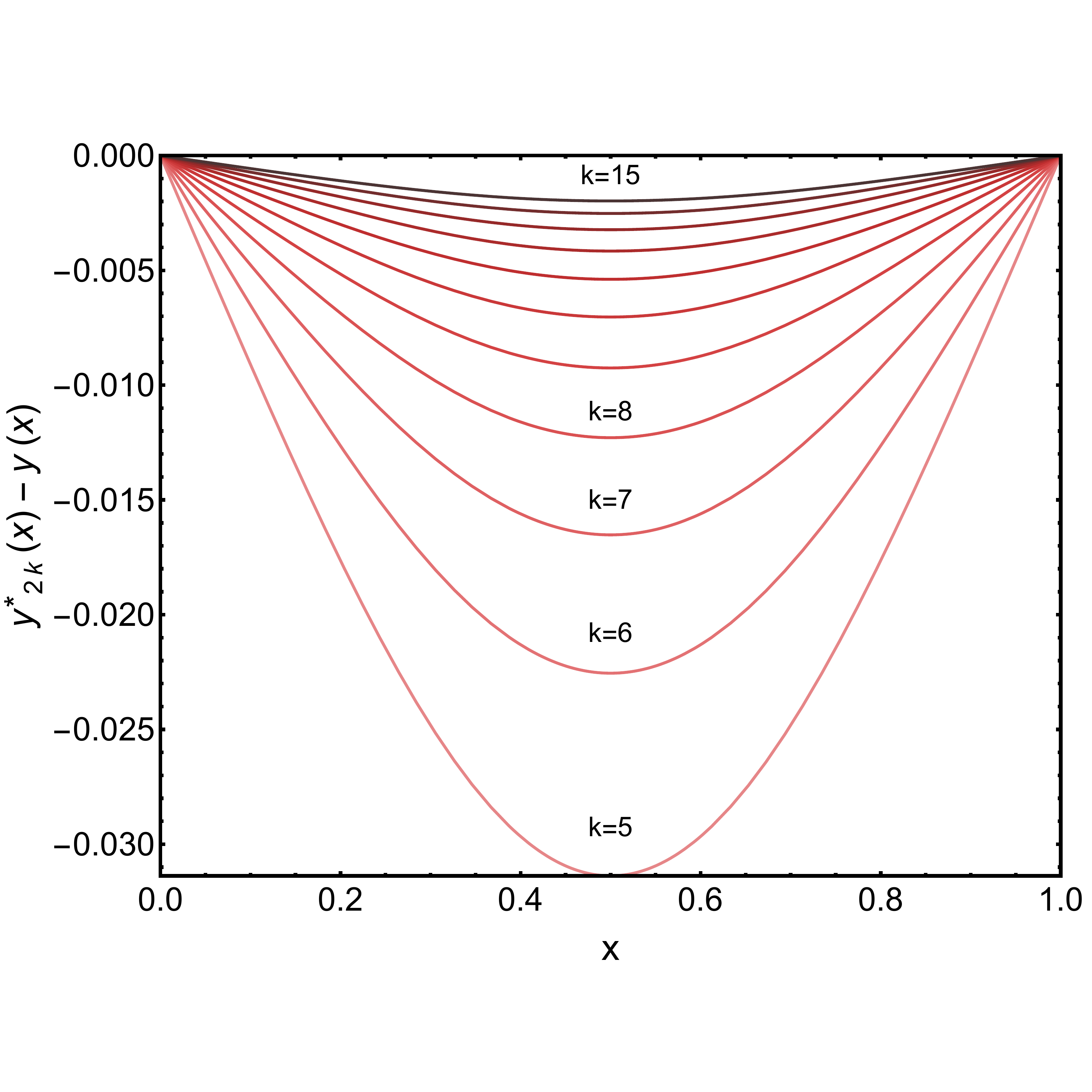}
\vskip-0.75cm
\captionof{figure}{Example 2, ODE eq.(\ref{odep2}): Aproximations to the solution $y(x;B_1^*)$. Top row left: $y_{2k+1}(x;p_0^{(2)*}(2k+1))$ with $k=2,..,7$. Top row right: $y_{2k+1}(x;p_0^{(2)*}(2k+1))-y(x;B_1^*)$ with $k=2,..,14$. Bottom row left: $y_{2k}(x;p_0^{*}(2k))$ with $k=3,..,15$. Bottom row right: $y_{2k}(x;p_0^{*}(2k))-y(x;B_1^*)$ with $k=5,..,15$. 
Increasing color curve intensities indicates larger values of $k$.
}\label{fig15}  
\end{center}
\begin{center}
\vskip -1.0cm
\includegraphics[height=7cm]{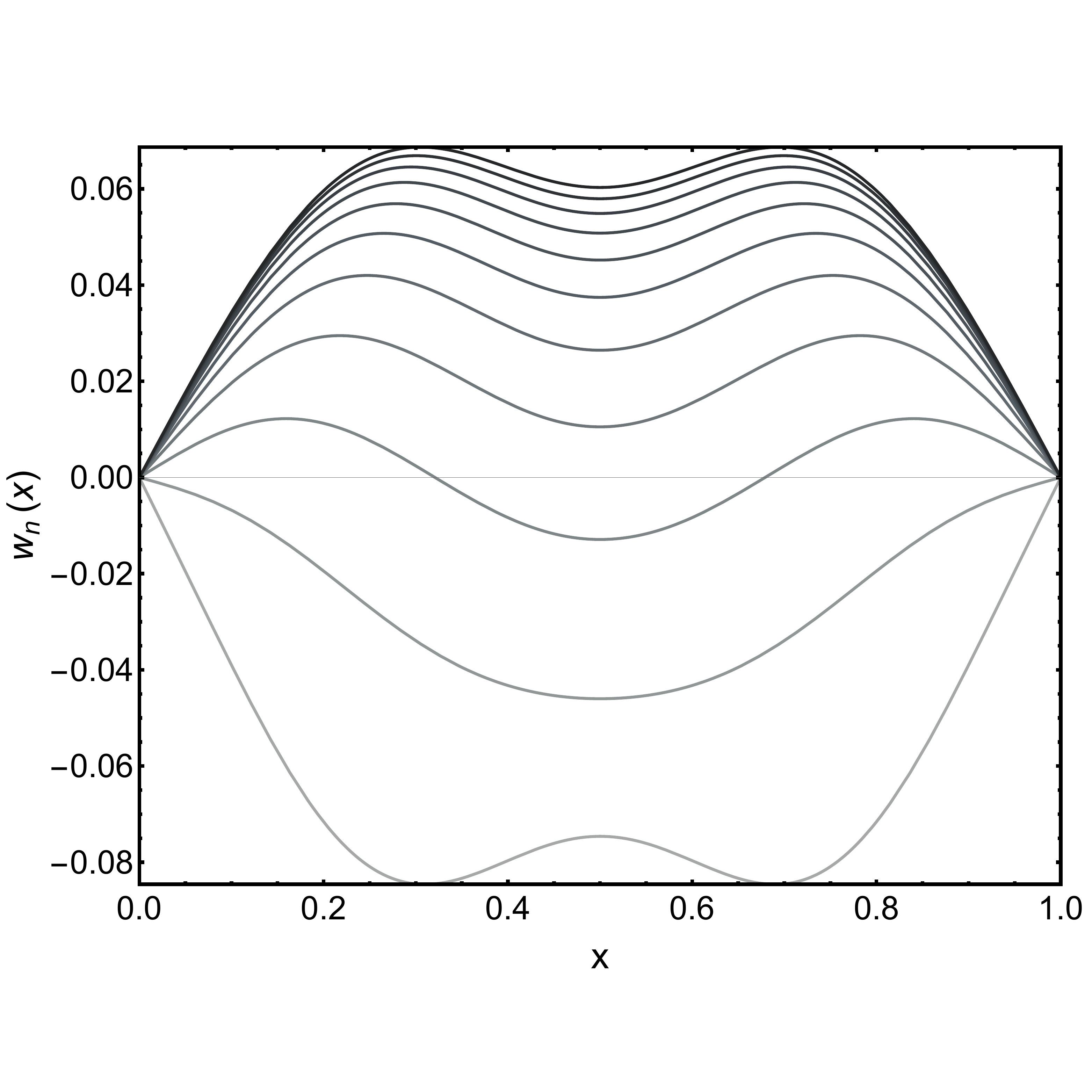}
\includegraphics[height=7cm]{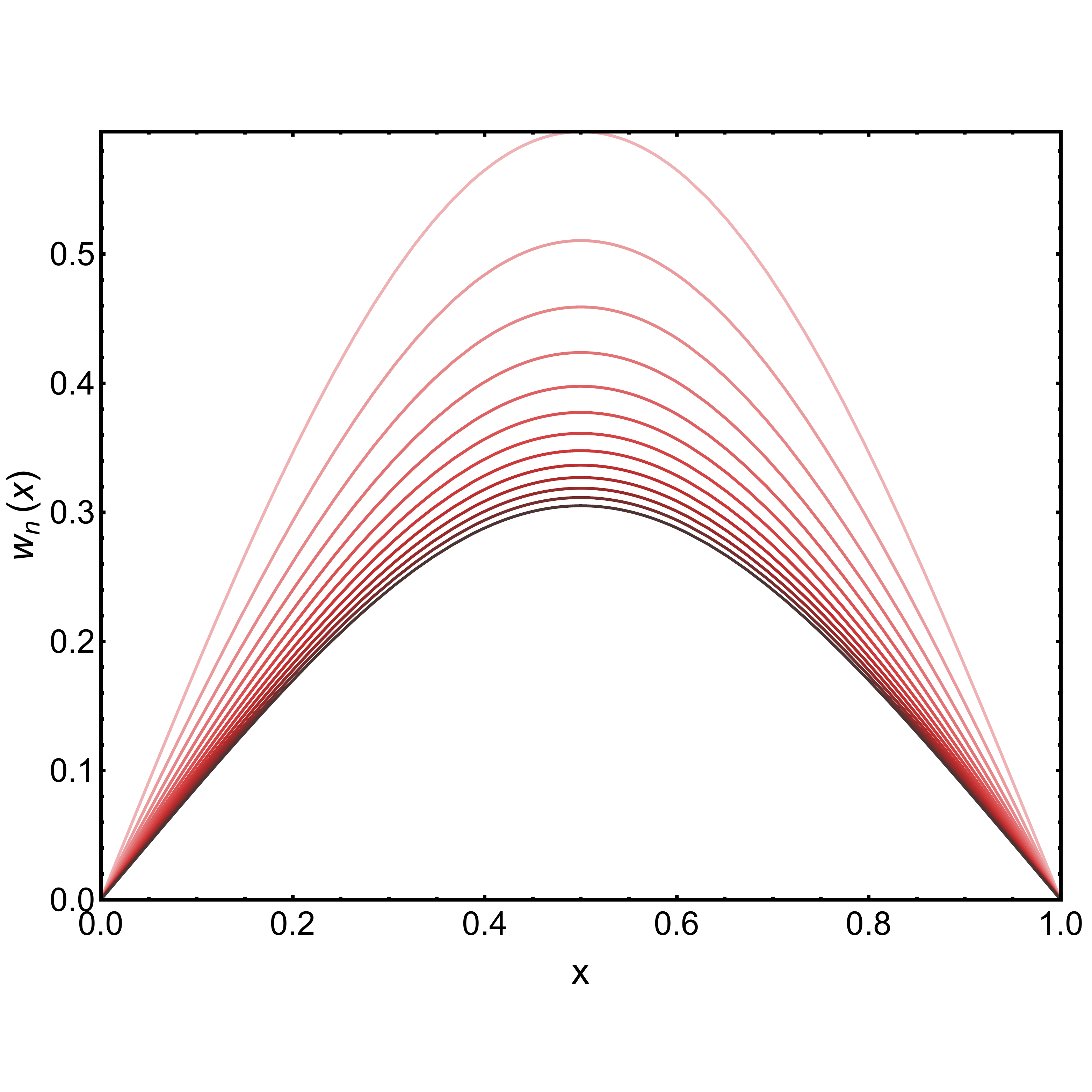}
\vskip-0.75cm
\captionof{figure}{Example 2, ODE eq.(\ref{odep2}): Ghost expansion terms $w_n(x;j)$ for the solution $y(x;B_1^*)$. $j=1$ is associated to the sequence $y_{2k+1}(x;p_0^{(2)*}(2k+1))$ (left figure) and $j=2$  for the sequence $y_{2k}(x;p_0^{*}(2k))$ (right figures).  
Increasing color curve intensities indicates larger values of $n$. $n=4,\ldots,14$ (left) and $n=3,\ldots,15$ (right).
}\label{fig16}  
\end{center}
\item {\bf The approximations to the solution {\boldmath $y(x;B_1^*)$}:} We see on figure \ref{fig15} how the sequences 
$y_{2k+1}(x;p_0^{(2)*}(2k+1))$ and $y_{2k}(x;p_0^{*}(2k))$ tend to the solution $y(x;B_1^*)$ monotonously. The even sequence needs a larger $k$ value than the odd one to get the same precision level. On both cases we can get explicitly the corresponding {\it Ghost Expansion} (\ref{ghost1}):
\begin{equation}
y(x;B_j^*)=\sum_{m=0}^\infty w_m(x;j)d_1^*(m;j)\quad j=1,2
\end{equation} 
whose first terms are explicitly given by:
\begin{itemize}
\item {\bf  j=1 (\boldmath $y_{2k+1}(x;p_0^{(2)*}(2k+1))$):}
\begin{eqnarray}
w_0(x;1)&=&\left[1+\frac{1}{2p}x(1-x)\right]\frac{1}{d_1^*(0;1)}\nonumber\\
\text{with}&& d_1^*(0;1)=0.134419\quad p=0.290468\nonumber\\
w_1(x;1)&=&-\frac{x(1-x)}{24 e^2 p^3 q d_1^*(1;1)}\left[12 e^2 p^3-36 e^2 p^2 q+36 e p q+e q x^2-e q x-e q-12 q\right]\nonumber\\
\text{with}&& d_1^*(1;1)=0.0224301\quad p=0.301550\quad q=0.290468 \nonumber
\end{eqnarray}
\item {\bf j=2 (\boldmath$y_{2k}(x;p_0^{*}(2k))$):}
\begin{eqnarray}
w_0(x;2)&=&\left[1+\frac{2ep-1}{2ep^2}x(1-x)\right]\frac{1}{d_1^*(0;2)}\nonumber\\
\text{with}&& d_1^*(0;2)=0.215464\quad p=0.367879\nonumber\\
w_1(x;2)&=&-\frac{x(1-x)}{480 e^3 p^4 q^2 d_1^*(1;2)}\biggl[480 e^3 p^4 q-240 e^2 p^4-960 e^3 p^3 q^2+1440 e^2 p^2 q^2\nonumber\\
&+&80 e^2 p q^2 x^2-80 e^2 p q^2 x-80 e^2 p q^2-960 e p q^2+2 e^2 q^2 x^4-4 e^2 q^2 x^3\nonumber\\
&+&e^2 q^2 x^2-60 e q^2 x^2+e^2 q^2 x+60 e q^2 x+e^2 q^2+60 e q^2+240 q^2\biggr]\nonumber\\
\text{with}&& d_1^*(1;2)=0.107717\quad p=0.385692\quad q=0.367879 \nonumber
\end{eqnarray}
\end{itemize}
We show in figure \ref{fig16} the behavior of $w_n(x;j)$ for increasing $n$-values. We observe that both sequences have different functional forms. In both cases, $w_n(x;j)$ converge to a limiting curve when $n\rightarrow\infty$.
\begin{center}
\vskip -1.0cm
\includegraphics[height=7cm]{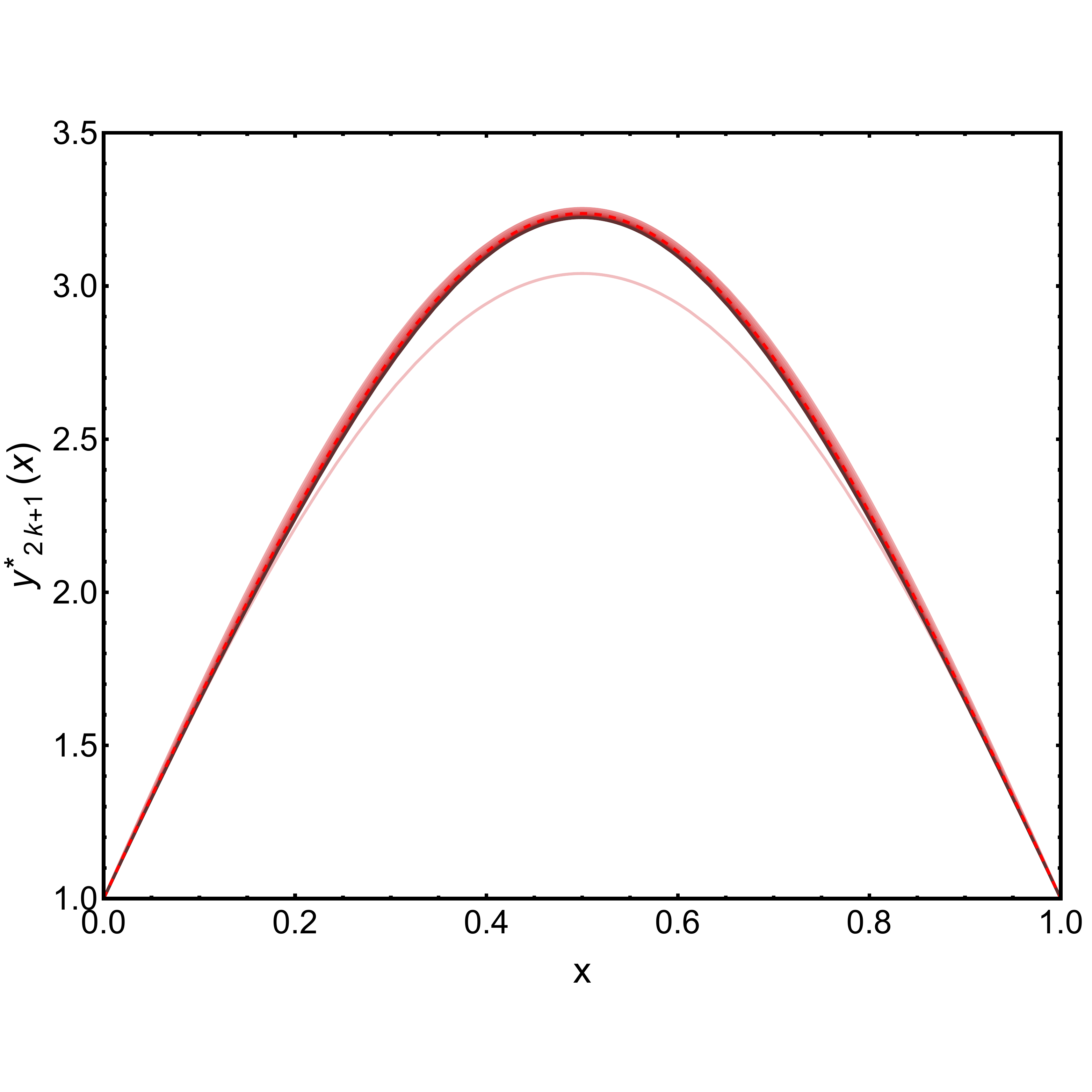}
\includegraphics[height=7cm]{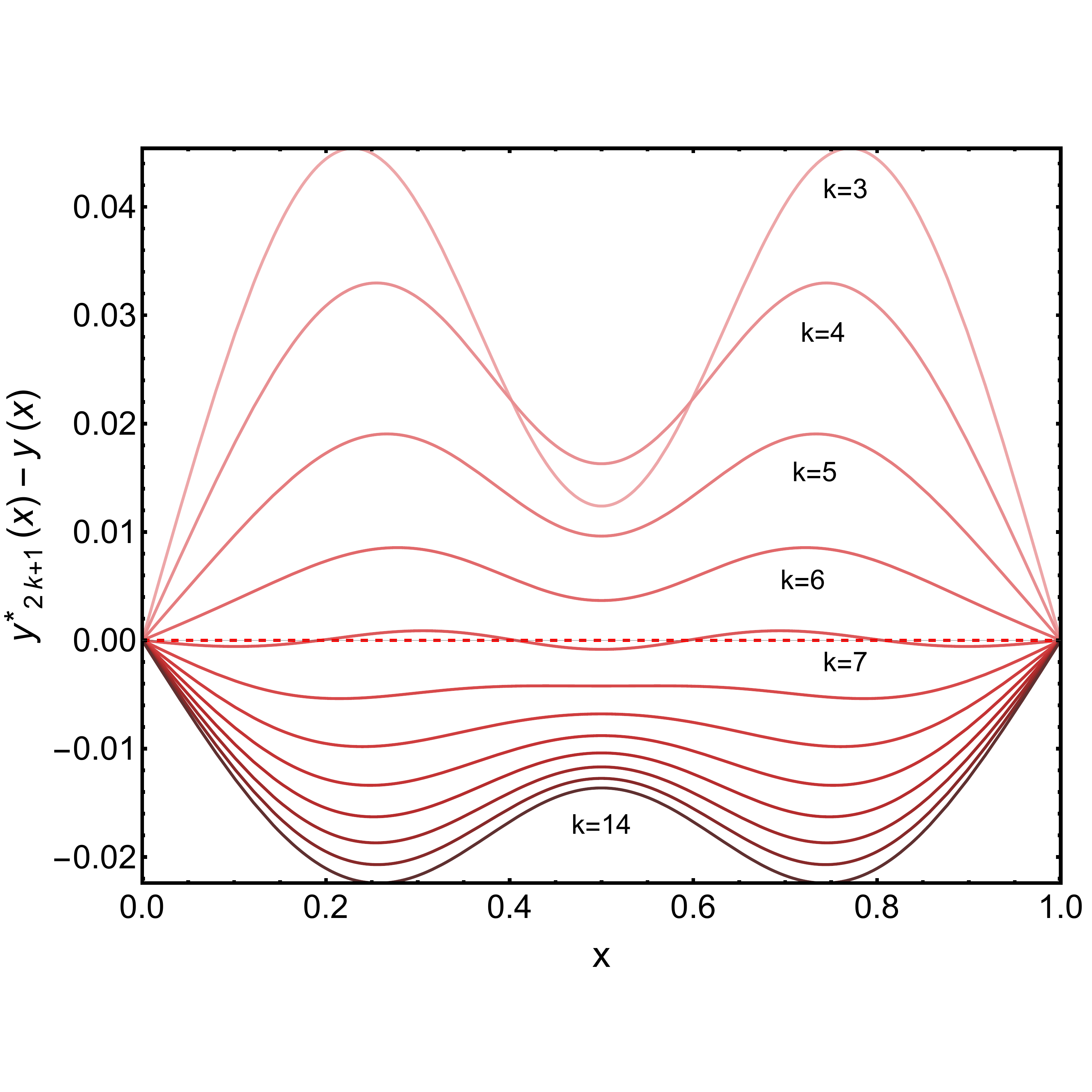}
\vskip-0.75cm
\captionof{figure}{Example 2, ODE eq.(\ref{odep2}): Aproximation to the solution $y(x;B_2^*)$. Left: $y_{2k+1}(x;p_0^{(1)*}(2k+1))$ with $k=2,..,7$. Right: $y_{2k+1}(x;p_0^{(1)*}(2k+1))-y(x;B_2^*)$ with $k=2,..,14$.  
Increasing color curve intensities indicates larger values of $k$.
}\label{fig17}  
\end{center}
\item {\bf The approximations to the solution {\boldmath$y(x;B_2^*)$}:} 
As we already commented above, the no-full convergence of the sequence $y_{2k+1}(x;p_0^{(1)*}(2k+1))$ to $y(x;B_2^*)$ introduce some new issues to this method. Figure \ref{fig17} explicitly shows such convergence where we observe that up to $k=7$, the convergence seem to be a normal one but, afterwards, the optimal configurations separate from the solution and they converge to a limiting form that is near the solution, but it is not the solution.
\begin{center}
\includegraphics[width=5cm]{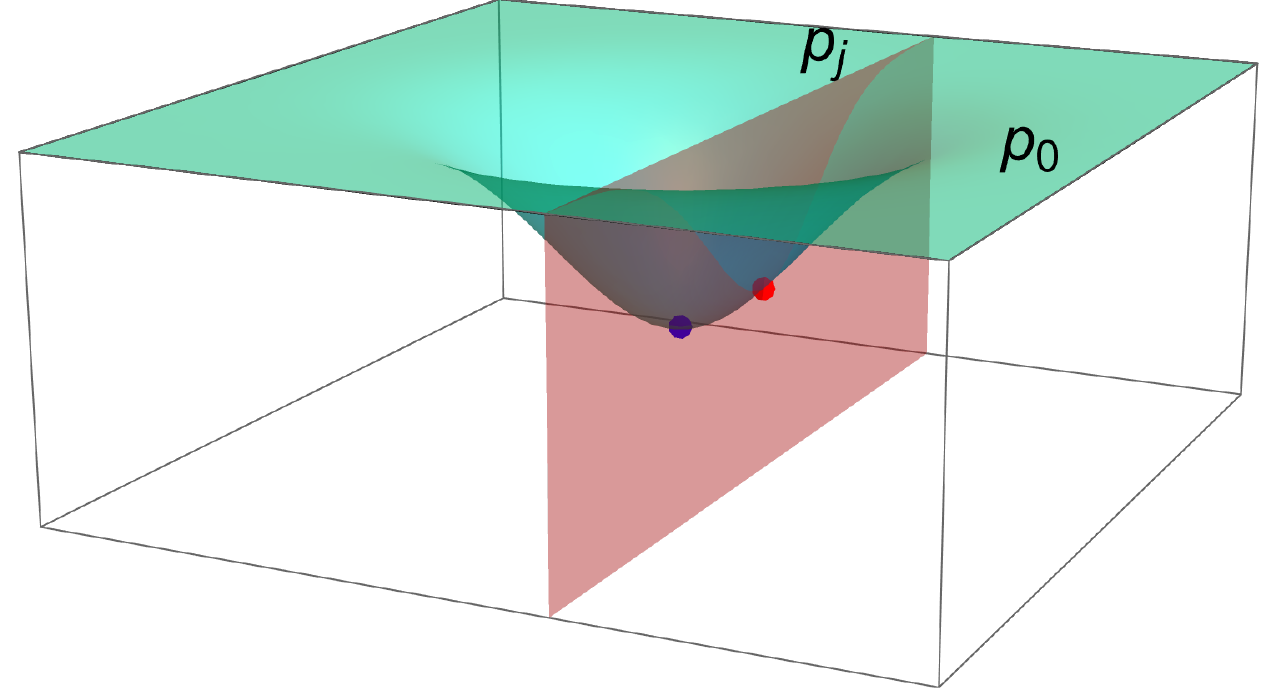}
\includegraphics[width=5cm]{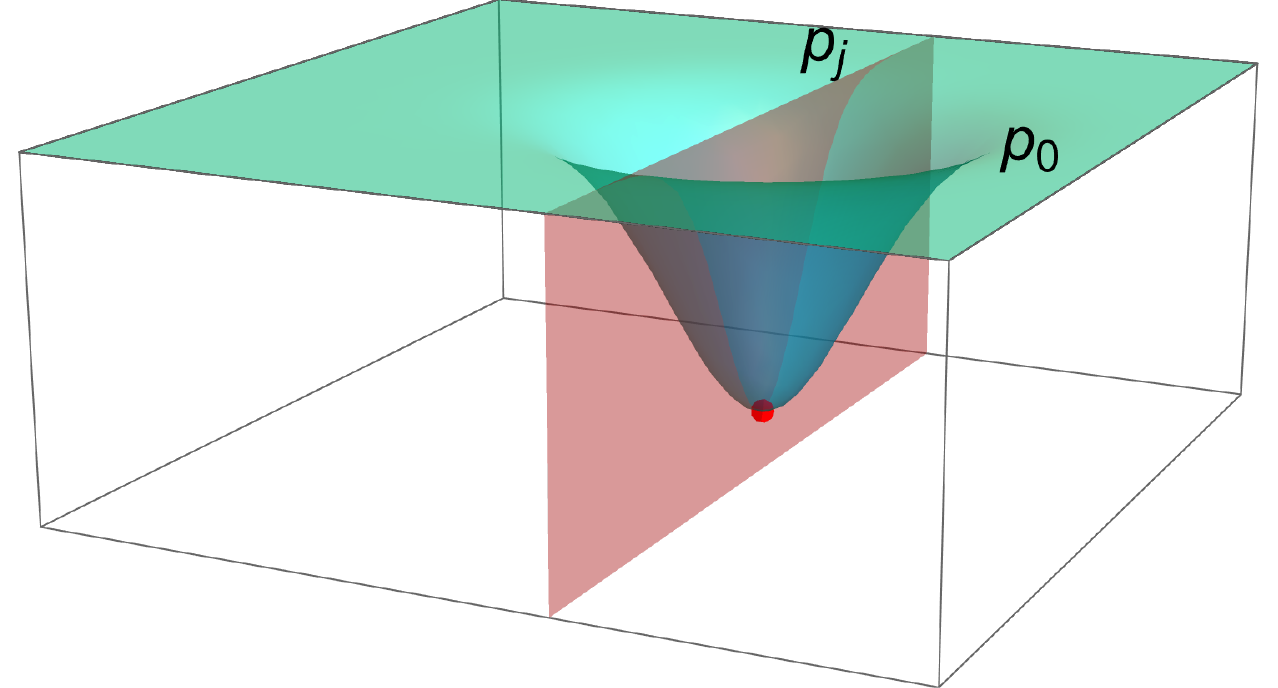}
\includegraphics[width=5cm]{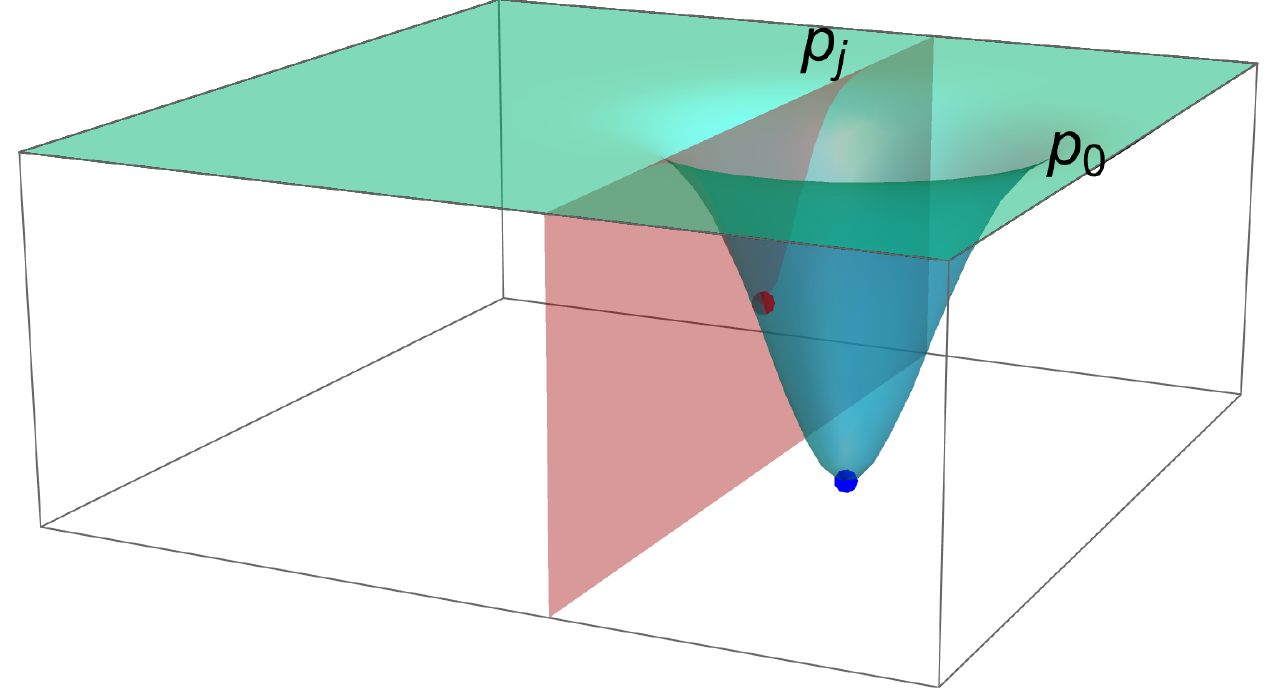}
\captionof{figure}{Example 2, ODE eq.(\ref{odep2}): Sketch that explains the behavior of $y_{2k+1}(x;p_0^{(1)*}(2k+1))$ while converging to $y(x;B_2^*)$. Left to Right are increasing $k$ values. Blue dots are the minimum of $d_1(2k+1,p_0,p_1,p_2,p_3)$. Red dots are the minimum of the restricted case $d_1(2k+1,p_0,0,0,0)$.
}\label{fig18}  
\end{center}
\begin{center}
\includegraphics[width=10cm]{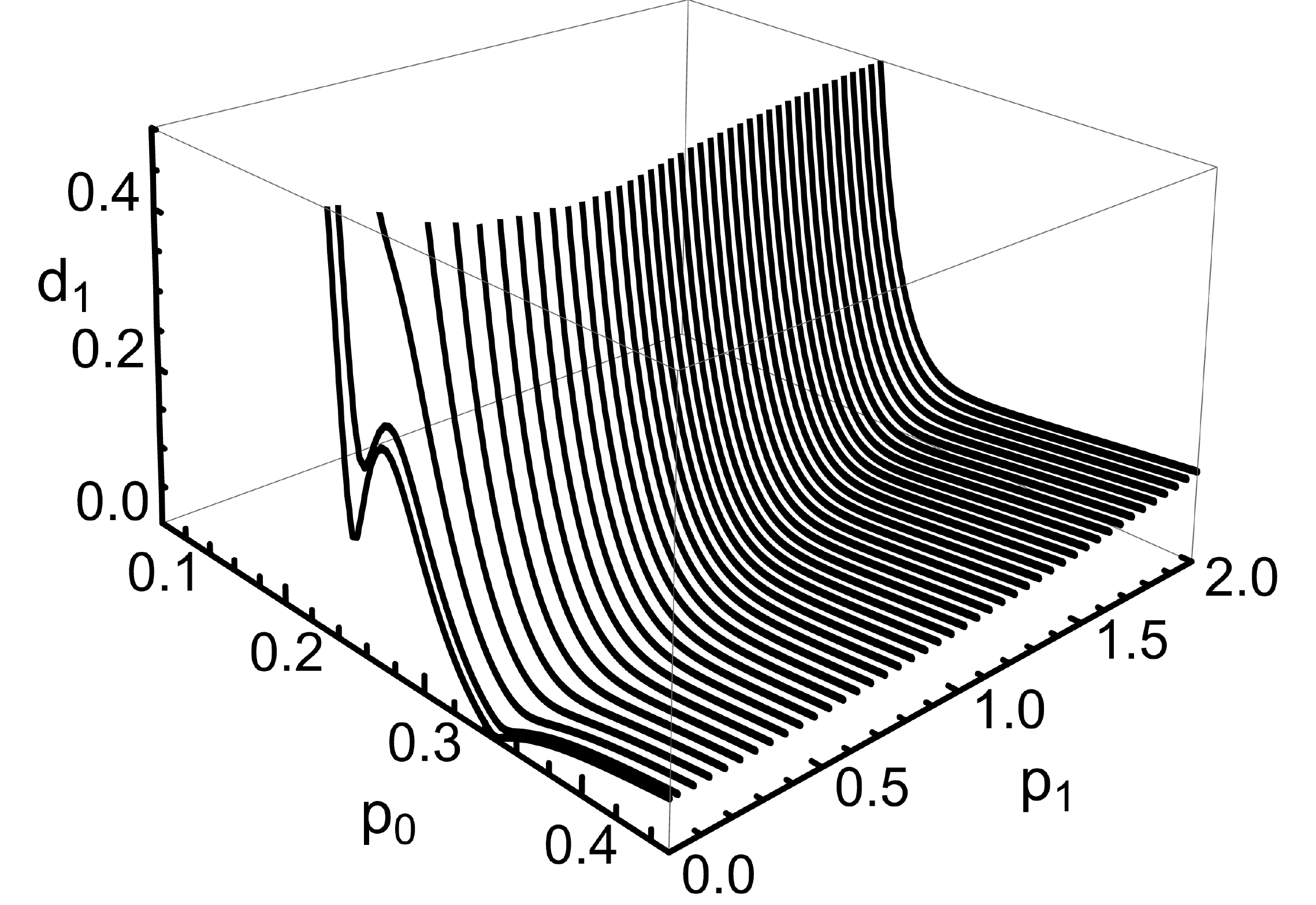}
\captionof{figure}{Example 2, ODE eq.(\ref{odep2}): $d_1(n=7;p_0,p_1)$ computed once obtained $y_n(x;p_0,p_1)$ from the general recurrence.
}\label{fig19}  
\end{center}
We have checked that this behavior has nothing to do with technical problems, such as precision when computing the minimums or the handling of the series expansions. We think that the problem is to use the assumption that  $p_1=p_2=p_3=0$ on the extended ODE (\ref{eode}). In general, the asymptotic minimum would be at $p^*=(p_0^*,p_1^*,p_2^*,p_3^*)$. However, when we restrict our method to live on the subspace $(p_0,0,0,0)$ the recurrence is unable to reach the real minimum and, at most, it is capable to be near it if $p_j^*$ $(j=1,2,3)$ are relatively near to zero. In Figure \ref{fig18} we show a sketch of the mechanism that may be the one responsible for the expansion's observed  conduct. For initial $k$ the distance have a minimum that becomes deeper as $k$ increases. Moreover, let us assume that it moves through (or near) the $\bar p\equiv(p_0,0,0,0)$ surface. During this part of the iteration (figure \ref{fig18} left), we would see how our restricted scheme also presents an increasingly deeper minimum. However, once the minimum $p^*$ has crossed the surface $\bar p$ (figure \ref{fig18} center) and tends to its limit, the restricted minimum worsens its distance to the solution (figure \ref{fig18} right). Assuming this picture we tried to go beyond $(p_0,0,0,0)$'s subspace but we had a relevant restriction when doing iteratively the algebraic integrals (it is out of the scope in this paper to implement an only-numerical scheme). We found that the ODE (\ref{odep2}) only permits such algebraic iteration for the subspace $(p_0,p_1,0,0)$. However, the ODE's symmetry makes that this extended subspace doesn't clarify that the proposed explanation is correct. In fact, we can prove that the minimum $(p_0^*(n),0,0,0)$ is locally stable on the extended $p_0,p_1,0,0)$ subspace:
\begin{equation}
\frac{\partial}{\partial p_1}d_1(n;p_0,p_1)\biggr\vert_{p_0=p_0^*;p_1=0}=0
\end{equation}
due to the fact that $y_n(x;p_0,p_1)=y_n(1-x;p_0,-p_1)$ and the structure of the ODE (\ref{odep2}). In figure \ref{fig19}  we show $d_1(n=7;p_0,p_1)$ where we see how the local minima $(p_0^*,0)$ are stable under $p_1$ perturbations and there are no trace of other minima.

We think that it is possible to design a numerical scheme to study the influence of parameters $p_2$ and $p_3$ on the convergence to the solution $B_2^*$. However, it is out of our goals in this paper, where we look just for algebraic methods. Let us remark that this method manages to get a first approximation to the solution. A similar practical problem arises when we have a meager convergence rate to the solution. This overall situation may appear in other ODEs. Therefore, we have two possible strategies to deal with this situation. If we need a precise algebraic approximation, we have to attempt a change of variables on the original ODE and apply the method to the transformed ODE. When we need a numerical, very precise result, we give below an approach that starts with the initial rough approximation and can improve it up to the desired precision.

\begin{center}
\includegraphics[width=7cm]{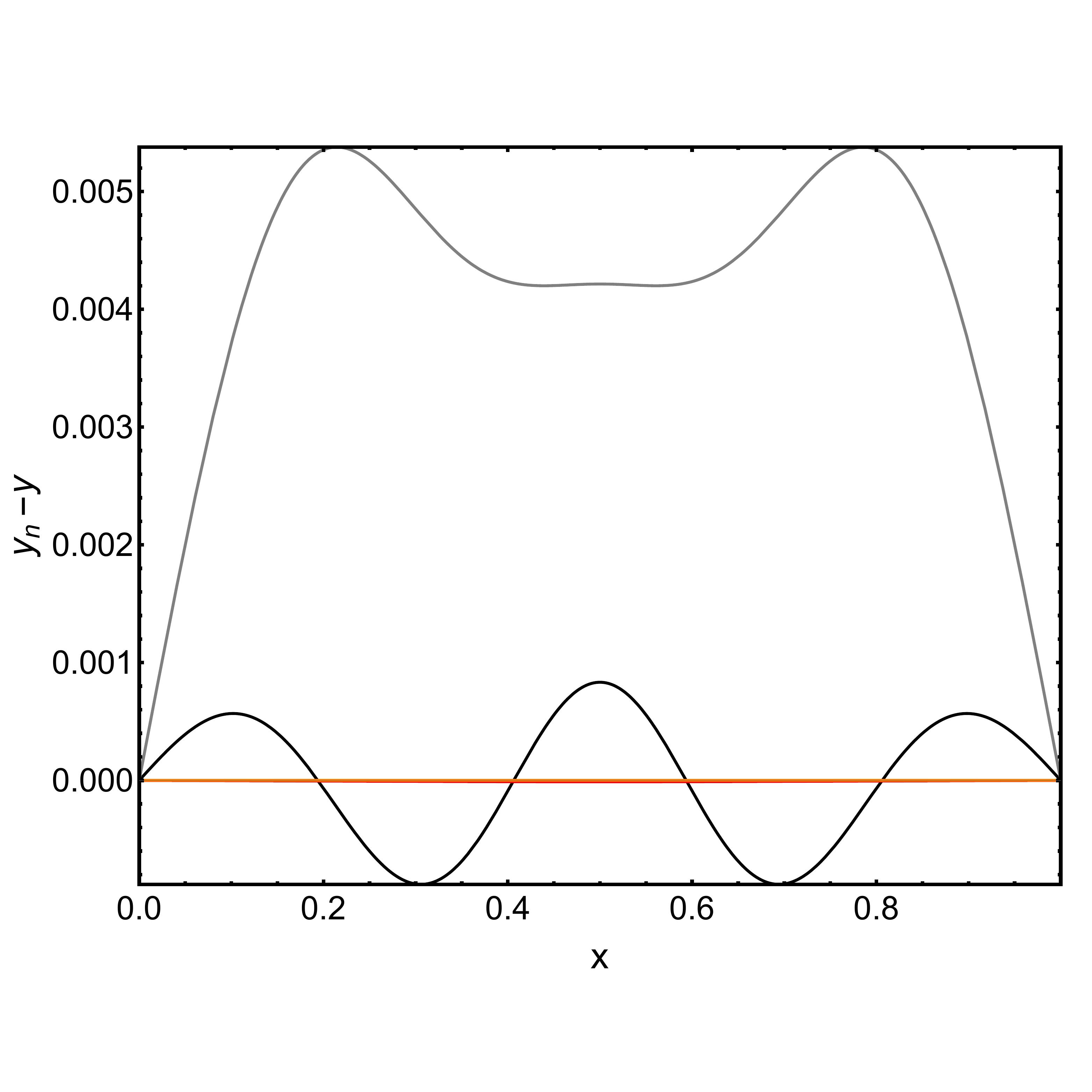}
\includegraphics[width=7cm]{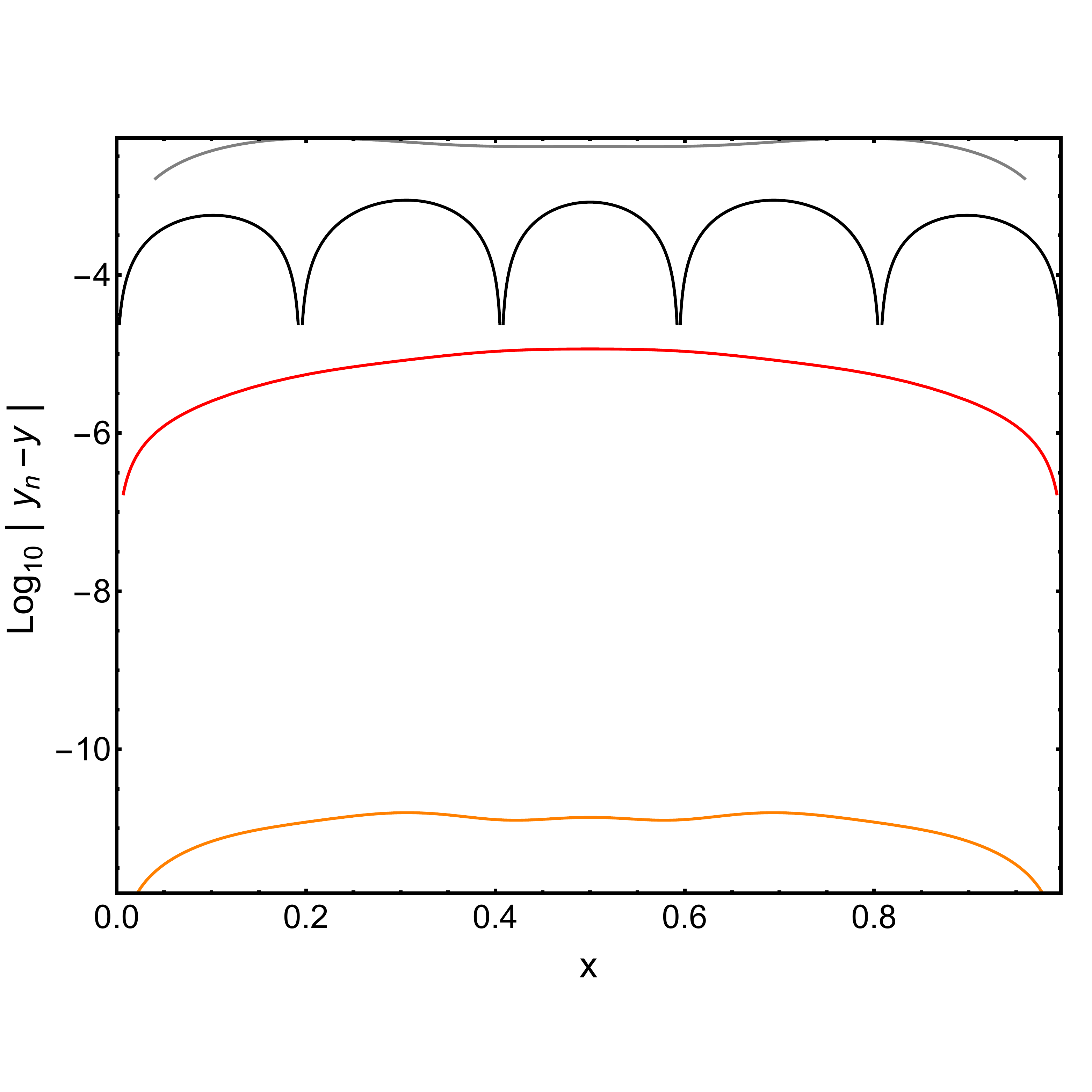}
\captionof{figure}{Example 2, ODE eq.(\ref{odep2}): Numerical convergence acceleration to $y(x;B_2^*)$: $y_{15}(x;p_0^{(1)*}(15))$ (black curve), $y_{17}(x;p_0^{(1)*}(17))$ (gray curve), $y_{15}(x;p_0^{(1)*}(15))+z_{1}(x)$ (red curve) and $y_{15}(x;p_0^{(1)*}(15))+z_{1}(x)+z_{2}(x)$ (orange curve). $z_{1,2}(x)$ are obtained by a numercial perturbation scheme (see main text).}\label{fig20}  
\end{center}
\item {\bf Accelerating the convergence to a solution:} We have seen that our method determines sequence of functions that converge to the ODE's solutions. However some times such convergence is slow as, for instance,  $y_{2k}(x;p_0^*(2k))$ where we need about $20$ iterations to decrease one order of magnitude the precision of our result. Moreover, the effect described for the convergence of $y_{2k+1}(x;p_0^{(1)*}(2k+1))$ towards $y(x;B_2^*)$, where there was an exponential convergence up to $n=15$ and an a  increment of the distance afterwards. We can systematically improve the approximate solution in either case by applying a straightforward, simple perturbation argument to the original ODE.

Let us assume that $y_n(x)$ is the approximate solution to $y(x)$. That is, it exists a function $z_n(x)$ such that
\begin{equation}
y(x)=y_n(x)+z_1(x)\quad,\quad z_1(0)=z_1(1)=0
\end{equation}
where it is assumed that $\equiv\vert z_1(x)\vert\leq \delta_1 <<\vert y(x)\vert$  $\forall x$. We substitute this decomposition on the original ODE (\ref{odep2}) and we get, up to order $\delta_1^2$ a second order differential equation for $z_1(x)$ :
\begin{equation}
z_1''(x)-y_n''(x)z_1(x)+y_n''(x)+e^{y_n(x)}=0
\end{equation}
This equation can be solved numerically by any simple routine as NDSolveValue in Mathematica. The resultant approximation $y_{n+1}(x)=y_n(x)+z_1(x)$ will have an error of order $\delta_1^2$. We can iterate the process by using $y_{n+1}(x)$ to improve the approximate solution. Let us remark that this straightforward algorithm stands on the existence of a good approximation whose distance to the exact solution is bounded uniformly on $x$ by a small parameter. Any approximation obtained by our scheme has this property. 
We show in figure \ref{fig20} the application of this perturbation iteration to the approximation $y_{15}(x;p_0^{(1)*}(15))$ that it is shown as a black curve. Just for sake of comparison we also show $y_{17}(x;p_0^{(1)*}(17))$ (gray curve) that separates from the exact value as we already commented above. The first perturbative correction, $z_1(x)$ is shown a a red curve. We see that this first correction  improves the precision of the approximation from about $10^{-4}$ up to $10^{-6}$. The second iteration, $z_2(x)$, improves much more the convergence up to about $10^{-11}$ for all $x$ values.

\item {\bf The effect of a change of variables:} The iterative scheme presented in this paper depends on the ODE's structure. Therefore, a change of variables may affect or not its rate of convergence or even its algebraic feasibility or simplicity. At this moment, we do not have any apriori argument to know what will happen after a change of variable. Let us take, for instance the ODE in this section (\ref{odep2}). We may think that our method's different convergence properties towards the $y(x;B_1^*)$ and $y(x;B_2^*)$ could be due to the solution's overall magnitude (for instance, the area between $y=1$ and the solution). It is a simple exercise to check that nothing new occurs if we do the change of variable $y\rightarrow 1/y$ to the original ODE and apply the scheme to the resulting ODE. 

We may apply other natural change of variables: $z=\exp(-y)$. The resulting ODE is in this case:
\begin{equation}
zz''-z-z'^2=0\quad,\quad z(0)=z(1)=e^{-1}
\end{equation}
The exact solutions (\ref{odep2e}) are re-writted:
\begin{equation}
z_{1,2}(x)=\frac{1}{2B_{1,2}^*}\cosh^2\left[B_{1,2}^*(x-\frac{1}{2})\right]\label{zsol}
\end{equation}
where $B_{1,2}^2$ have the values presented in figure \ref{fig10}. It seems that this ODE is, apparently, more complex that the original one. However, one sees inmediatly that the algebraic iteration of our method can be done for any set of $p$'s. In fact, our method determines the exact result for $n=0$ where we find two minima that make the distance $d_1=0$: 
$$(p_0^*,p_1^*,p_2^*,p_3^*)=\left(p_0,0,4p_0B_{1,2}^*,\frac{2p_0B_{1,2}^*}{\cosh(B_{1,2}^*/2)}\right)$$  
In other words, the differential equations (\ref{basicODE}):
\begin{equation}
\tilde z_0''+4B_{1,2}^*\tilde z_0+\frac{2B_{1,2}^*}{\cosh(B_{1,2}^*/2)}=0\quad,\quad \tilde z_0(0)=\tilde z_0(1)=e^{-1}
\end{equation}
have (\ref{zsol}) as solutions. This is a particular case because the $n=0$ iteration contains enough structure to have the same solution as the non-linear ODE we wanted to solve. Let us stress that we need the extended linear operator with four parameters to describe the non-linear equation correctly. A description with only the $p_0$ parameter won't be able to reproduce such a result. Finally, this example illustrates how a simple change of variables may simplify, in this case, the ODE's resolution.
\begin{center}
\includegraphics[width=5cm]{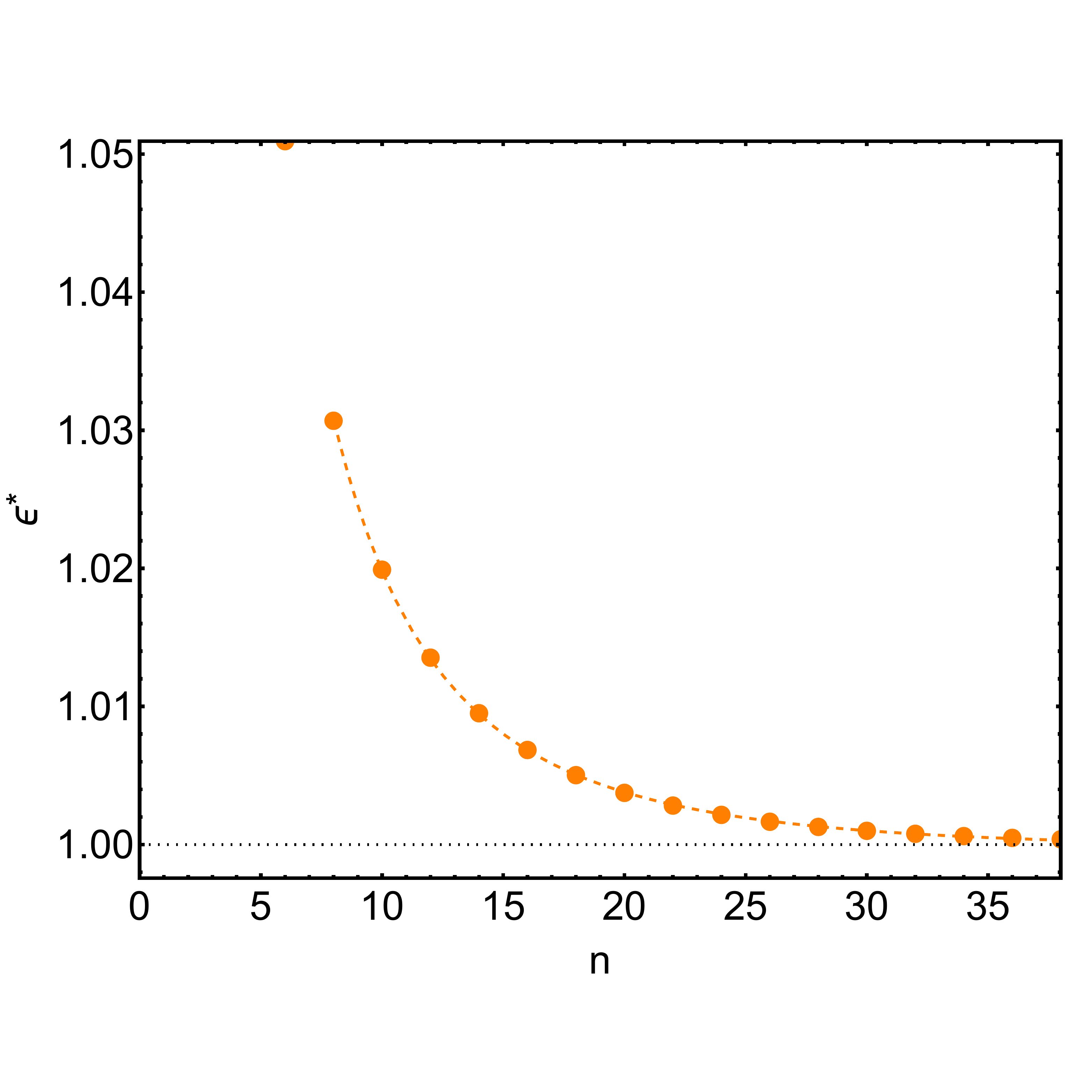}
\includegraphics[width=5cm]{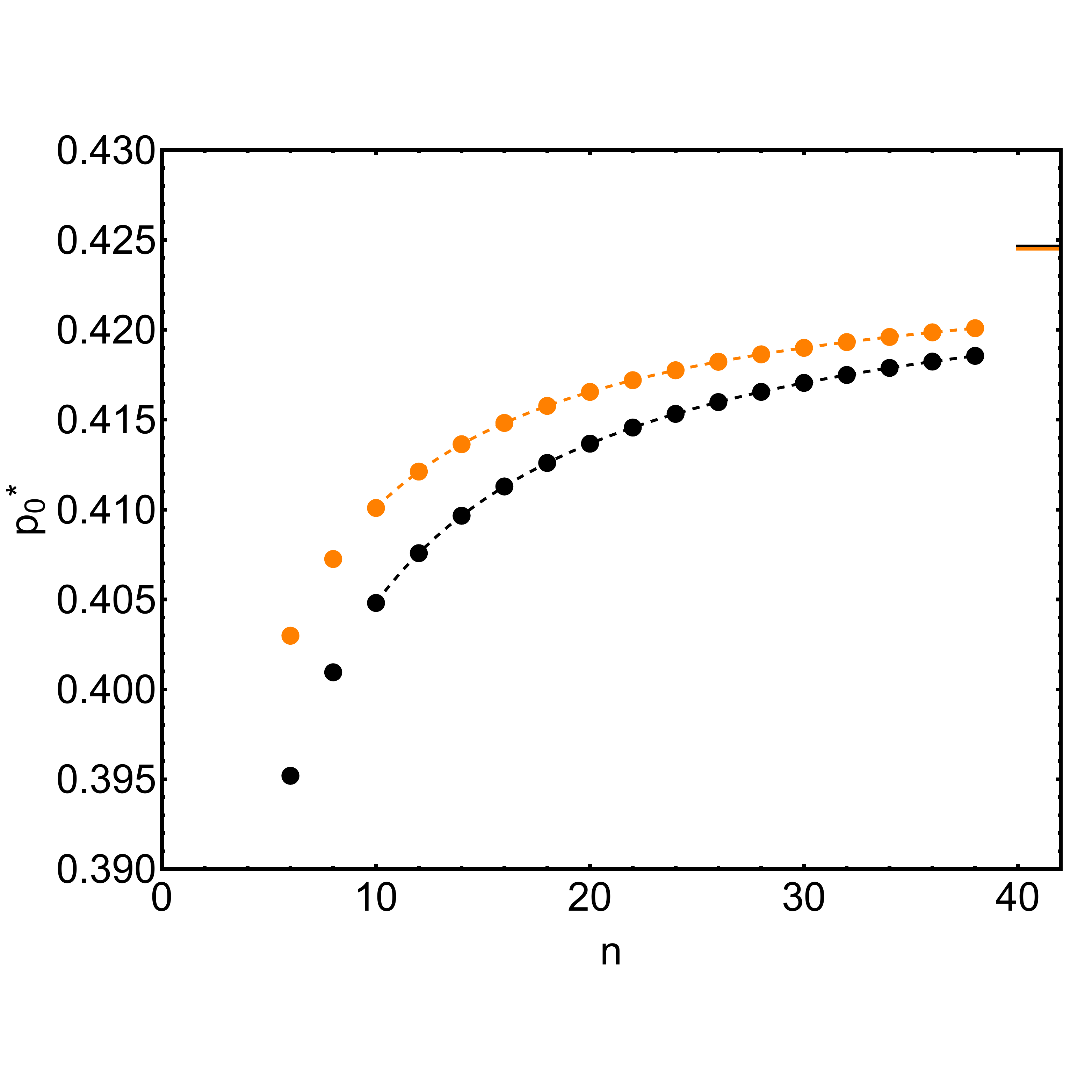}
\includegraphics[width=5cm]{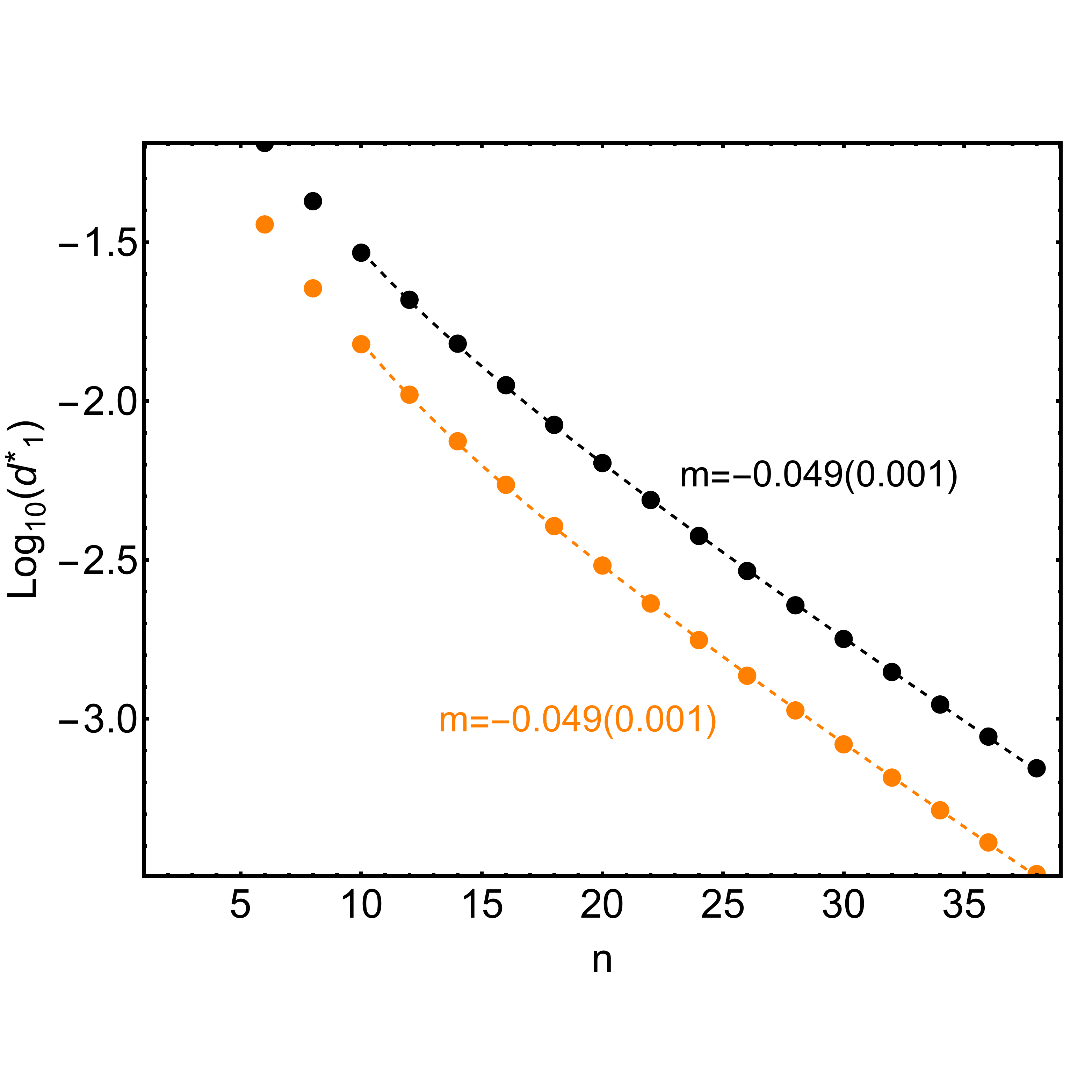}
\captionof{figure}{Example 2, ODE eq.(\ref{odep2}): Asymptotic behavior with the $n$-th even perturbative approximation of the values where the minima of $d_{1}(n,\epsilon,p_0)$ are located, $(\epsilon^*(n),p_0^*(n))$ (left and center figure respecitvely) and the distance's values at the minimum, $d_{1}^*$. Black dots in center and right figures are the results already presented in figure \ref{fig13} for the case $\epsilon=1$. Dashed lines are fits explained in the main text.}\label{fig21}  
\end{center}
\item {\bf \boldmath The $(p_0,\epsilon)$ minimization approach (Strong Conjecture):}
Let us comment here on what happens if we use $\epsilon$ as an extra parameter to be minimized at each iteration step in our scheme. As a typical example, we show in figure \ref{fig21} the sequence of minima values of the parameters $(\epsilon,p_0)$ for the case $n$-even (where we know that there is a unique minimum). We observe the expected behavior:
\begin{itemize}
\item $\epsilon\rightarrow 1$ as $n\rightarrow\infty$: We have fitted to the data the function $a_0+a_1/n+a_2n^2+a_3/n^3$ and we got $a_0=1.0006(0.0005)$. 
\item $p_0^*(n)$ converges to the same value when $n\rightarrow\infty$ as in the case with $\epsilon=1$: $p_0(\infty)=0.42451(0.00004)$ compared with $p_0(\infty)=0.42465(0.00002)$  for the $\epsilon=1$ case. 
\item The distance at the minimum, $d_1(n,\epsilon^*(n),p_0^*(n))$ decreases exponentially fast with $n$ with the same rate as in the $\epsilon=1$ case (see figure \ref{fig21}).
\end{itemize}
We see that the use of $\epsilon$ as an extra parameter to be minimised improves the precision of the approximation from the beginning of the iteration. Still, it does not improve the convergence rate. That it, it does not introduce new behavior to our scheme. Therefore, we may use it if the total CPU time to get the desired precision is improved that, from our point of view, it will depend on the problem we are solving. An initial better accuracy could be quickly compensated or not with the time needed to find a two dimensional minimum for large $n$-values.
\end{itemize}
\begin{center}
\vskip -0.5cm
\includegraphics[height=9cm]{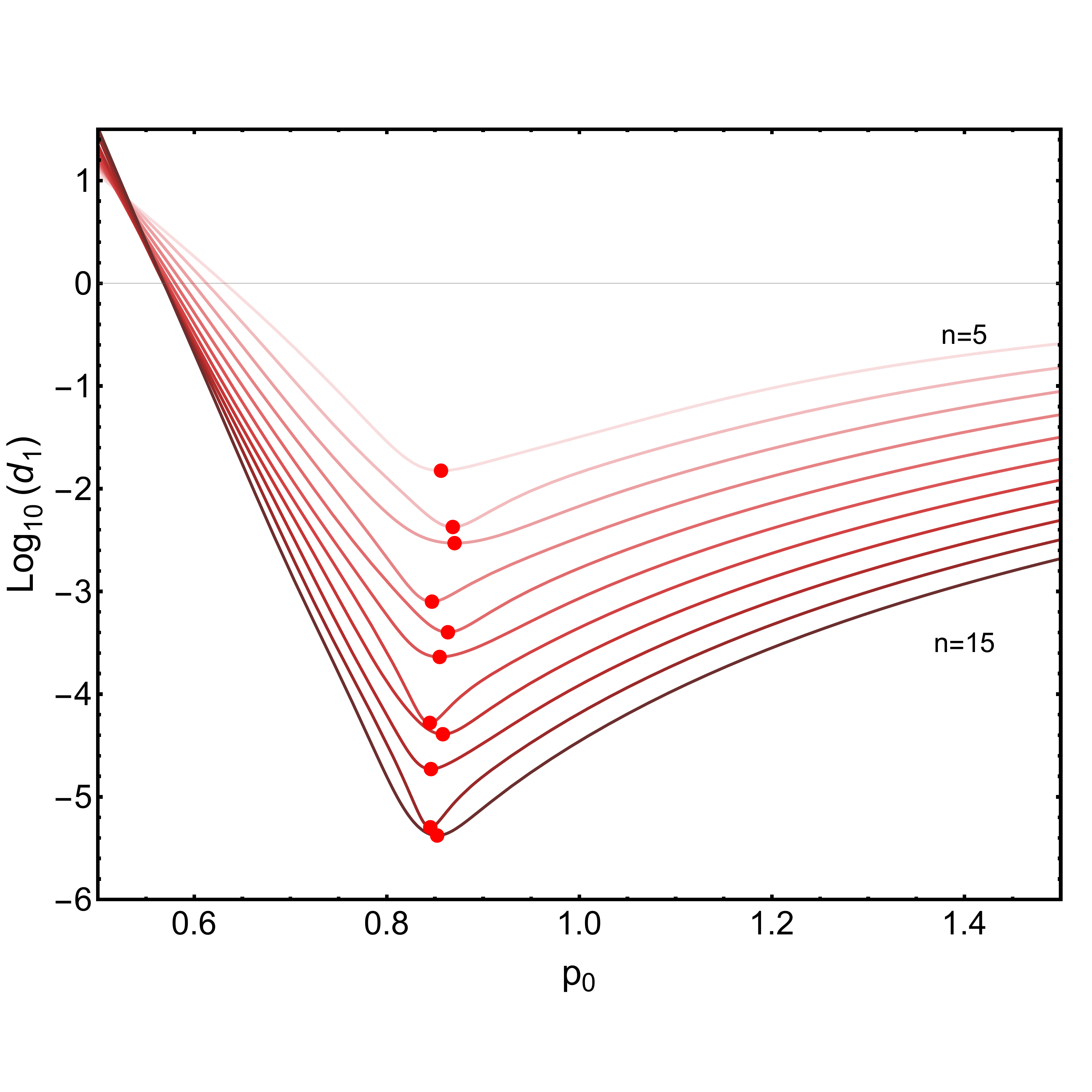}
\vskip-0.75cm
\captionof{figure}{Example 3, ODE eq.(\ref{odep3}): Decimal logarithm of distance $d_{1}(n;\epsilon=1,p_0)$ defined by eq.(\ref{mea1}) versus $p_0$ for each $n$th-perturbation approximation. The $y_n$ are obtained from the eODE's perturbative expansion (\ref{eode}). Cherry-tones curves increase intensity with $n$. We plot only $n=5,\ldots ,15$. Red dots shows the unique minima.}\label{fig22}  
\end{center}
\begin{center}
\vskip -0.5cm
\includegraphics[height=7cm]{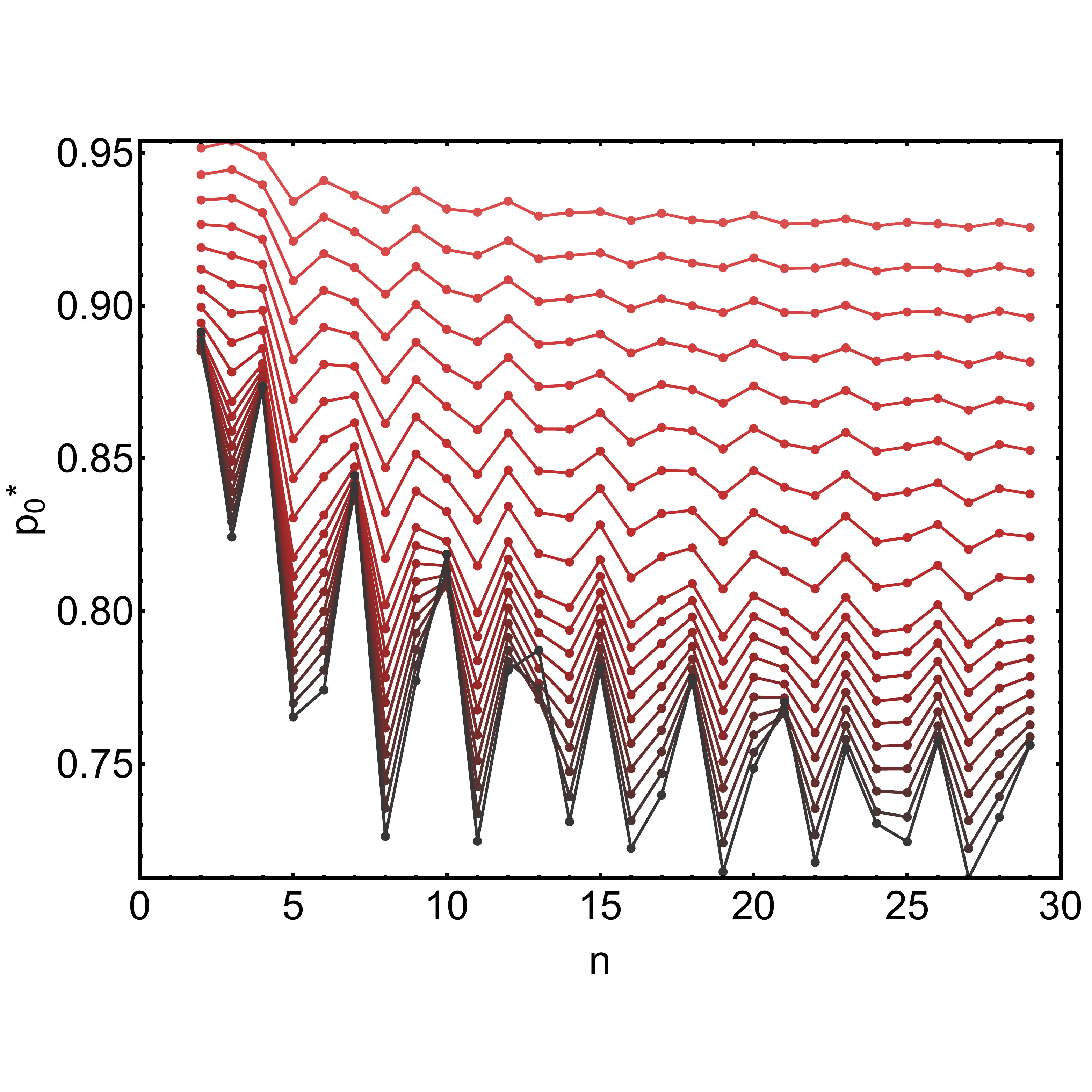}  
\includegraphics[height=7cm]{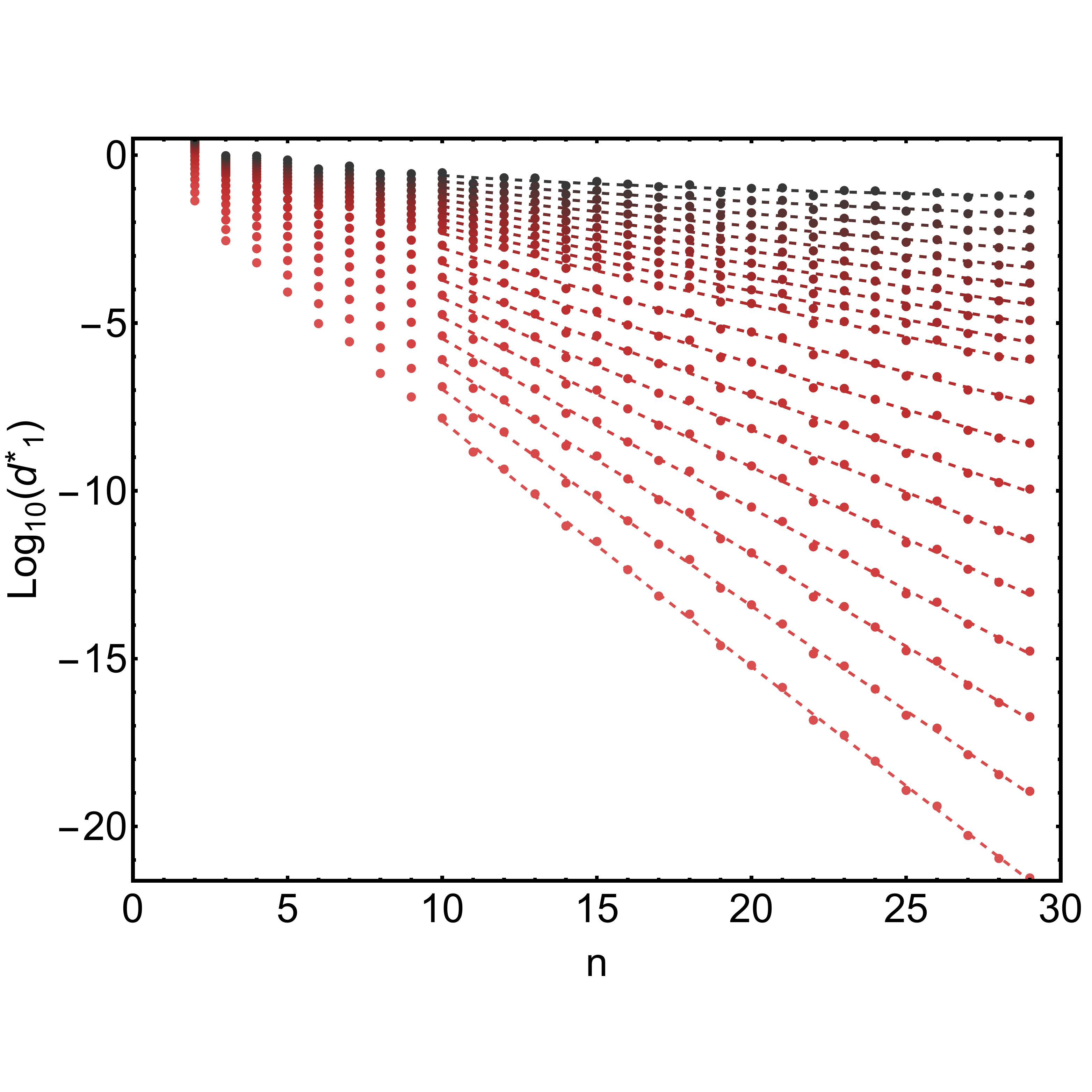}
\vskip -0.5cm
\captionof{figure}{Example 3, ODE eq.(\ref{odep3}): Asymptotic behavior with the $n$-th perturbative approximation of the values where the minima of $d_{1}$ are located, $p_0^*$ (left figure) and the values of such distances at the minima, $d_{1}^*$ (right figure). Each cherry tone represents a value of $\xi=1,1.2, 1.4, 1.6,1.8, 2, 2.2, 2.4, 2.6, 2.8, 2.9, 3.0, 3.1, 3.2, 3.3, 3.4, 3.5, 3.6, 3.7$ from red to dark red. Dashed lines are asymptotic fits explained in the main text.} \label{fig23}  
\end{center}
\section{EXAMPLE 3: \boldmath $y''+\xi (y'+y^2)=0$ (BVP)}

The differential equation we study is:
\begin{equation}
y''+\xi(y'+y^2)=0\quad ,\quad y(0)=0\quad, y(1)=1\label{odep3}
\end{equation}
where there are no known exact analytical solutions. Moreover, It seems from numerical computations that there is one solution for $\xi\in[0,\xi_c\simeq 3.73]$ and none when $\xi>\xi_c$. We want to check the behavior of our method for this type of parametric transition. 

 We focus on the restricted conjecture ($\epsilon=1$) and the free parameter's space $(p_0,p_1,p_2,p_3)=(p_0,0,0,0)$. We show on figure \ref{fig22} the form of the distance function $d_1(\epsilon=1,p_0)$ as a function of $p_0$ for $\xi=2.2$ and iteration orders $n=5,6,\ldots,15$. We observe the existence of an unique minimum that gets deeper as we increase $n$. However, the values where the minima are located, $p_0^*(n)$, do not follow a clean and systematic convergence pattern. That is confirmed on figure \ref{fig23} (left) where we plot $p_0^*(n)$ vs $n$ for different $\xi$ values. We see how they seem to fluctuate around an asymptotic value. Such ``fluctuation'' increases as $\xi$ increases, but it is always bounded (observe in the figure that the maximum fluctuation for $\xi=3.7$ is of order $0.04$ around $p_0=0.74$). Nevertheless, the value of the distance at the minimum, $d_1^*(n)$ has a more normal behaviour.  We see in figure \ref{fig23} that for each $\xi$ it decreases exponentially fast for large enough $n$. We show there by dashed lines the fits of the data to the function $m(\xi) n+a_1(\xi)+a_2(\xi)/n+a_3(\xi)/n^2$ . We see how the curves become more horizontal as we increase $\xi$ from $1$ (the fastest decay) up to $\xi=3.7$. We checked that for $\xi=3.8$ and beyond, the distance function $d_1(n;p_0)$ doesn't develop any minimum (for finite $p_0$-values). Therefore, the parameter $m(\xi)$ is helpful to locate the critical value $\xi_c$ where it separates the region with one ODE's solution to none. In figure \ref{fig24} we show the computed behavior of $m(\xi)$ vs $\xi$. From it, we can estimate the critical $\xi$ by fitting a third-order polynomial to the data and then solving the equation $m(\xi_c)=0$. We obtain $\xi_c\simeq 3.76818$ that is coherent with the numerical computations we have done. We may conclude that our scheme permits a systematic analysis of the existence of solutions of a nonlinear ODE depending on their parameter values. 
 \begin{center}
\vskip -0.5cm
\includegraphics[height=9cm]{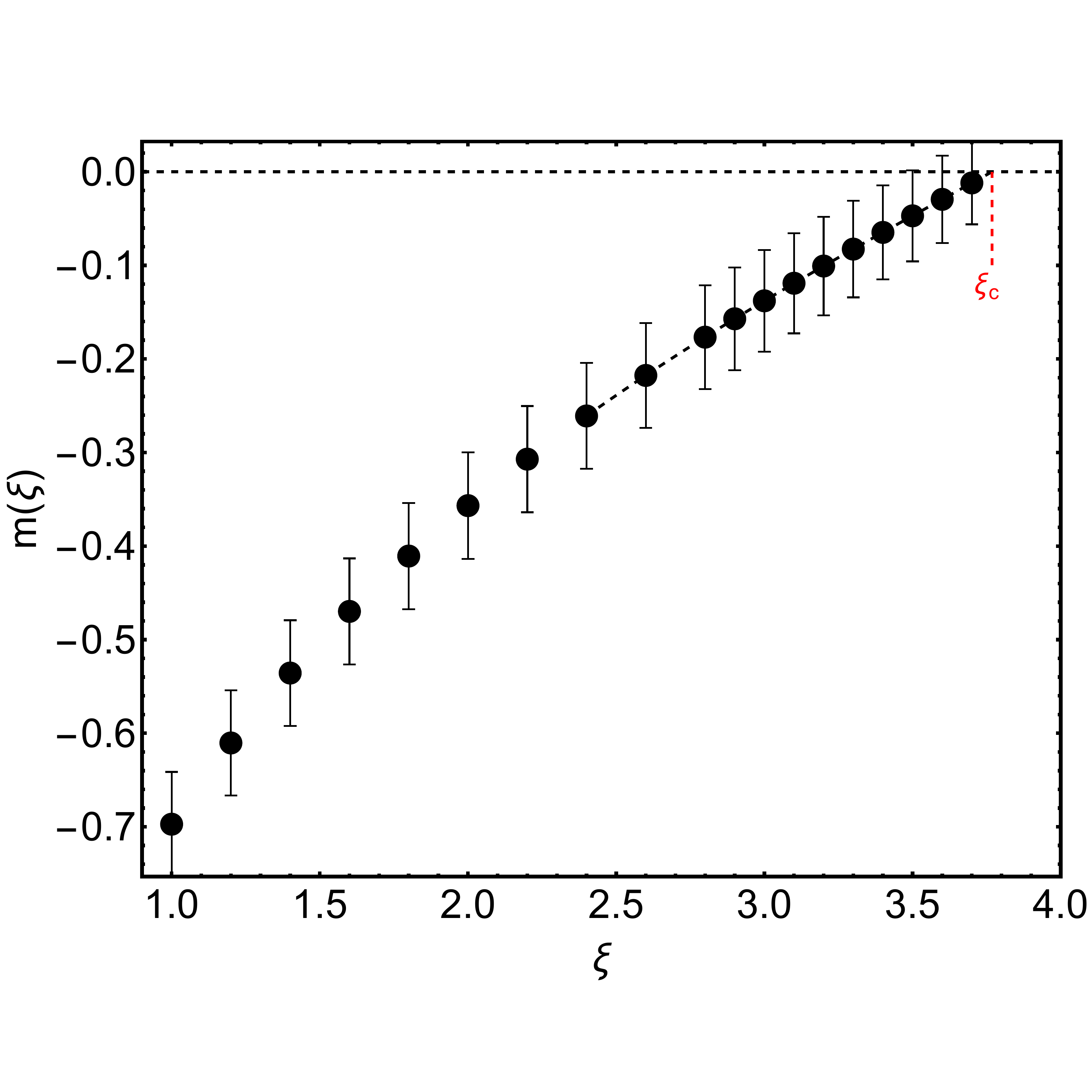}  
\vskip -0.5cm
\captionof{figure}{Example 3, ODE eq.(\ref{odep3}): Asymptotic behavior of $d_1^*(n;\xi)\simeq 10^{-m(\xi)n}$ for $n$ large enough. It is plotted $m(\xi)$ vs $\xi$.  Dashed line is a polynomial fit to the data that is used to estimate the critical $\xi$ value by solving $m(\xi_c)=0$. We obtain $\xi_c\simeq 3.76818$.} \label{fig24}  
\end{center}
  \begin{center}
  \begin{tabular}{ l | c | c }
    \hline
    \boldmath $\xi$&{\bf NDSolveValue}&{\bf This work}\\ \hline
    1.0&-1.3917&-21.5365\\ \hline
    2.0&-0.9504&-11.4222\\ \hline
    2.6 &-0.7554&-7.2920 \\ \hline
       3.0 &-0.6360&-4.9195 \\ \hline
       3.7 &-0.4031&-1.1837 \\
    \hline
  \end{tabular}
  \captionof{table}{Example 3, ODE eq.(\ref{odep3}): $Log_{10}(d_1)$ for the solution obtained using NDSolveValue Mathematica's routine and when using the scheme proposed in this paper with $n=29$.} \label{table}  
\end{center}
Finally, we use this example to compare our method with a well-known numerical routine quickly.  
 We show in figure \ref{fig25} the approximate solutions for $n=29$ when $\xi=2.2$ and $3$. We compare them with the numerical solutions obtained by using NDSolveValue routine from Mathematica. At a glance, there are no differences. They appear when we measure the distance $d_1$ to both approximations (see Table \ref{table}). There,  we observe that the precision of our result is systematically the best, with a difference of twenty orders of magnitude for the $\xi=1$ case. We could ask the Mathematica's routine to improve its precision or use our method with larger $n$, but these results give us some taste of how our scheme behaves.  Finally, let us remark that the exponential decay of $d_1^*(n)\simeq \delta^n$ implies that a ghost expansion of the solution exists as we explicitly commented in the previous examples.

\begin{center}
\vskip -0.5cm
\includegraphics[height=9cm]{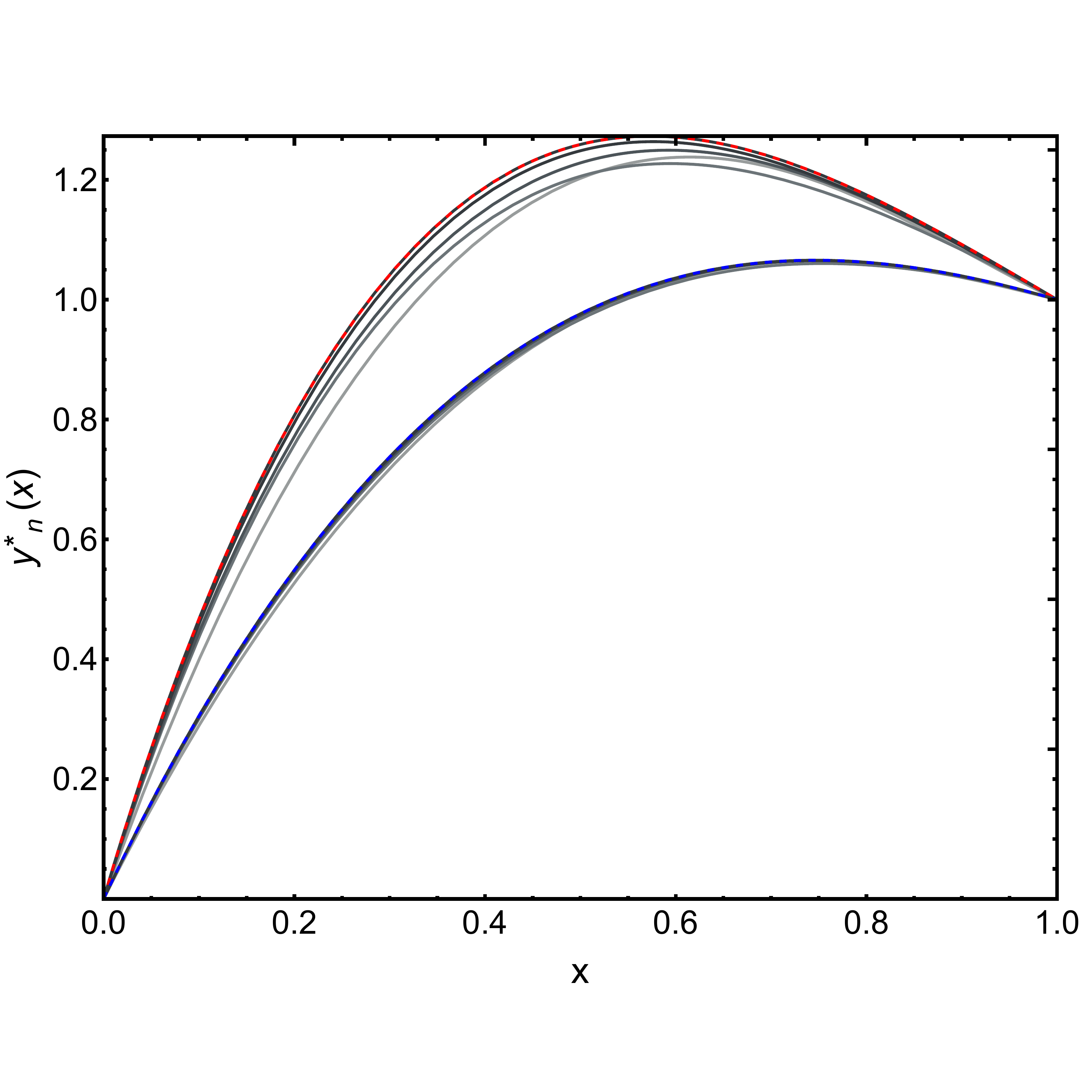}  
\vskip -0.5cm
\captionof{figure}{Example 3, ODE eq.(\ref{odep3}): Approximate solutions $y^*(x)\equiv y_{n}(x;p_0^*(n))$ vs $x$ for $\xi=2.2$ (bottom curve) and $\xi=3$ (top curve). Black curve correspond to $n=29$. Some gray curves are also shown that represent $n=2,3,4,5$ for each case. Red and blue curves are solutions computed with Mathematica's NDSolveValue routine.} \label{fig25}  
\end{center}
\section{Example 4:  The Lane-Emdem equation \boldmath $y''+2y'/x+y^m=0$ (IVP)}
Our method to build a sequence of approximated functions can be applied to Initial Value Problems (IVP) (see eq.\ref{solIVP}). The computational effort to compute $y_n$ is the same as in the BVPs, and it depends only on the ODE's structure. 
However, we now have to introduce an arbitrary fix interval $T$   where $x\in[0,T]$.  $T$ is a relevant external parameter. Our method attempts to find an approximation to the solution on all the ``time interval'' $T$ by looking for the best $p^*(n;T)$  that minimises the distance $d_1$ on such $T$-interval at each $n$-iteration level.  As we will see, the convergence rate to the exact solution typically depends on $T$  and therefore, we could need more iterations as $T$ increases to reach a certain precision level. 

We use the Lane-Emdem differential equation to understand the behavior our our method when applied to this type of problems. This equation is given by:
\begin{equation}
y(x)''+\frac{2}{x}y(x)'+y(x)^m=0\quad,\quad y(0)=1\quad,\quad y'(0)=0\label{odep4}
\end{equation}
with $m=0,1,\ldots$. There are known exact solutions for this equation for the $m$-values: 
\begin{eqnarray}
m=0&:& y(x)=1-\frac{x^2}{6}\nonumber\\
m=1&:& y(x)=\frac{\sin x}{x}\nonumber\\
m=5 &:& y(x)=\frac{1}{\sqrt{1+x^2/3}}
\end{eqnarray}
that we will use to check the behavior of our method. 

We first compute the sequence of approximations given by eq. (\ref{y1pert2}) with (\ref{solIVP}) for $m=0, 1,\ldots, 6$. For simplicity we will restrict ourselves to the case $\epsilon=1$ (restricted conjecture) and $p=(p_0,0,0,0)$. 
\begin{center}
\vskip -0.5cm
\includegraphics[height=7cm]{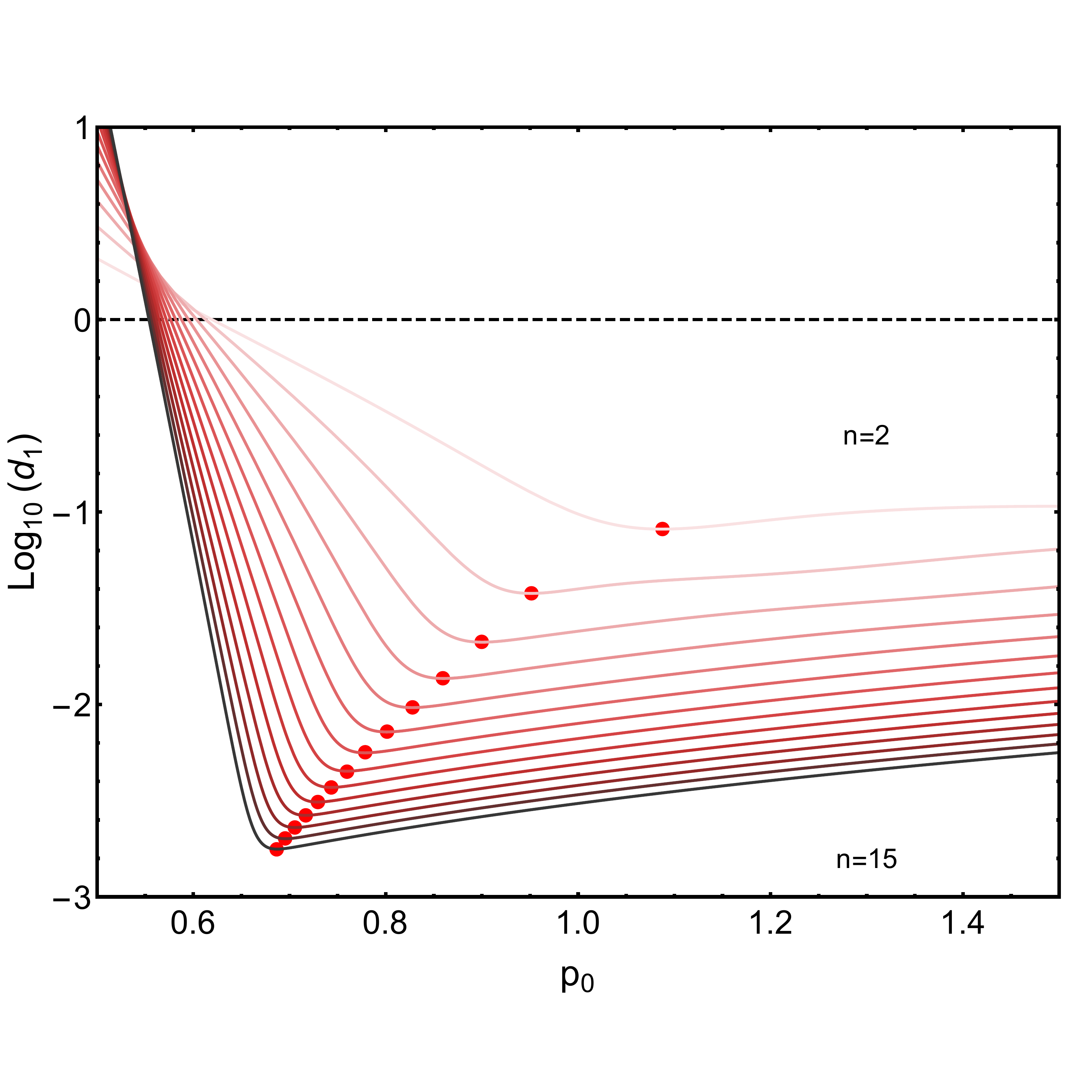}  
\includegraphics[height=7cm]{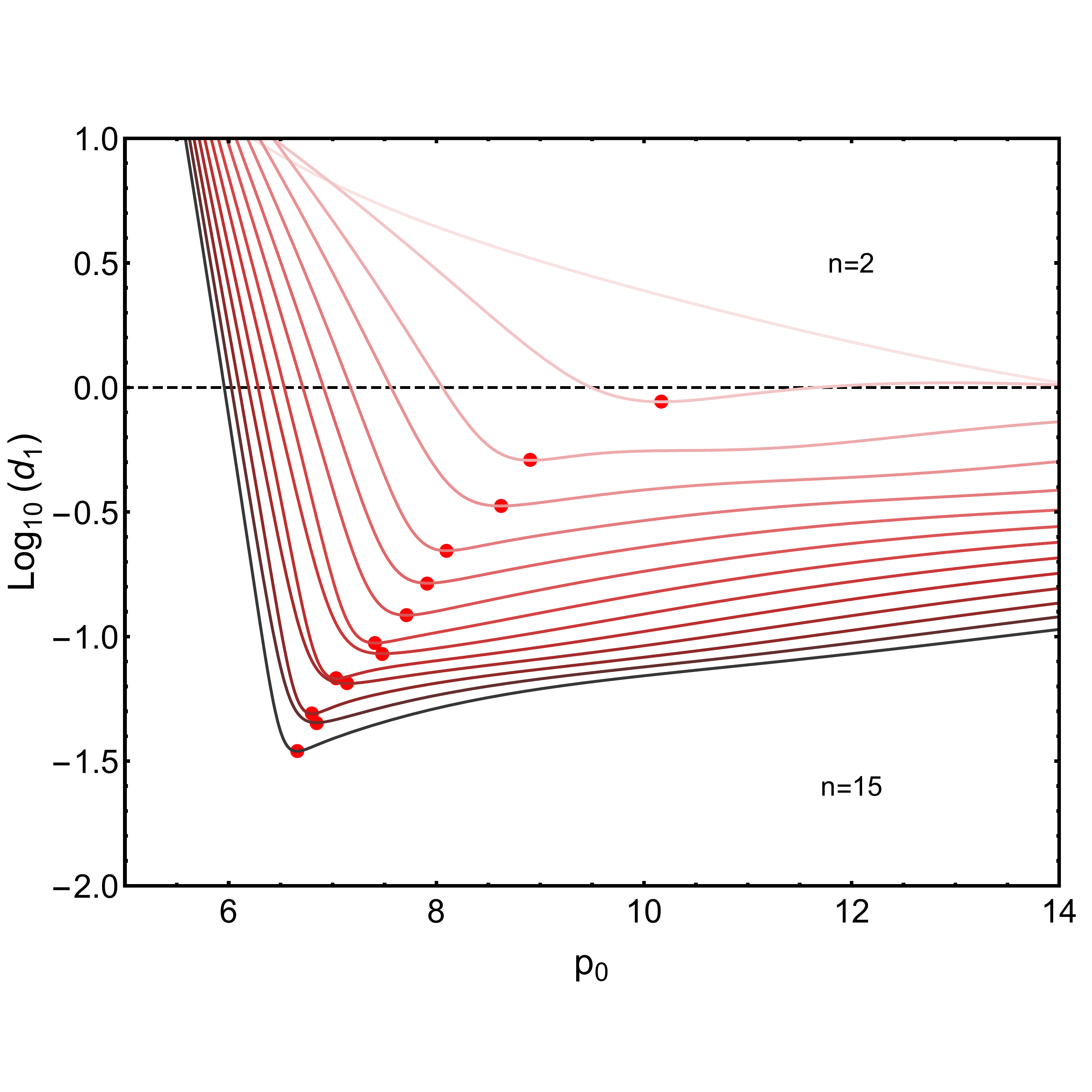}  
\vskip -0.5cm
\captionof{figure}{Example 4, ODE eq.(\ref{odep4}): $\log_{10}d_1(n;p_0)$ vs. $p_0$ for $n=2, 3,\ldots, 15$ and $m=0$, $T=1$ (left) and $m=2$, $T=5$ (right).} \label{fig26}  
\end{center}
 \begin{center}
\vskip -0.5cm
\includegraphics[height=7cm]{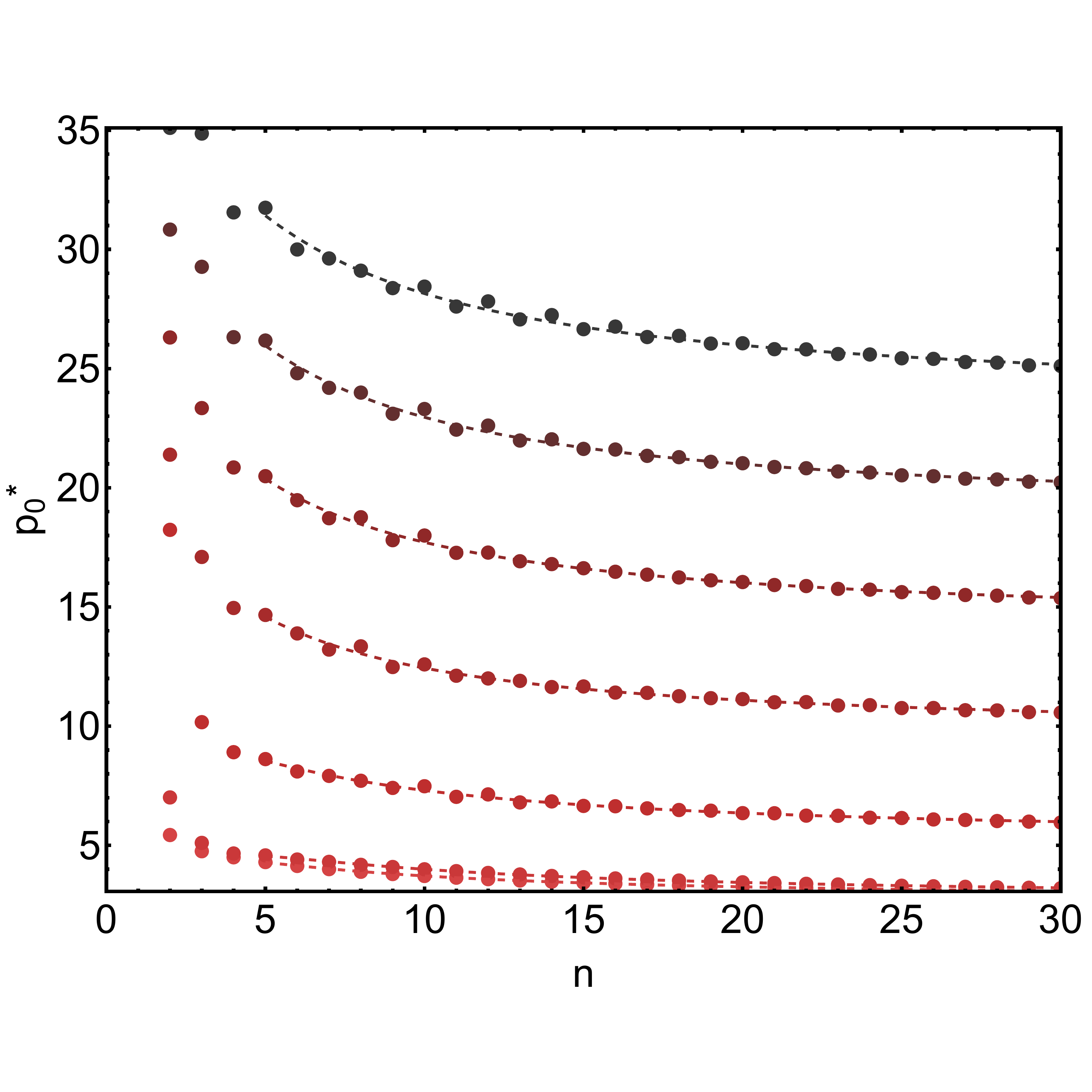}  
\includegraphics[height=7cm]{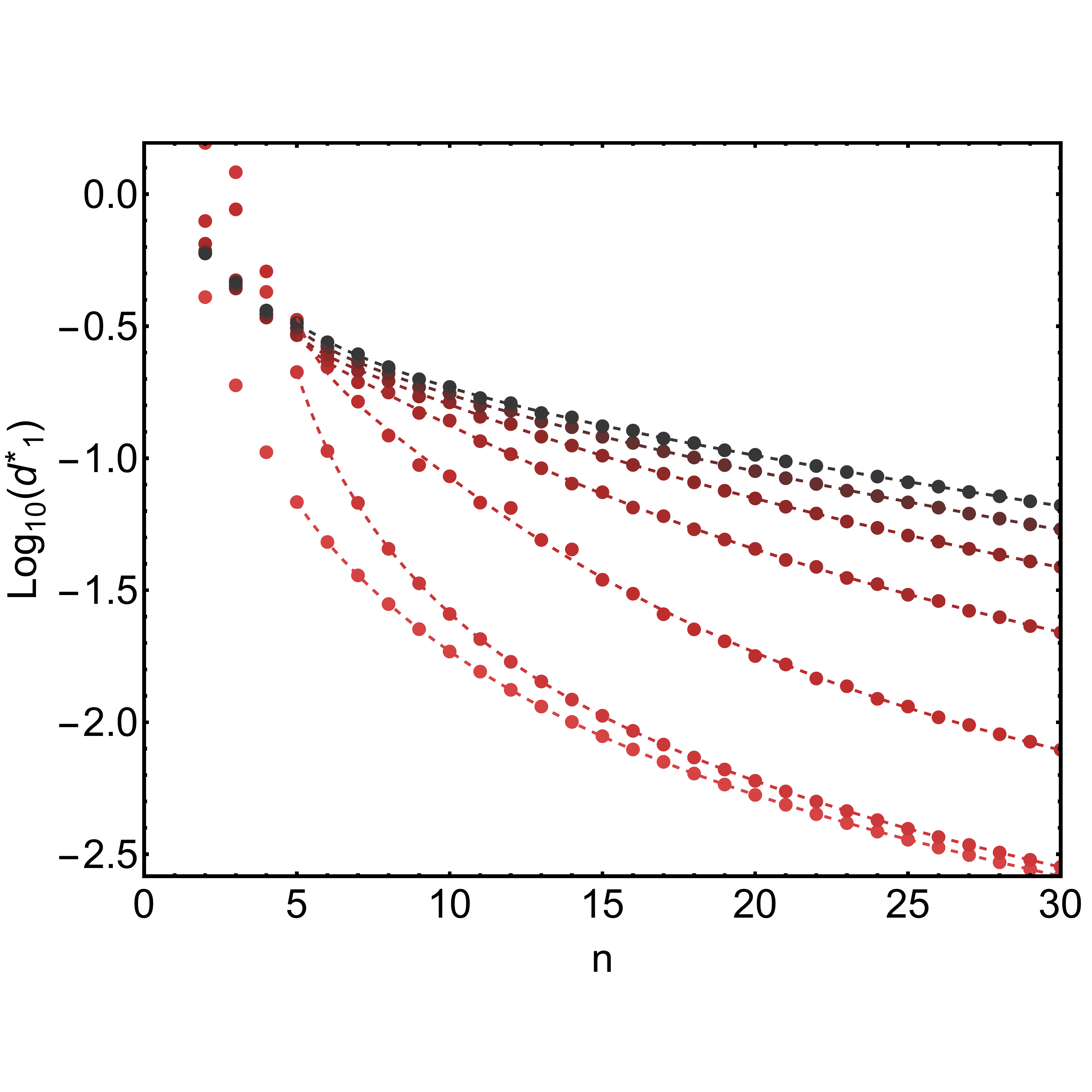}  
\vskip -0.5cm
\captionof{figure}{Example 4, ODE eq.(\ref{odep4}): $p_0^*(n,T)$ vs. $n$ (left) and $\log_{10}d_1^*(n,T)$ vs $n$ (right) for $m=0, 1, 2,\ldots, 6$ (curves from cherry tone to black) for $T=5$. Dahsed lines are fits (see text).} \label{fig27}  
\end{center}
\begin{itemize}
\item {\boldmath $d_1(n;p_0)$\bf , the minima and the serie's convergence:} We observe in figure \ref{fig26} the typical behaviour for the distance $d_1(n;p_0)$ vs. $p_0$ for two cases: $m=0$ and $T=1$ and $m=2$ and $T=5$. We see how there is a well defined local minimum that is getting deeper as $n$ increases. 
We compute the expansion, the minima $p_0^*(n,T)$ and the distance at the minima, $d_1^*(n,T)$ for $m=0, 1,\ldots,6$ and $T=5$ as we show in figure \ref{fig27}. We did the following fits to the data for each $m$: $p_0^*(n,T)=p_0^*(T)+a_1(T)/n+a_2(T)/n^2$ and $d_1^*(n,T)=b_0(T)n+b_1(T)+b_2(T)/n+b_3(T)/n^2$ for $n\in[10,30]$.  The distance's exponential decay rate, $\delta=10^{b_0(5)}$ such that $d_1^*(n)\simeq \delta^n$, and the asymptotic ($n\rightarrow\infty$) $p_0^*(5)$ for each ODE's $m$ are:
  \begin{center}
  \begin{tabular}{| c || c | c |}
    \hline
    {\boldmath $m$}&{\boldmath $\delta$}&{\boldmath $p_0^*$}\\ \hline
    0&0.970(0.002)&2.69(0.02)\\ \hline
    1&0.970(0.004)&2.71(0.03)\\ \hline
    2 &0.99(0.04)&5.2(0.2) \\ \hline
       3&0.97(0.01)&9.6(0.3)\\ \hline
       4&0.966(0.009)&14.1(0.3) \\ \hline
         5&0.964(0.007)&18.7(0.2) \\ \hline
           6&0.966(0.007)&23.5(0.5)\\ \hline
    \hline
  \end{tabular}
 \end{center}
That is, our method works, and it exists a sequence of functions that converges to the exact solution on the given $T$-interval. Moreover, the scheme has an exponential decay that permits the definition of a ghost expansion of the solution for a given $T$. Nevertheless, we observe that the distance decay rate, $\delta$, is near one in almost all the cases studied, making it difficult to get good practical approximations of the solutions. One way to improve the decay rate is to attempt a change variables.  As we already commented, we do not have any preliminary recipe that guarantees any better behaviour. We think much more rigorous formal work is needed to understand the mechanisms that may accelerate the convergence under a change of variables as the control of singular points on the ODE ($x=0$ in this example).

\begin{center}
\vskip -0.5cm
\includegraphics[height=7cm]{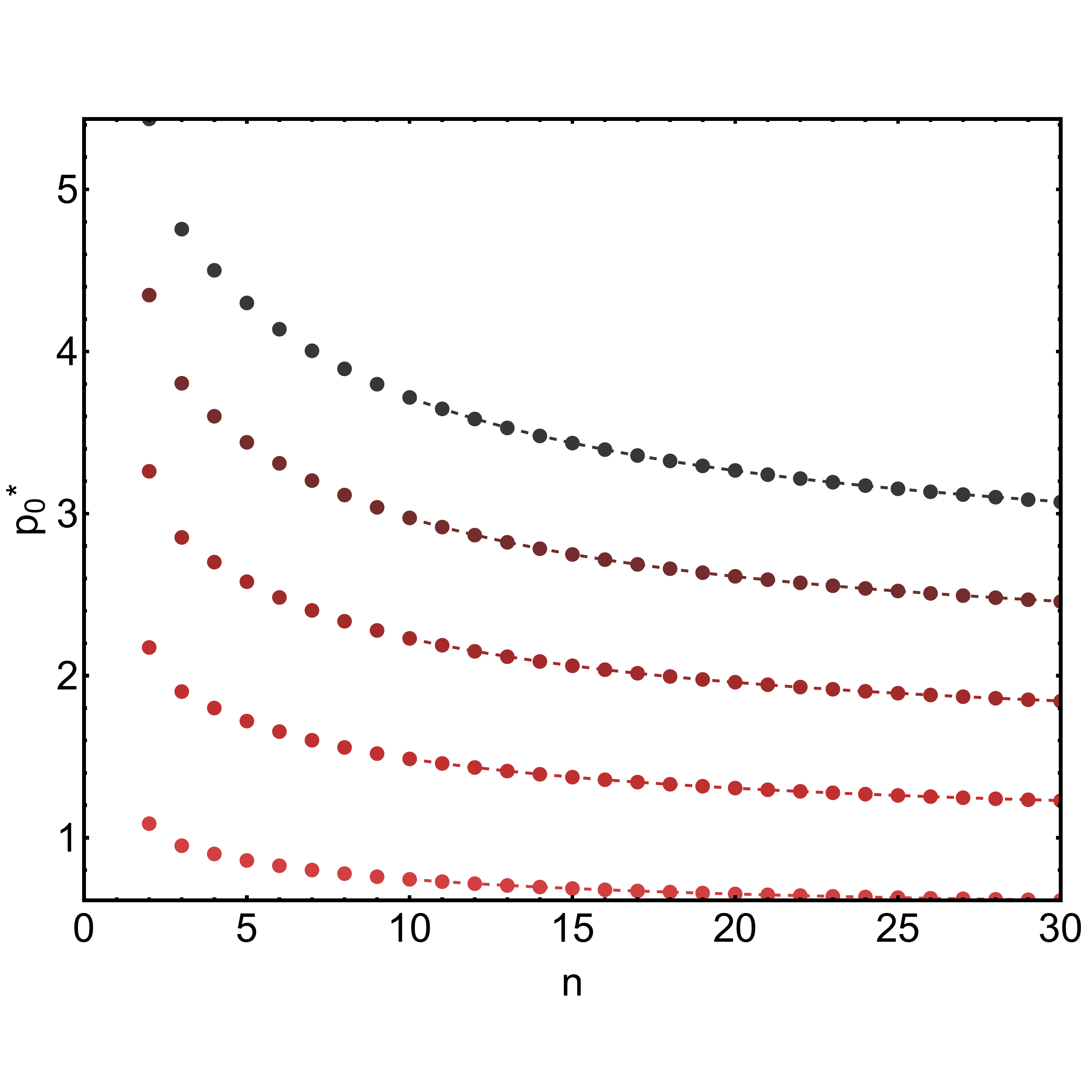}  
\includegraphics[height=7cm]{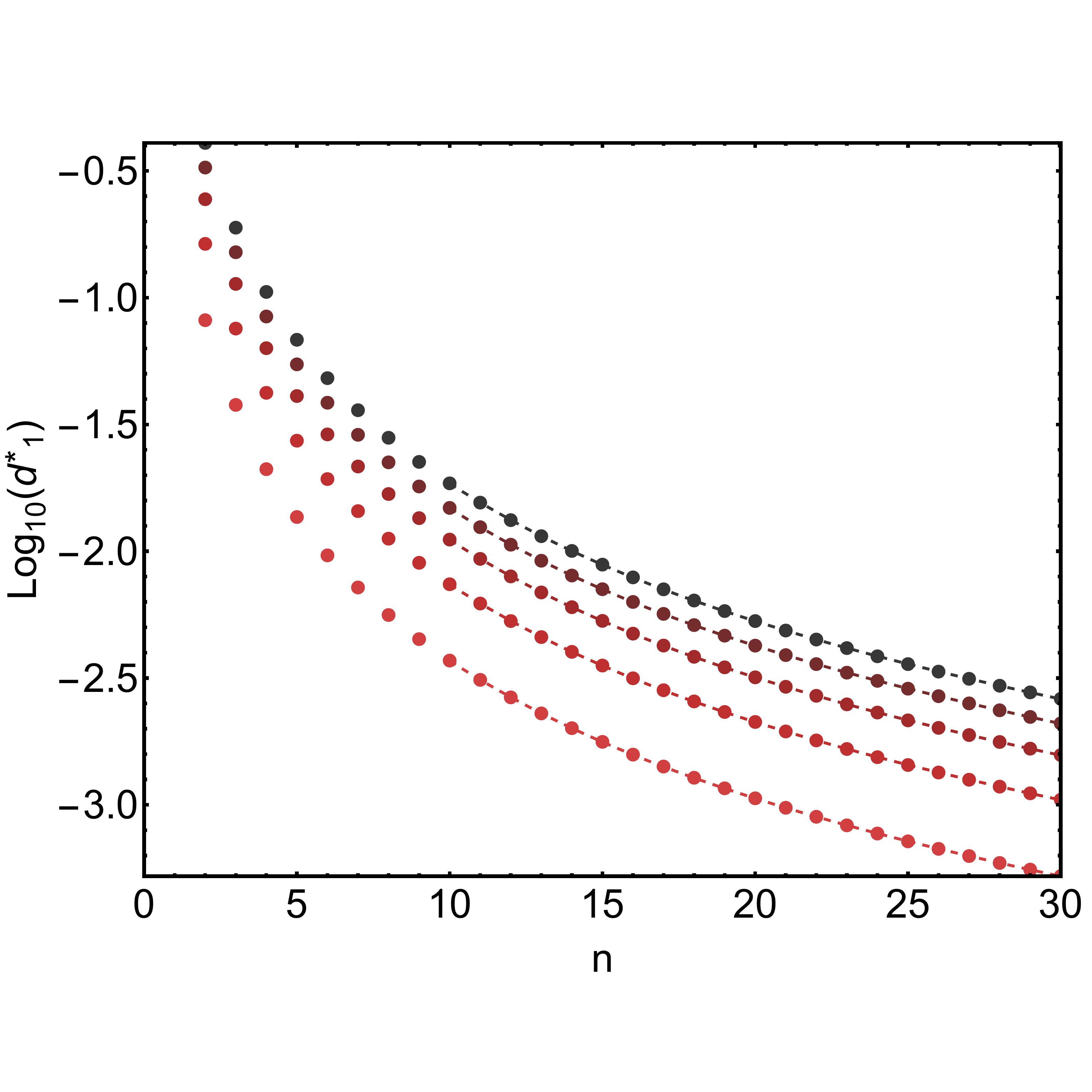}  
\vskip -0.5cm
\captionof{figure}{Example 4, ODE eq.(\ref{odep4}): $p_0^*(n,T)$ vs. $n$ (left) and $\log_{10}d_1^*(n,T)$ vs $n$ (right) for $T=1, 2,\ldots, 5$ (curves from cherry tone to black) for ODE (\ref{odep4}) with $m=0$. Dashed lines are fits (see text)} \label{fig28}  
\end{center}
\item {\bf The effect of the T-interval:} Our scheme works for any $T$-interval. However, the number of iterations needed to get a given precision depends on it. In general, there are two main aspects to remark. First, the minima's convergence rate slows down as we increase $T$ and second, the minima values, $p_0^*(n,T)$, seem to be proportional to $T$. These ideas are confirmed when we analyze $p_0^*(n,T)$ and $d_1^*(n,T)$ for $m=0,1$. For instance we show in figure \ref{fig28} their behaviour for the $m=0$ case (the results for the $m=1$ case are very similar qualitatively and quantitatively). We observe an exact scaling with $T$: $p_0^*(n,T)=Tp_0^*(n,1)$ and $d_1^*(n,T)=Td_1^*(n,1)$ for all $n$'s and $T$'s. For $m=1$ such scaling is not exact but it is dominant for large values of $n$. This permits us to predict how many iterative orders, $n(p,t)$ are needed for a given $T=10^t$ to reach a precision $d_1^*=10^{-p}$.  $n(p,t)$ is solution of the equation:
\begin{equation}
b_0(1)n(p,t)^3+(p+t+b_1(1))n(p,t)^2+b_2(1)n(p,t)+b_3(1)=0
\end{equation}
where $b_i(1)$ are the coefficients for the fitted function to the $d_1^*(n,T=1)$ data (see above). As an example we have computed $n(p,t)$ for $m=0$:
 \begin{center}
  \begin{tabular}{| c || c | c | }
    \hline
    {\boldmath $t$}&\multicolumn{2}{|c|}{\boldmath $n(p,t)$} \\ \hline
    &\phantom{  }{\boldmath $p=3$}\phantom{  }&\phantom{  } {\boldmath $p=8$}\phantom{  }\\ \hline
      \phantom{  }0\phantom{  } &21&365\\ \hline
    \phantom{  }1\phantom{  } &68&442\\ \hline
    \phantom{  }2\phantom{  } &138&519 \\ \hline
     \phantom{  }3\phantom{  } &213&595 \\ 
          \hline
  \end{tabular}
 \end{center}
 where we have used the values from the fitting: $b_0(1)=-0.0129$, $b_1(1)=-3.309$, $b_2(1)=13.672$ and $b_3(1)=-36.061$. Let us remark that: (1) for a given interval $T=10^t$, any increase of the precision needs of a large number of iterations due to the near one decay rate in this problem, and (2) for a given precision $p$, the number of iterations to reach such precision increases linearly with $t$. In fact, $n(p,t)\simeq t/|b_0(1)|$ when $t\rightarrow\infty$ for any precision $p$. This result is consistent with the idea that we need to get infinite terms from the recurrence to describe the solution along the real line with infinite precision. Finally, we see that the value of the decay coefficient $|b_0(1)|$ (that it is an ODE-depending parameter) is very relevant from a practical point of view for finite $T$-values. For instance, a  low value would imply many terms to get the prescribed precision or the reverse. In our example, $|b_0(1)|$ is very small, and after $30$ iterations, we get a precision of order $10^{-3}$ that is not very good when we compare with the results we obtained in the BVP examples above. In any case, even if the basic decay rate for $T=1$ is good enough, we need to design an algorithm to deal with problems that look for large values of $T$. 

\item {\bf Change of variables:} We introduce the change of variables $u(x)=x y(x)$ into the ODE (\ref{odep4}). The transformed ODE to be solved is:
\begin{equation}
u(x)''+x^{1-m}u(x)^m=0\quad,\quad u(0)=0\quad,\quad u'(0)=1\label{cv4}
\end{equation}
The extended ODE for $p_1=p_2=p_3=0$ is now written:
\begin{equation}
\left(p_0+\epsilon\left(1-p_0\right)\right)u''+\epsilon x^{1-m}u(x)^m=0\label{cv4ex}
\end{equation}
that we solve perturbatively with the same boundary conditions as in eq.(\ref{cv4}). 

First, we can show that for any $m$, the iterative solution of eq.(\ref{cv4ex}) obtained by our method when $\epsilon=1$ and $p_0=1$ (without any minimization) is just the Taylor's expansion of ODE's solution (\ref{cv4}) around $x=0$.   Let us prove this curious property. Let us take the original ODE (\ref{cv4}) and introduce the following transformation: $x=\epsilon^{1/2}\bar x$ and $u(x)=\epsilon^{1/2}\bar u(\bar x,\epsilon)$. The ODE becomes
\begin{equation}
\bar u''+\epsilon\bar x^{1-m}\bar u^m=0\quad,\quad \bar u(0)=0\quad,\quad\bar u'(0)=1
\end{equation}
that is just the extended ODE (\ref{cv4ex}) when $p_0=1$. This ODE can be solved order by order in $\epsilon$ and the perturbative solution is of the form
\begin{equation}
\bar u(\bar x,\epsilon)=\sum_{n=0}^\infty\tilde u_n(\bar x)\epsilon^n\label{per}
\end{equation}
Therefore the solution of the original equation should be:
\begin{equation}
u(x)=\epsilon^{1/2}\sum_{n=0}^\infty\tilde u_n(\epsilon^{-1/2}x)\epsilon^n\label{ux}
\end{equation}
for any $\epsilon$. That implies that the coefficients $u_n$ should be of the form: $\tilde u_n(\bar x)=a_n \bar x^{2n+1}$ with $a_n$ some constants in order to get rid of the $\epsilon$ dependence on the right hand side of eq. (\ref{ux}). Therefore, the  solution of the ODE (\ref{cv4}) can be written:
\begin{equation}
u(x)=\sum_{n=0}^\infty a_n x^{2n+1}
\end{equation}
that it is just its Taylor expansion around $x=0$. Finally, we observe that $\lim_{\epsilon\rightarrow 1}\bar u(\bar x,\epsilon)=u(x)$. Therefore the solution of our extended ODE (\ref{cv4ex}) with $p_0=1$ corresponds to the Taylor expansion of the solution around the origin in the limit $\epsilon=1$. We can easily compute the expansion (\ref{per}) using our scheme for $p_0=1$ and, for example, we obtain the following expressions:
\begin{eqnarray}
m=1& :& u(x)=x-\frac{x^3}{6}+\frac{x^5}{120}-\frac{x^7}{5040}+\frac{x^9}{362880}-\frac{x^{11}}{39916800}+O(x^{13})\nonumber\\
m=2& :& u(x)=x-\frac{x^3}{6}+\frac{x^5}{60}-\frac{11 x^7}{7560}+\frac{x^9}{8505}-\frac{97 x^{11}}{10692000}+O(x^{13})\nonumber\\
m=3&:& u(x)=x-\frac{x^3}{6}+\frac{x^5}{40}-\frac{19 x^7}{5040}+\frac{619 x^9}{1088640}-\frac{17117 x^{11}}{199584000}+O(x^{13})\label{Taylor}
\end{eqnarray}
As we will see, except for the $m=1$'s case, $p_0=1$  is not the value that minimizes the distance $d_1(n;p_0)$ at each $n$ iteration level. Therefore it is not the best approximation to the real solution for any $n$th-iteration and finite $T$-value.
\begin{center}
\vskip -0.5cm
\includegraphics[height=9cm]{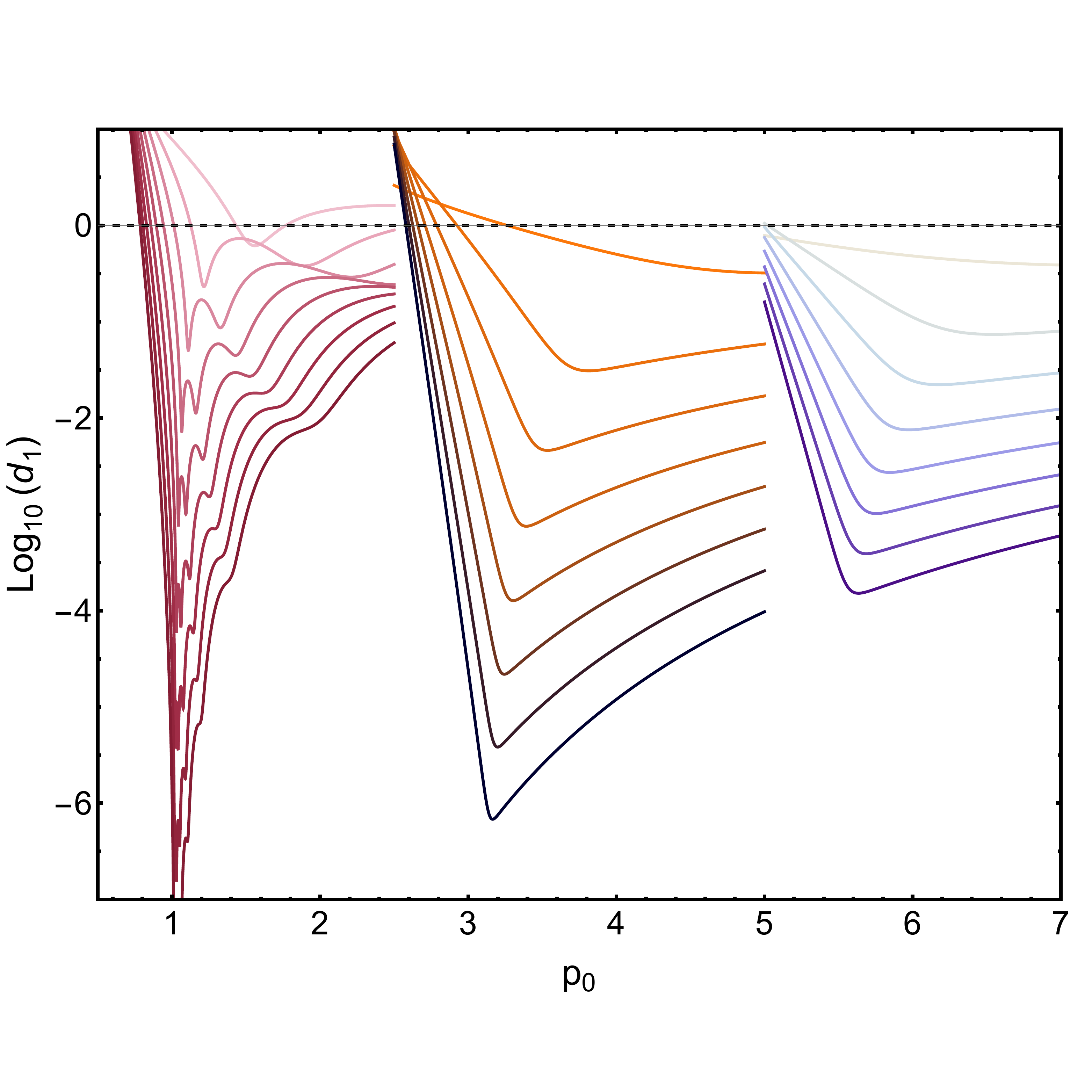}  
\vskip -0.5cm
\captionof{figure}{Example 4, ODE eq.(\ref{odep4}):  $\log_{10}d_1(n,p_0)$ vs. $p_0$ for the scheme applied to the ODE (\ref{cv4ex}) with $\epsilon=1$ and $T=5$. Set of curves at the left of figure (Valentine Tones):  $m=1$ and $n=2,3,\ldots 10$ from top to bottom. Set of curves at the center of figure (Rust Tones): $m=3$ and $n=2, 6,\ldots 30$ from top to bottom. Set of curves at the right of figure (Blue Tones):  $m=5$ and $n=2, 6,\ldots 30$ from top to bottom.} \label{fig29}  
\end{center}
We show in figure \ref{fig29} the behavior of $\log_{10}d_1(n;p_0)$ vs. $p_0$ for $m=1$, $3$ and $5$. We immediately observe that the minima are now deeper compared with the corresponding results for the original ODE (see for instance figs. \ref{fig26} and \ref{fig27}): around  to four times higher precision for a given iteration $n$. Let us remark that the distance for $u$'s is the same that the one for $y$'s, and therefore, any $u_n(x)$ having a smaller distance than the corresponding $y_n(x)$ is a better approximation.  That is, in this case, the change of variables improves very much the convergence behavior of our scheme. We do not show results for $m=0$ because it gives the exact result for $n=1$ and $p_0(n)^*=1$ with zero corrections to the exact result when $n>1$.
Finally, we show in figure \ref{fig29} how the distance for the case $m=1$ has many local minima that converge to the $p_0=1$ value as $n$ increases. This implies that the Taylor expansion is becoming a good approximation to the solution. On the contrary, in all other cases, the minima do not converge to one as we show in figure \ref{fig30}. 
\begin{center}
\vskip -0.5cm
\includegraphics[height=7cm]{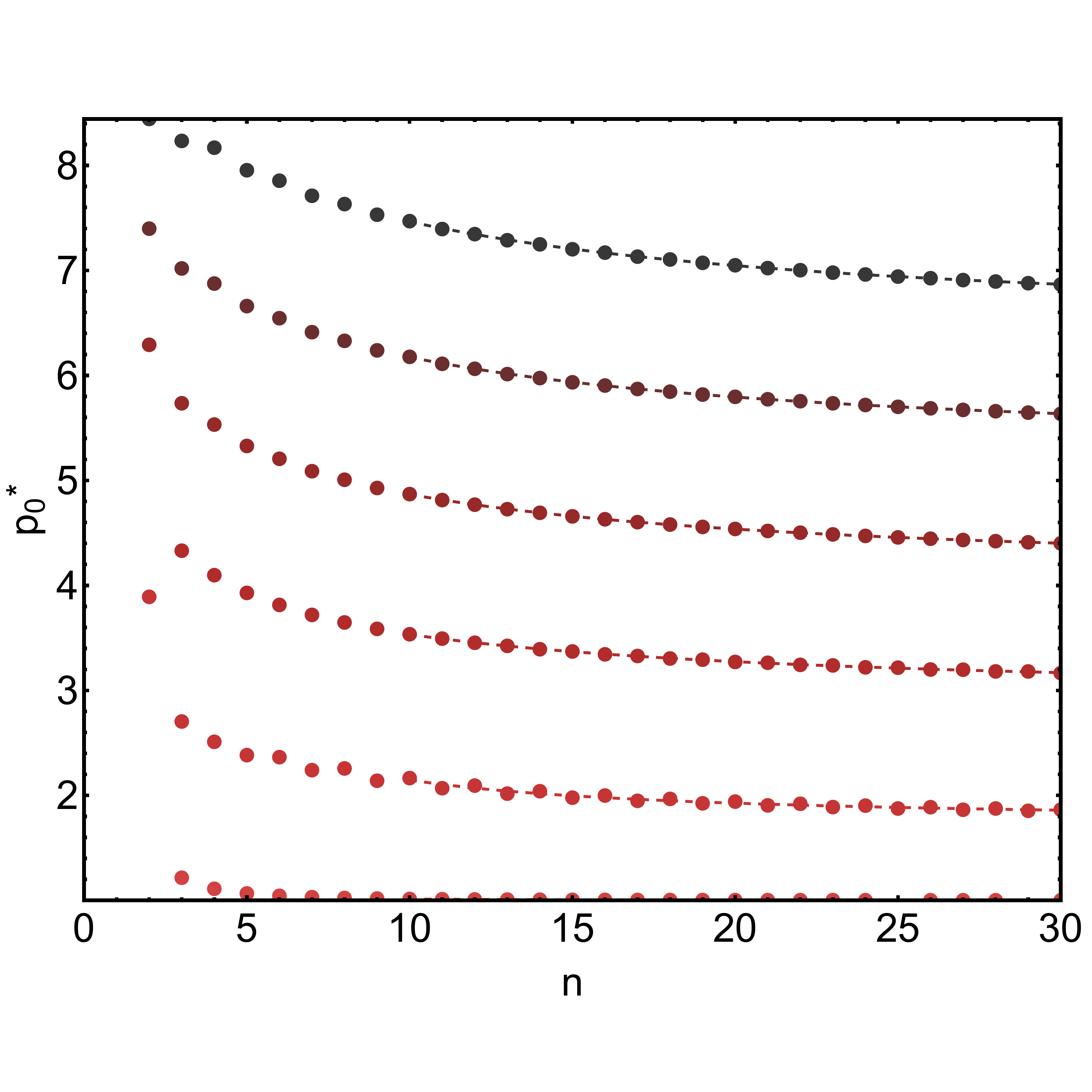}  
\includegraphics[height=7cm]{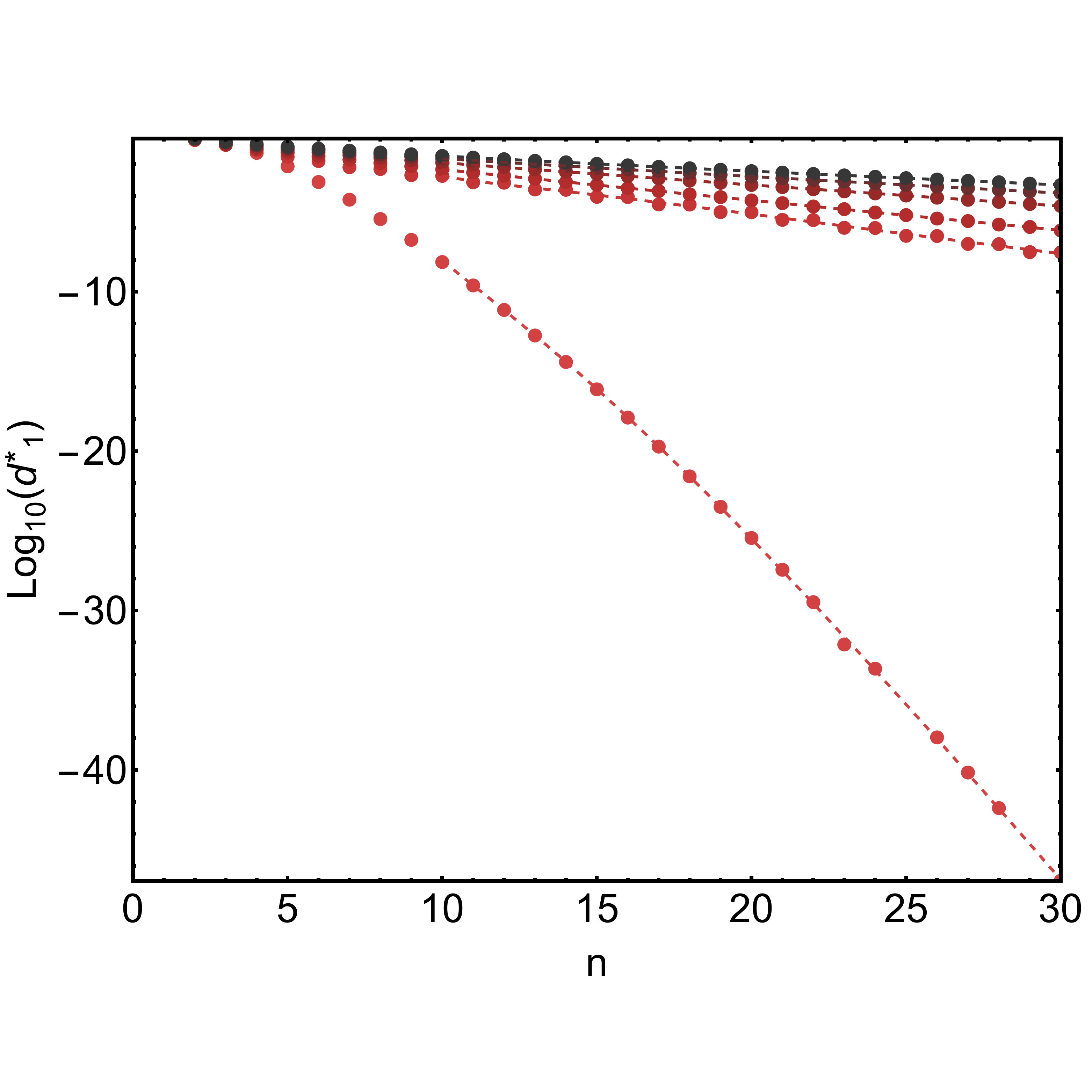}  
\vskip -0.5cm
\captionof{figure}{Example 4, ODE eq.(\ref{odep4}): $p_0^*(n,T)$ vs. $n$ (left) and $\log_{10}d_1^*(n,T)$ vs $n$ (right) for $m=1, 2,\ldots, 6$ (curves from cherry tone to black) for ODE (\ref{cv4}) with $T=5$. Dashed lines are fits (see text).} \label{fig30}  
\end{center}
We have fitted the data of figure \ref{fig30} in the same way as we did for the scheme applied to the original ODE (\ref{odep4}) (see above). The rate of convergence, $d_1^*(n)\simeq\delta^n$ and the asymptotic $p_0^*$ are now:
 \begin{center}
  \begin{tabular}{| c || c | c |}
    \hline
    {\boldmath $m$}&{\boldmath $\delta$}&{\boldmath $p_0^*$}\\ \hline
       1&0.003(0.003)&1.001(0.001)\\ \hline
    2 &0.7(0.3)&1.74(0.04) \\ \hline
       3&0.65(0.04)&2.93(0.01)\\ \hline
       4&0.749(0.002)&4.096(0.003) \\ \hline
         5&0.803(0.003)&5.272(0.005) \\ \hline
           6&0.838(0.003)&6.454(0.007)\\ \hline
    \hline
  \end{tabular}
 \end{center}
we see how the rate of convergence is higher than one for the original ODE (see similar table above). Now they range from almost zero for $m=1$ to $\simeq 0.8$ for $m=5$ or $6$ compared with $\simeq 0.9$ for all cases before. 
\begin{center}
\vskip -0.5cm
\includegraphics[height=7cm]{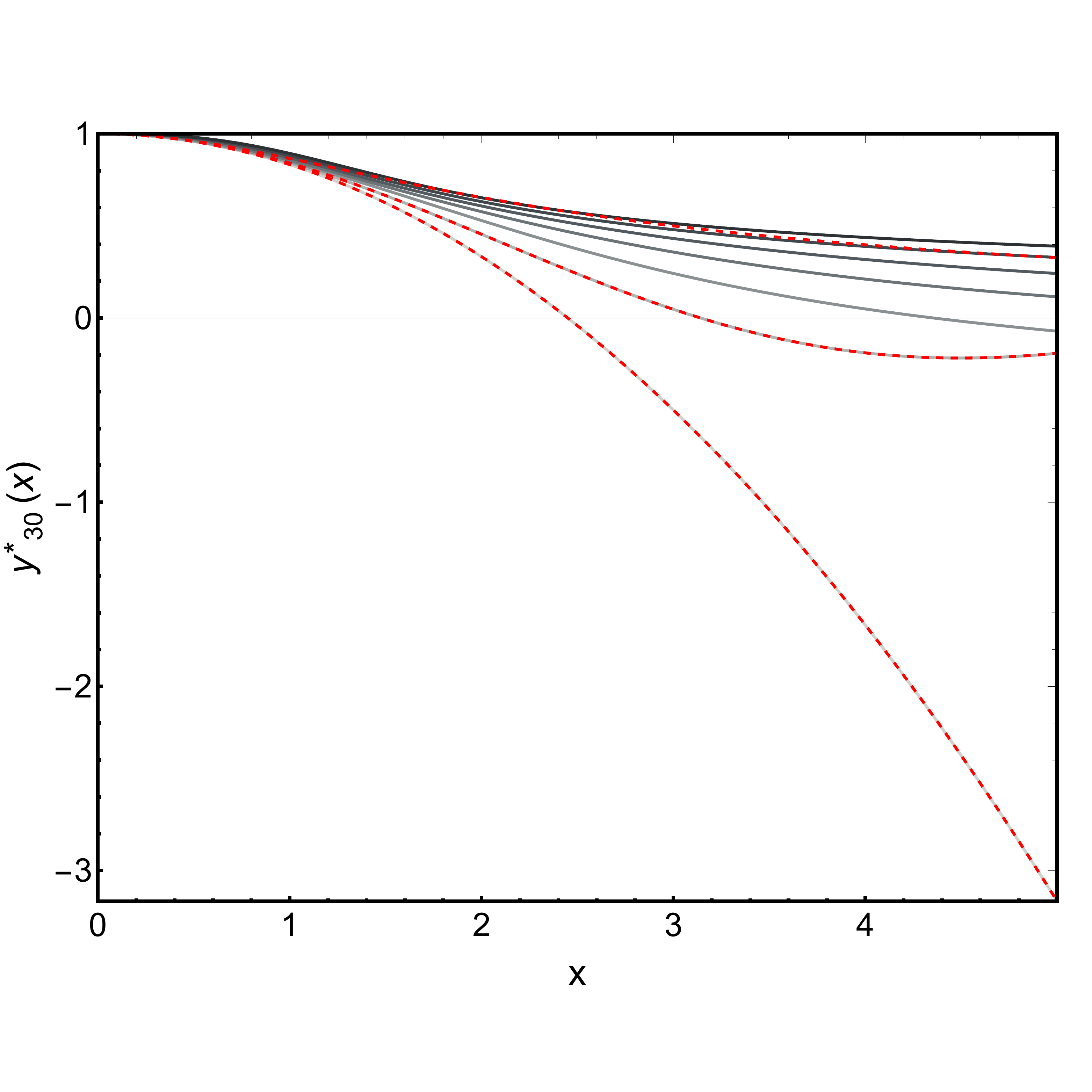}  
\includegraphics[height=7cm]{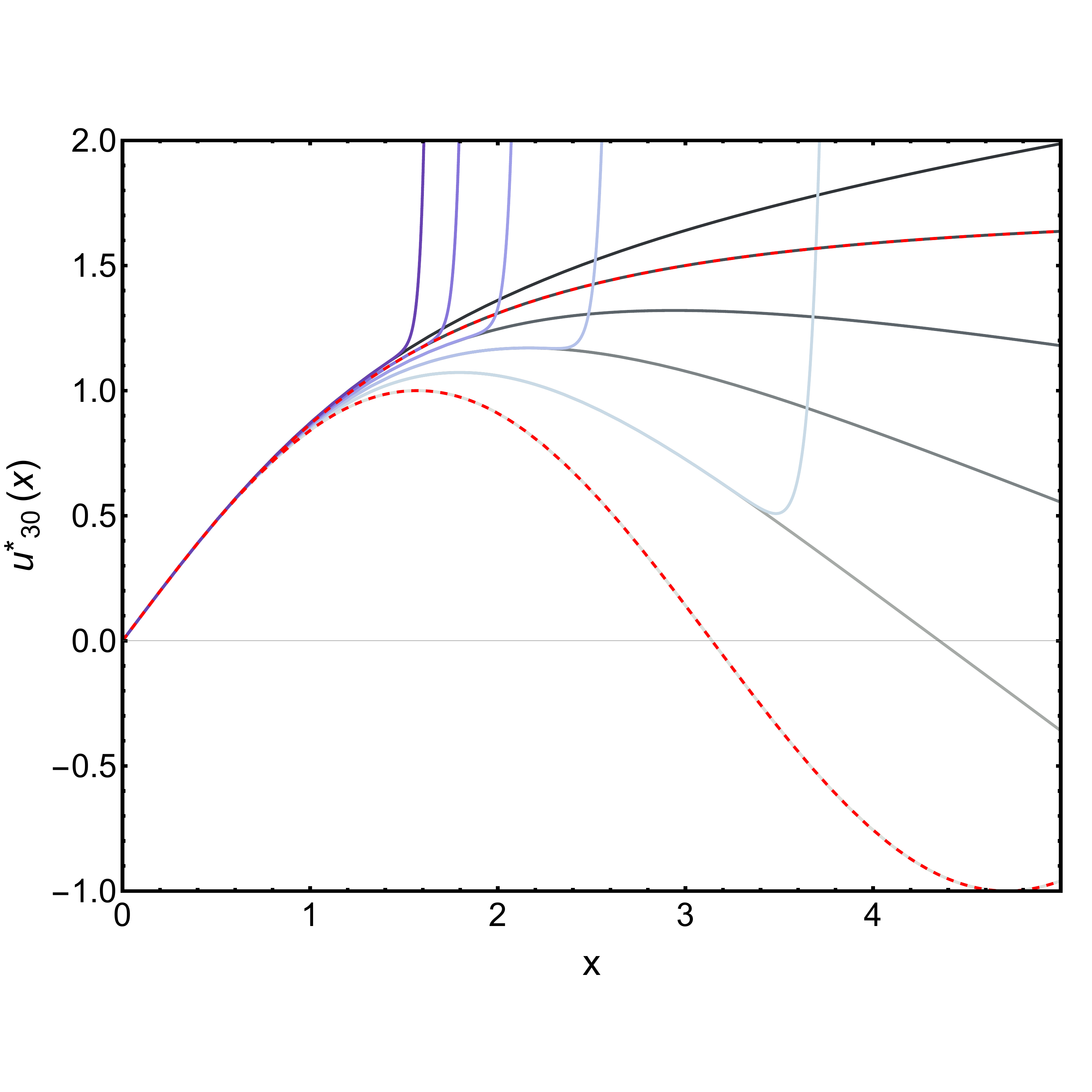}  
\vskip -0.5cm
\captionof{figure}{Example 4, ODE eq.(\ref{odep4}). Left: Approximate solution to the ODE (\ref{odep4}) $y_{30}(x)$ for $m=0,1,2,\ldots, 6$ from bottom to top. Right: Approximate solution to the ODE (\ref{cv4}) $u_{30}(x)$ for $m=1,2,\ldots, 6$ from bottom to top. In both cases $x\in[0,5]$. Dahsed red lines are the known exact solutions for $m=0, 1$ and $m=5$. Blue curves on the right figure are the Taylor expansions around $x=0$ of each solution up to order $x^{60}$. The Taylor expansion for the $m=1$ case is also drawn but it superimpose graphically the exact solution in this $x$-interval.} \label{fig31}  
\end{center}
We see in figure \ref{fig31} the approximate solutions at the iterative level $n=30$  for $m=0, 1,\ldots, 6$ when we use our scheme for the ODE (\ref{odep4}) and after the change of variables for the ODE (\ref{cv4}). All the curves have a smooth behavior and they follow the exact known results for $m=0,1$ and $5$. We only appreciate some deviations from the exact results when $m=5$ and $y_{30}(x)$ due to the slow convergence of our series in such case ($\delta\simeq 0.964$) with a distance $d_1^*\simeq 10^{-1}$ that corresponds to an average error on the solution of $\simeq 10^{-2}$. Nevertheless, after the change of variables such effect is much smaller because $\delta\simeq 0.838$, $d_1^*\simeq 10^{-4}$ with an average error with respect to the exact solution of around $10^{-8}$ when $m=5$. Finally, we plot the Taylor's expansion of the solutions around zero up to $x^{60}$  (\ref{Taylor}) that we found when $p_0=1$. We observe in figure \ref{fig31} how they systematically deviate from the solution for large enough $x$-values in the interval $[0,5]$ (except for the case $m=1$). This reflects the fact that traditional perturbation theories around the initial condition are typically exponentially unstable for large distances from it. 
 \item {\bf Ghost expansions:}
 We have seen again in this ODE  (\ref{odep4}) with parameter $m$ that their distance $d_1^*(n;m)$ decreases to zero exponentially fast in $n$. Therefore we may obtain the corresponding Ghost Expansion (\ref{ghost1}):
 \begin{equation}
 y(x;m)=x+\sum_{n=1}^\infty w_n(x;m)d_1^*(n;m)\label{gh1}
 \end{equation}
The first terms when $T=5$ are:
 \begin{eqnarray}
 w_1(x;m)&=&-\frac{x^3}{6p_0^*(1;m)d_1^*(1;m)}\nonumber\\
 w_2(x;m)&=&x^3\frac{20 p_0^*(2;m)^2-40 p_0^*(2;m) p_0^*(1;m)+20 p_0^*(1;m)+m p_0^*(1;m) x^2}{120 p_0^*(2;m)^2 p_0^*(1;m)d_1^*(2;m)}  \nonumber\\
 w_3(x;m)&=&-\frac{x^3}{15120 p_0^*(3;m)^3 p_0^*(2;m)^2d_1^*(3;m)} \biggl[2520 \biggl(p_0^*(3;m)^3 (1-2  p_0^*(2;m))\nonumber\\
 &+&3 p_0^*(3;m)^2  p_0^*(2;m)^2-3 p_0^*(3;m)  p_0^*(2;m)^2+ p_0^*(2;m)^2\biggr)\nonumber\\
 &+&126 m x^2 \left(p_0^*(3;m)^3-3p_0^*(3;m) p_0^*(2;m)^2+2  p_0^*(2;m)^2\right)\nonumber\\
 &+&m (8 m-5)  p_0^*(2;m)^2 x^4\biggr]
 \end{eqnarray}
 with 
 \begin{eqnarray}
 \{p_0^*(1;i)\}_{i=1}^m&=&(1.55999,4.20161,4.77566,5.26814,5.70278,6.09350)\nonumber\\
 \{p_0^*(2;i)\}_{i=1}^m&=&(1.21599,3.89134,4.33081,6.29169,7.39812,8.44295)\nonumber\\
 \{p_0^*(3;i)\}_{i=1}^m&=&(1.11073,2.70399,4.09761,5.73501,7.01895,8.23354)
 \end{eqnarray}
 and
 \begin{eqnarray}
 \{d_1^*(1;i)\}_{i=1}^m&=&(0.616388,0.699358,0.640784,0.602187,0.572778,0.548822)\nonumber\\
 \{d_1^*(2;i)\}_{i=1}^m&=&(0.231361,0.322503,0.158973,0.356377,0.383203,0.401424)\nonumber\\
 \{d_1^*(3;i)\}_{i=1}^m&=&(0.0503811,0.207376,0.0882926,0.191017,0.221962,0.245027)
 \end{eqnarray}
 Finally we show in figure \ref{fig32} the behavior of $w_n(x;m)$ for larger $n$ values (from $9$ up to $30$). 
 \begin{center}
\vskip -0.5cm
\includegraphics[height=5cm]{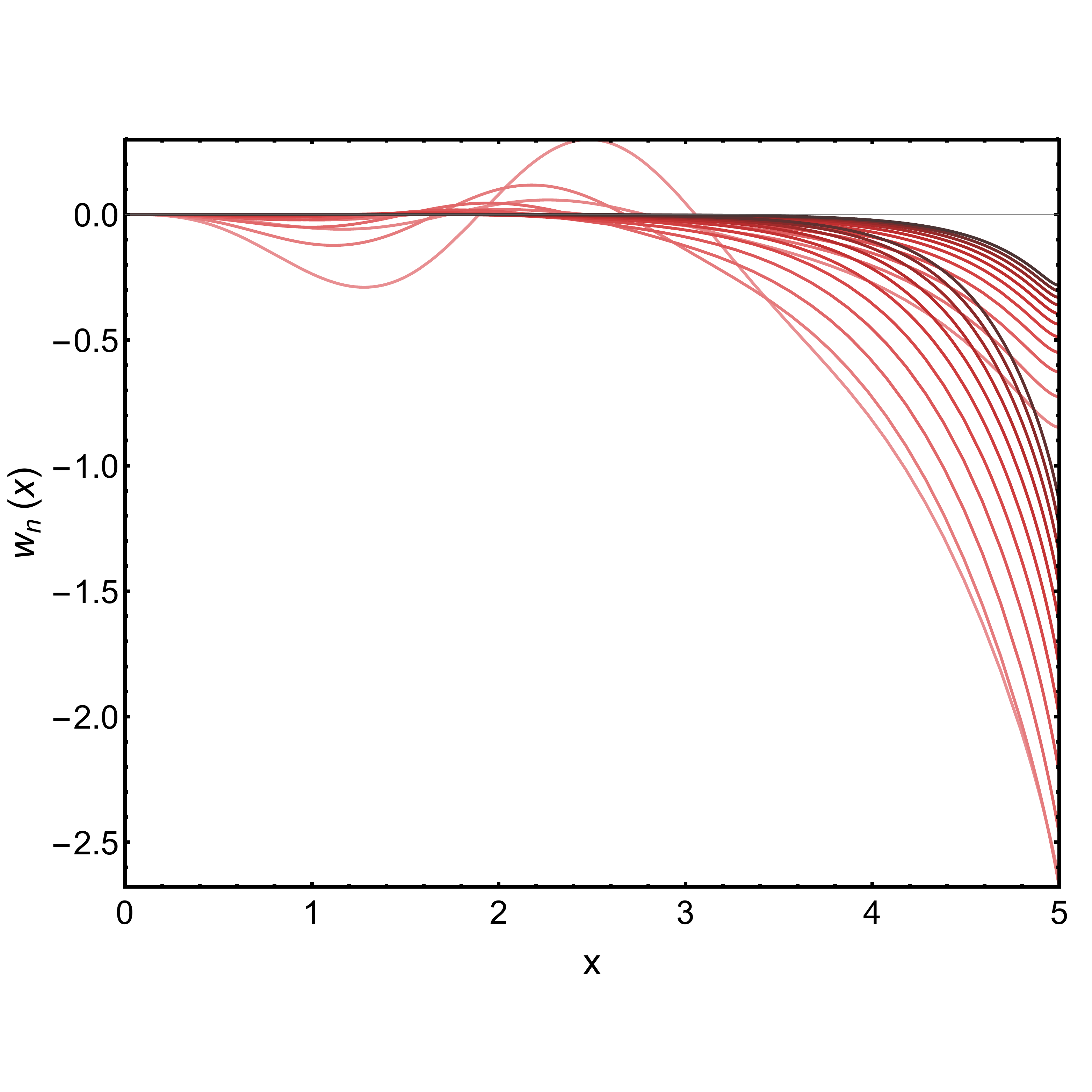}  
\includegraphics[height=5cm]{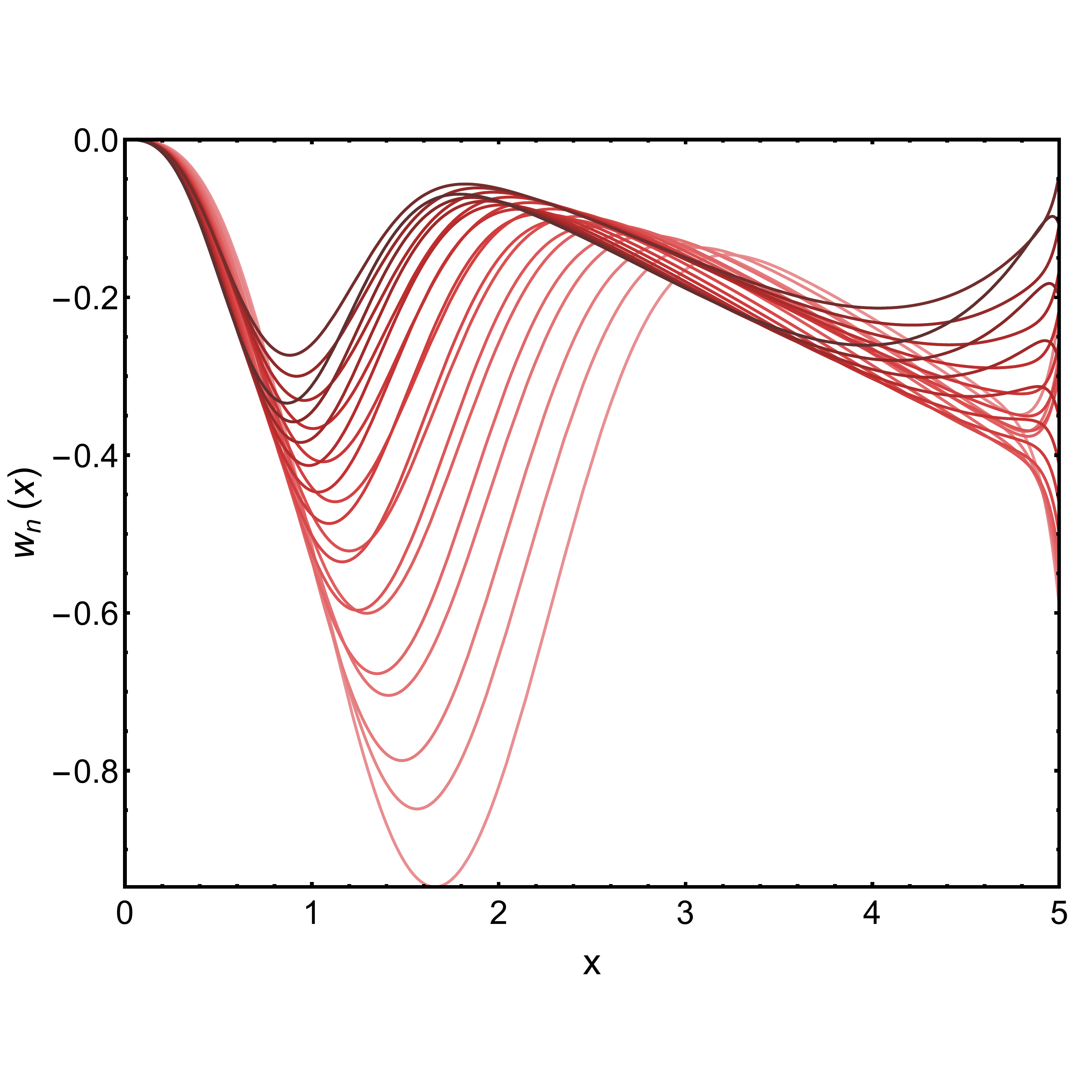}  
\includegraphics[height=5cm]{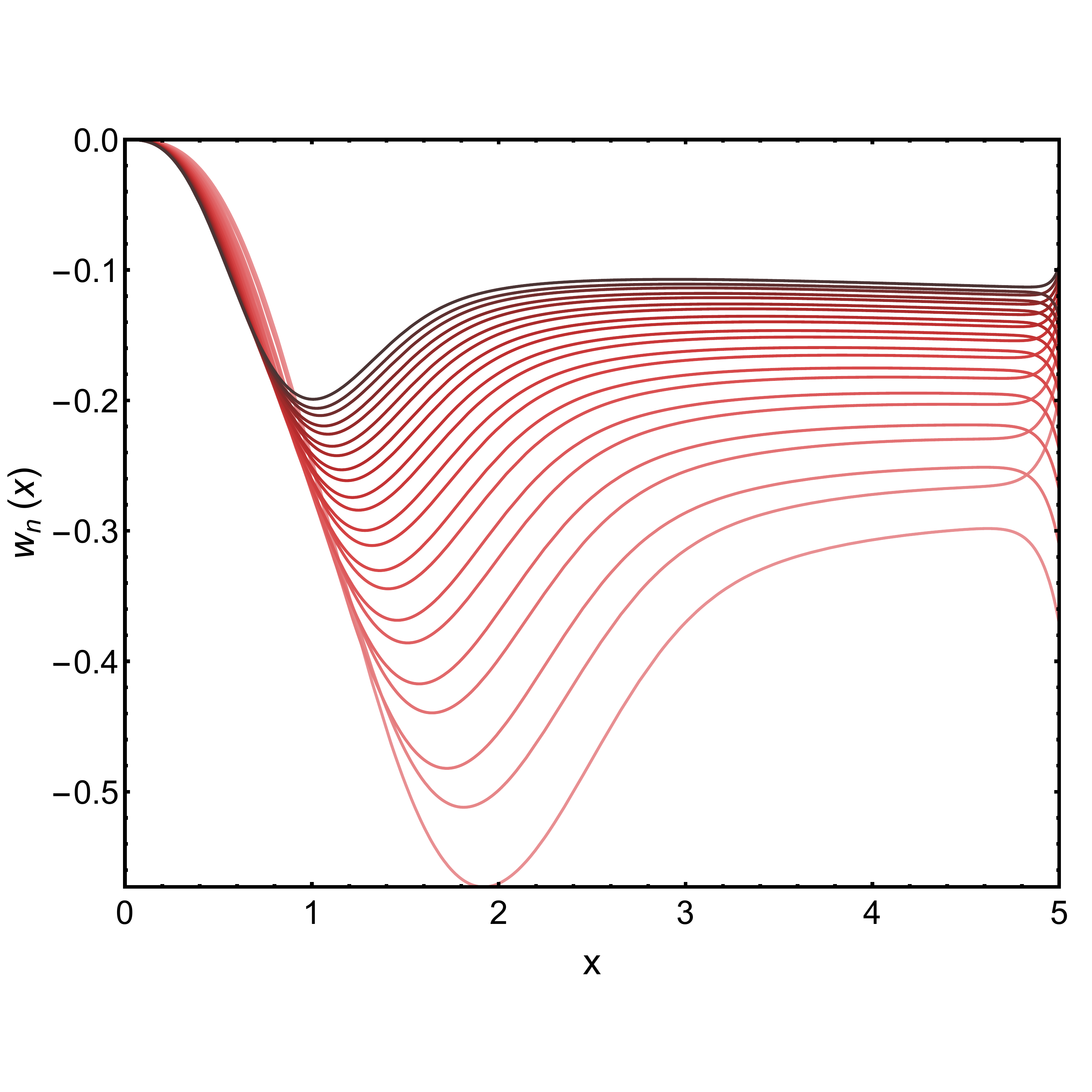}  
\vskip -0.5cm
\captionof{figure}{Example 4, ODE eq.(\ref{odep4}): Functions $w_n(x;m)$ corresponding to the ghost expansion (\ref{gh1}) from $n=9$ to $n=30$ (light to dark cherry-tones respectively) and $m=2$, $m=3$, $m=4$ (left, center and right figures).} \label{fig32}  
\end{center}
 \item {\bf An algorithm for large intervals:} We commented above that the distance to the exact solution behaves as $d_1(n;T)\simeq T\delta^n$ where $T$ is the length of the $x$-inteval where the $n$'th approximation is computed. From a practical point of view one can deal in a computer with a maximum finite $n_{max}$ (actually $n_{max}\simeq 50$  in a laptop). Therefore, our method can reach a given precision up to a maximum $T_{max}\simeq d_1\delta^{-n_{max}}$ and not beyond it. We can surpass this limitation by designing a multiple-interval algorithm in which we apply the method recurrently in time intervals smaller than $T_{max}$. We propose the following very simple algorithm:
 \begin{itemize}
 \item (0) Obtain the $n$-th approximation, $y_n(x;p,x_c,\bar y_0,\bar y_1)$ with generic boundary conditions $y(x_c)=\bar y_0$ and $y'(x_c)=\bar y_1$ for a chosen $n$ that is fixed all over the algorithm. This is an algebraic step that is done only once.
 \item (1) Fix the initial interval $T$, the desired maximum target distance $d_{max}$ and the initial boundary conditions at $x_c=0$.
 \item (2) Compute the optimal set of parameters $p^*$ that minimize the distance $d_1(p,T)$ for the given boundary conditions at $x_c$.
 \item (3) If $d_1(p^*,T)>d_{max}$ (the distance is larger than our precision goal) then we reduce $T$: $T=3T/4$ and go to step (2).
 \item (4) We have reach the precision goal on the interval. Therefore we look for the initial conditions for the next iteration from the known algebraic solution $y_n(x;p^*,x_c,\bar y_0,\bar y_1)$. It is convenient, if possible, to look for $x^*$ nearest to $x_c+T$ such that $N(x)y_n(x;p^*,x_c,\bar y_0,\bar y_1)\vert_{x=x^*}=0$. That guaranties that the derivatives at $x^*$ are correctly related by the original ODE. Then, for the next iteration: $\bar y_0=y_n(x^*;p^*,x_c,\bar y_0,\bar y_1)$, $\bar y_1=y_n'(x^*;p^*,x_c,\bar y_0,\bar y_1)$ and  $x_c=x^*$. We also increment a little the $x$-interval: $T=9T/8$. Repeat the sequence by going to step (2) until it is reached the  desired total $x$-interval.
 \end{itemize}
 In this way we get an algebraic piecewise approximate solution with a distance smaller than $d_{max}$ to the true solution. Observe that the overall solution is continuous with first and second continuous derivatives. The unique source of accumulative error is on the computation of the initial condition at each interval. How it propagates is an open question that should be studied. In any case, we think that there is a large space for the improvement of this algorithm by applying aready well known optimizing strategies. 
 
 We have applied this algorithm to the transformed Lane-Emdem ODE (\ref{cv4}). 
 We see in figure \ref{fig33} how the algorithm behaves for a large $x$-interval ($T=20$) using $u_5(x)$ as approximate functions. The algorithm adapts the local interval to the average desired distance $d_{max}$. That is, each component of the picewise function has a $d_1^*$ distance smaller than $d_{max}$. For large $x$-variations of $u(x)$, the local $T$-interval is smaller to reach the desired precision and vice versa. Moreover, the number of intervals grows for smaller values of $d_{max}$, and it is expected that they diminish when increasing $n$. Finally, we see how the obtained piecewise solutions differ from the known exact one. For $d_{max}=10^{-6}$ and $10^{-8}$ the difference is maintained around $10^{-8}$ and $10^{-10}$ when $x<10$ in both cases respectively. However, they begin to grow for $x>10$. Similar behavior is observed when we solve numerically the ODE by using the standard DSolverValue routine from Mathematica with default parameters. It is remarkable that the case $d_{max}=10^{-12}$ manages to maintain the precision around $10^{-13}$ all over the interval. Similar behavior is found for the cases $m=2$ and $m=4$ where the exact solution is not known (see figure \ref{fig34})
  
 \begin{center}
\vskip -0.5cm
\includegraphics[height=5cm]{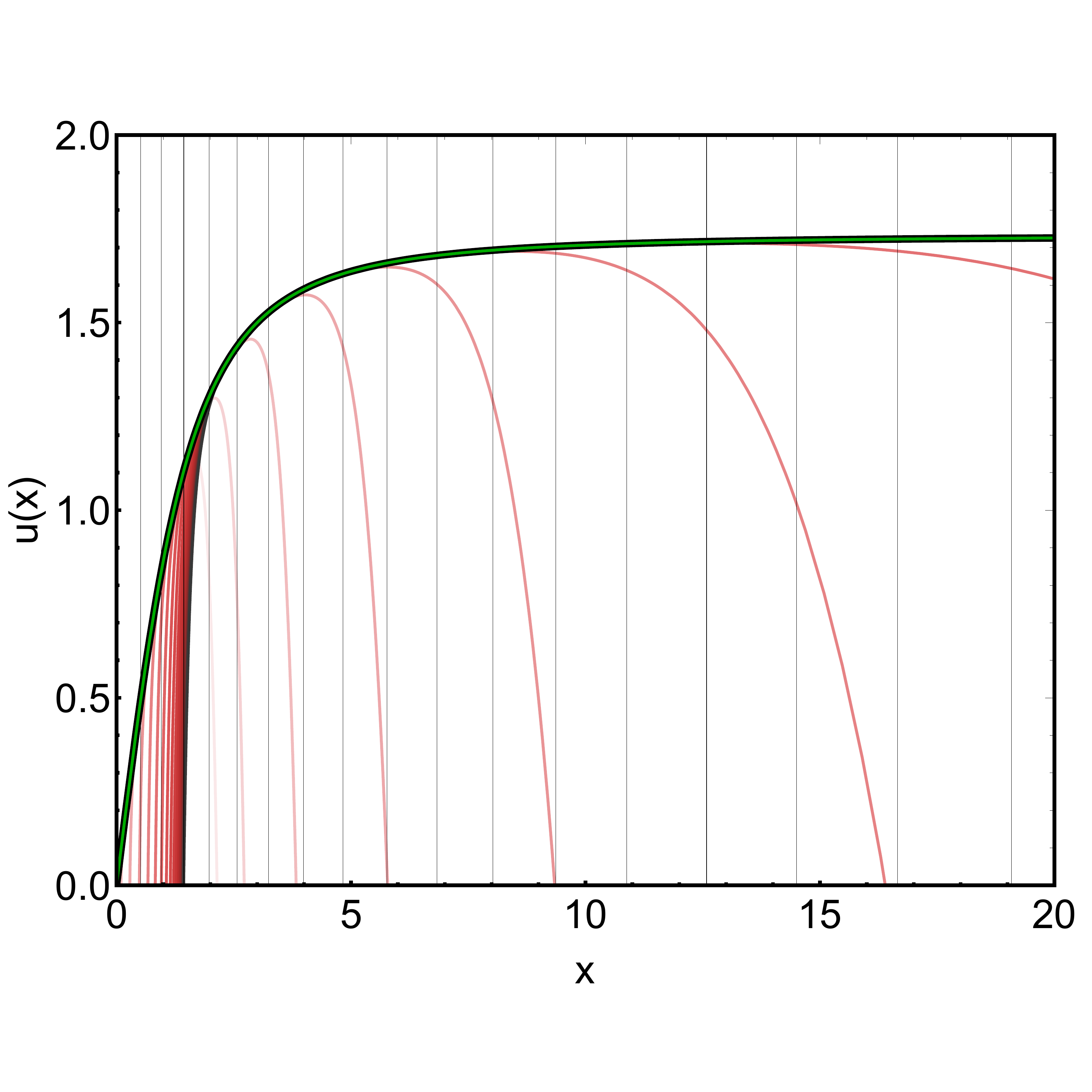}  
\includegraphics[height=5cm]{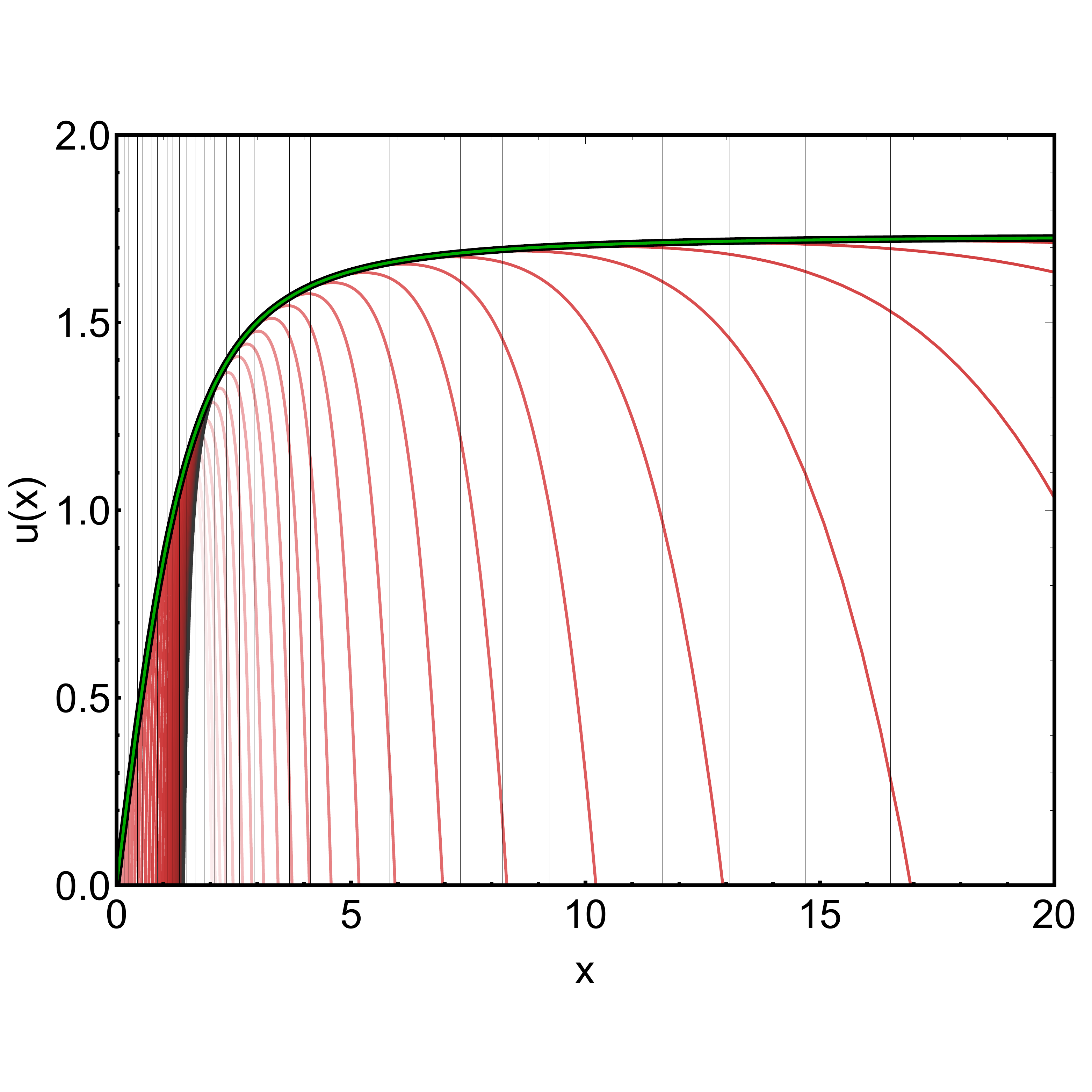}  
\includegraphics[height=5cm]{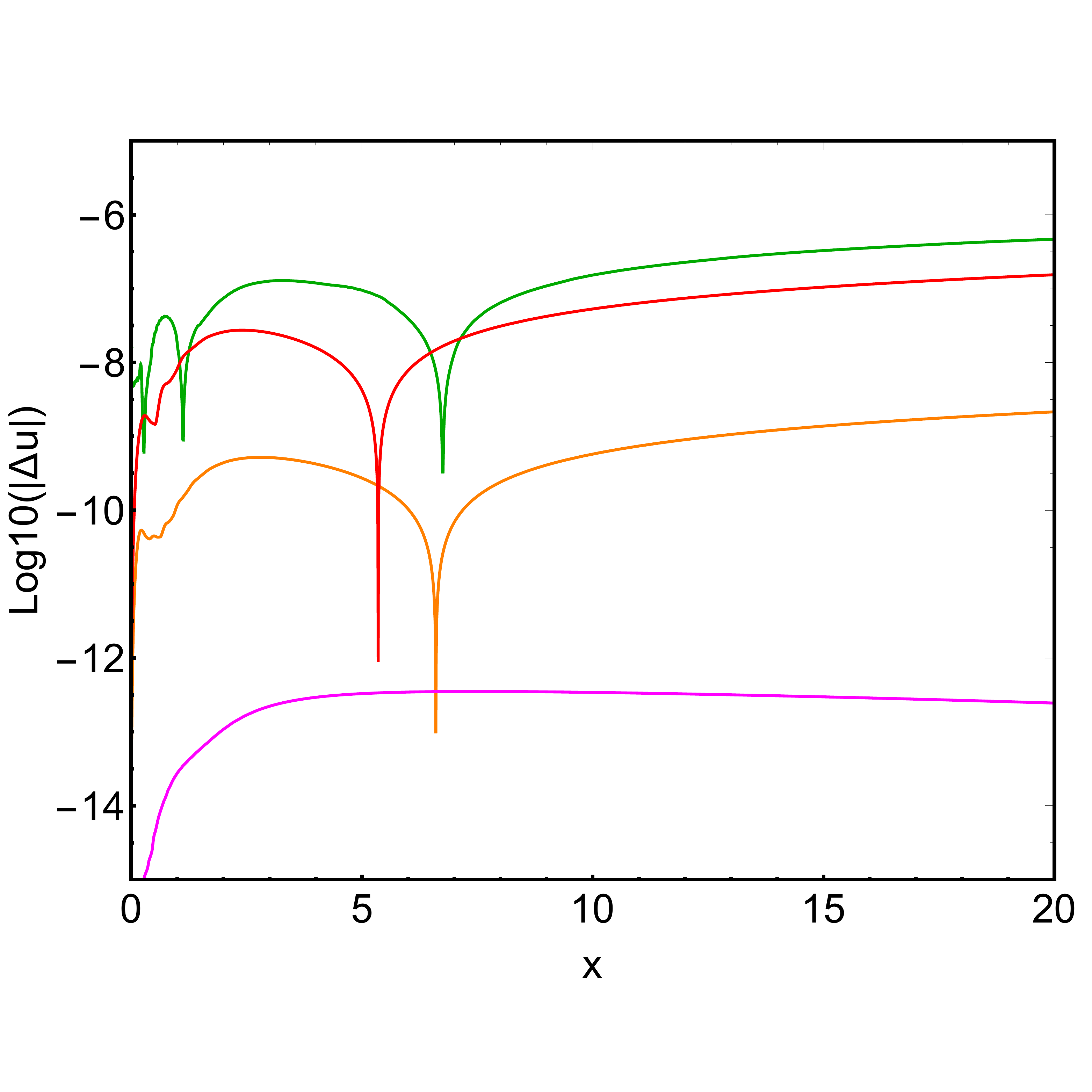}  
\vskip -0.5cm
\captionof{figure}{Example 4, ODE eq.(\ref{odep4}). Left: Solution of the transformed Lane-Emdem ODE (\ref{cv4}) with $m=5$ by using the multiple-interval algorithm with  described in the main text with $n=5$. Black-thick curve is the approximation obtained.  Vertical thin lines define the intervals defining the piecewise solution from our algorithm and the Cherry-tone curves are the extended local functions from the piecewise approximation. Green curve is the numerical solution obtained by  Mathematica's DSolveValue routine. Left: $d_{max}=10^6$, Center: $d_{max}=10^{-12}$, Right: $\log_{10}\vert u_5(x)-u(x)\vert$ where $u(x)=x/(1+x^2/3)^{1/2}$ is the known exact solution. Red, Orange and Magenta curves are for $d_{max}=10^{-6}$, $10^{-8}$ and $10^{-12}$ respectively.} \label{fig33}  
\end{center}

 \begin{center}
\vskip -0.5cm
\includegraphics[height=7cm]{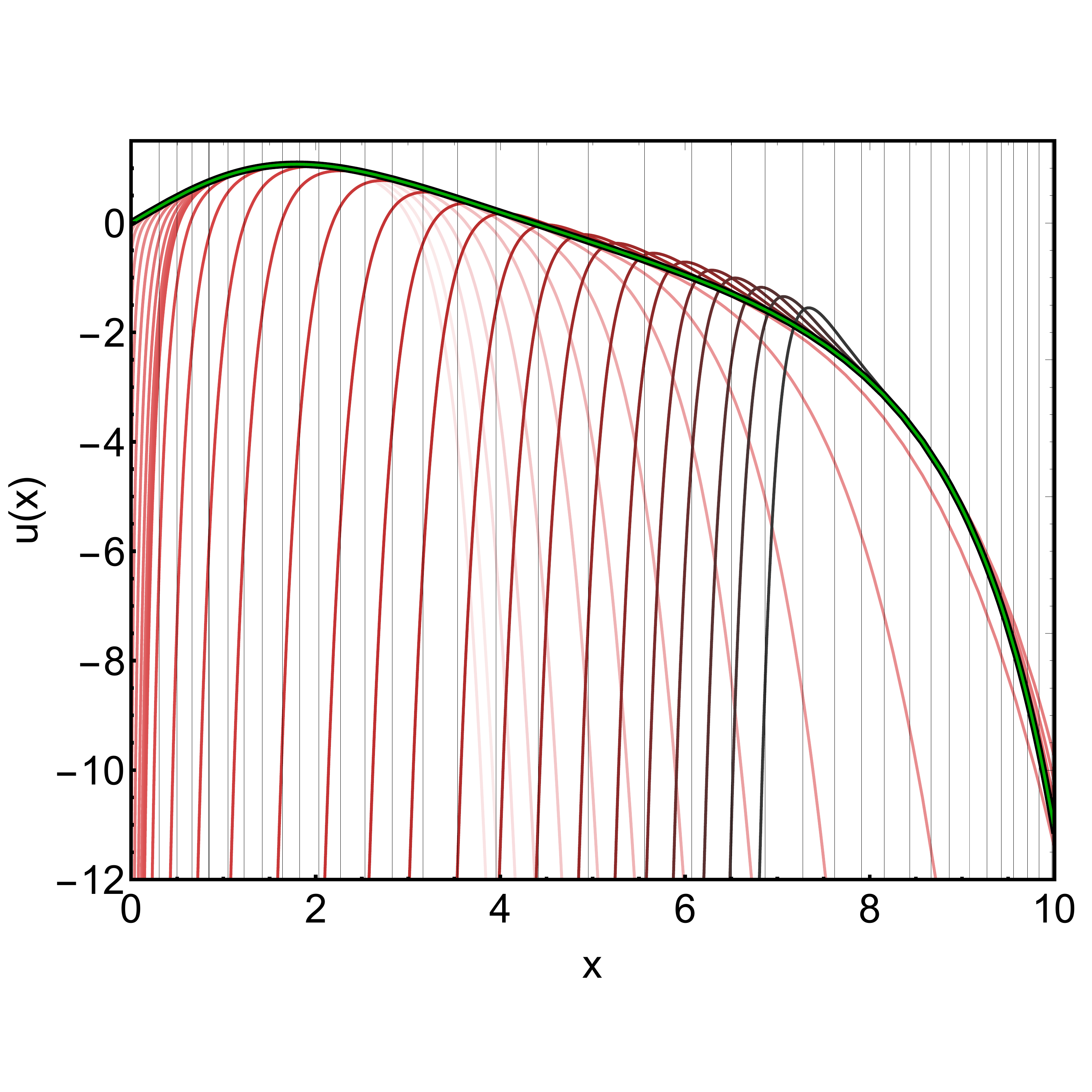}  
\includegraphics[height=7cm]{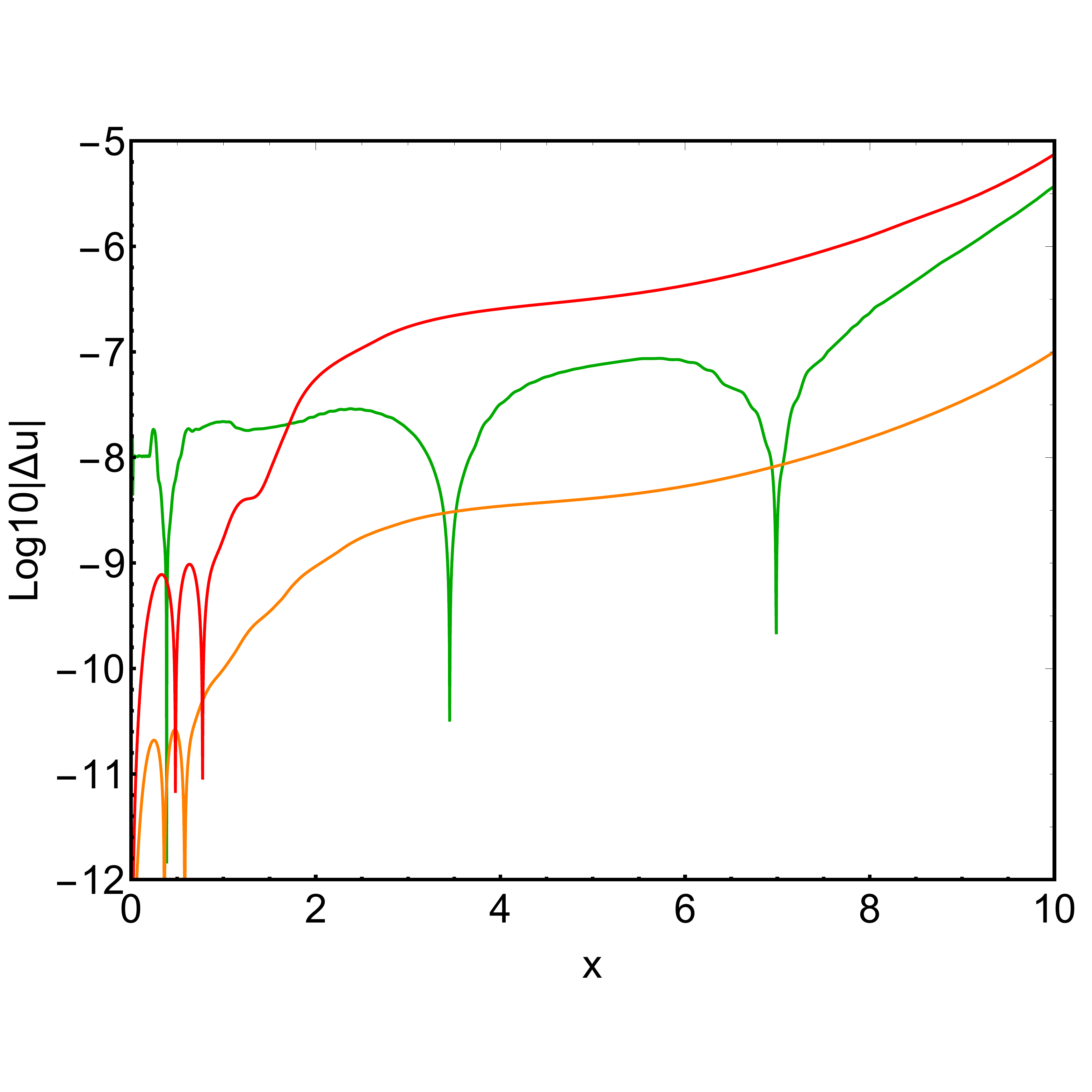}  
\includegraphics[height=7cm]{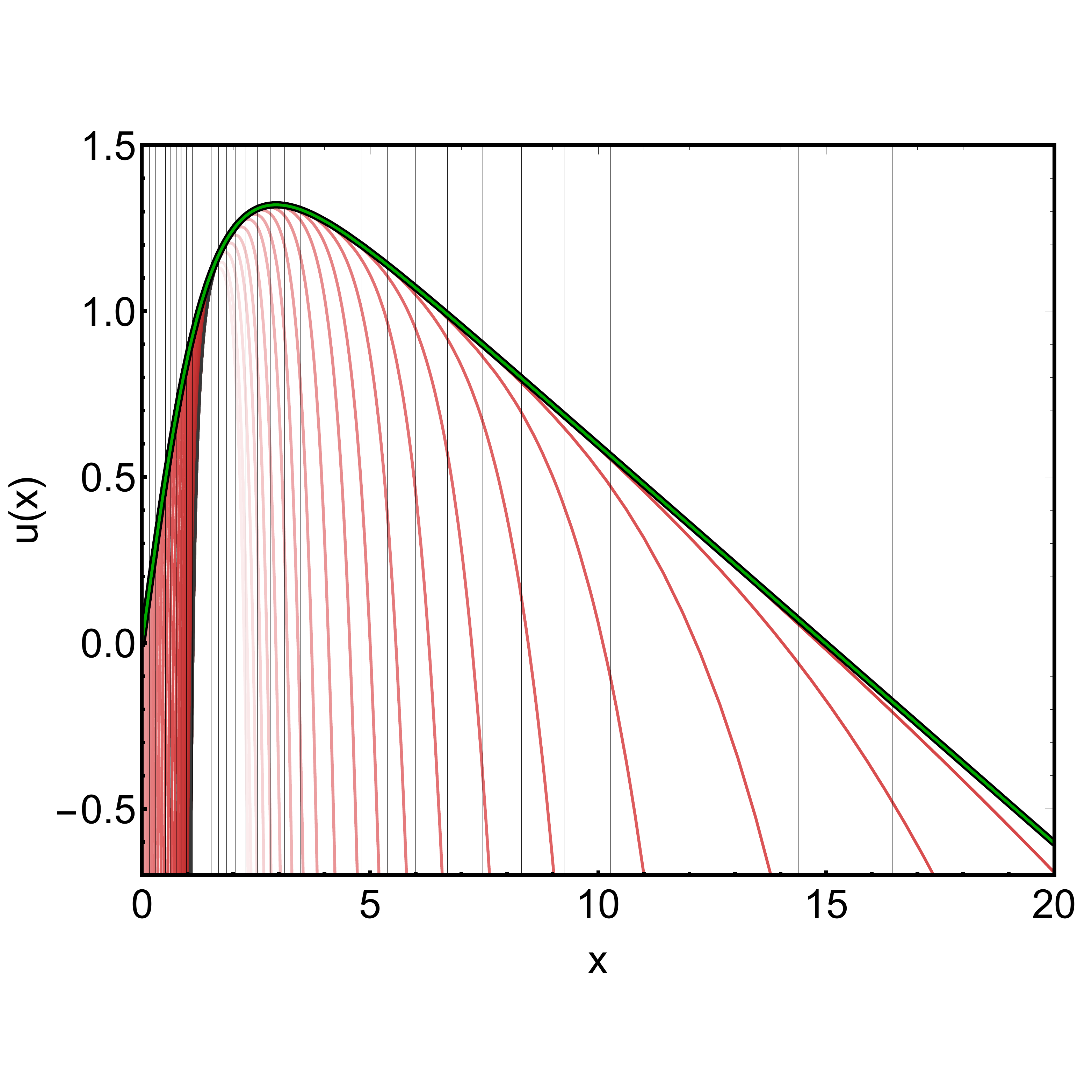}  
\includegraphics[height=7cm]{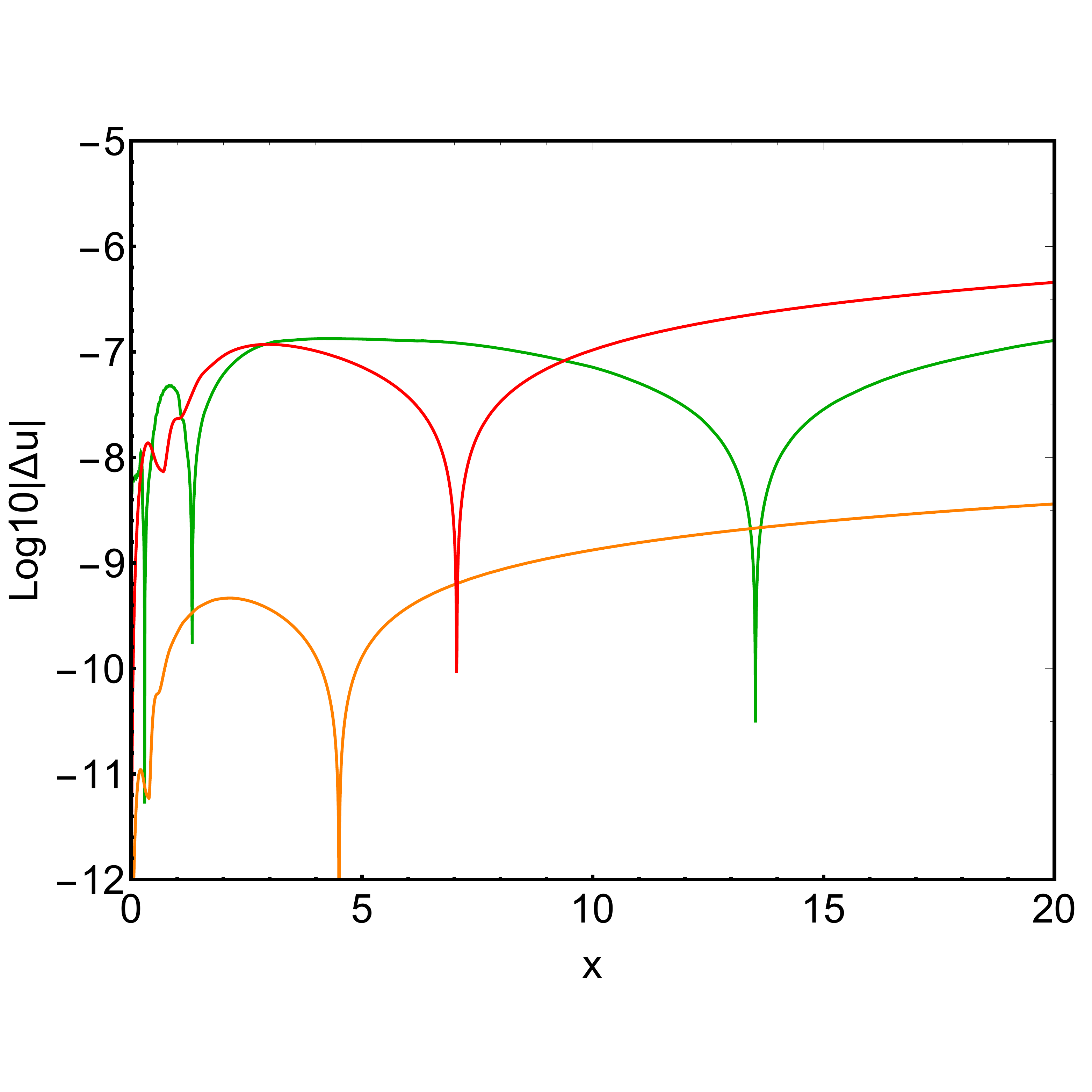}  
\vskip -0.5cm
\captionof{figure}{Example 4, ODE eq.(\ref{odep4}). Left: Solution of the transformed Lane-Emdem ODE (\ref{cv4}) with $m=2$ and $m=4$ (first and second row respectively) by using the multiple-interval algorithm with  described in the main text with $n=5$. Black-thick curve is the approximation obtained.  Vertical thin lines are the intervals defining the piecewise solution from our algorithm and the Cherry-tone curves are the extended local functions from the piecewise approximation. Green curve is the numerical solution obtained by  Mathematica's DSolveValue routine. Left:  $d_{max}=10^{-12}$, Right: $\log_{10}\vert u_5(x)-u_{ref}(x)\vert$ where $u_{ref}(x)$ is the approximation when $d_{max}=10^{-12}$. Red, and Orange  are for $d_{max}=10^{-6}$ and $10^{-8}$ respectively.} \label{fig34}  
\end{center}
\end{itemize}

\section{Conclusions}

We propose an algebraic method to approach the solution/s of ODEs with any boundary conditions. It is based on three main elements: (1) The definition of an extended ODE is composed of a linear generic differential operator that depend on few free parameters plus an $\epsilon$ formal perturbation of it formed by the original ODE minus the same linear term. (2) The assumption of a formal $\epsilon$ expansion of the solution and its application to the extended ODE to solve, order by order in $\epsilon$ the corresponding linear differential equation and (3) The fixing of the best parameter set by minimising a specified distance to the exact solution. As we commented in the introduction, there are other algebraic algorithms with some similarities to the one presented here. However, this is the first one that introduces parameters on the Linear Operator. This fact is crucial because it permits that the sequence of approximations is very sensitive to their values. Therefore, minor variations on them make the possibility that the approximations explore the space of possible solutions more efficiently. 

We have shown a set of typical examples where we check the feasibility of this scheme. Moreover, we highlight some exciting properties associated with our method: (1) The algorithm is algebraic and, therefore, is free of numerical errors. Its unique limitation is the necessity to make explicitly the integrals that define the algorithm's recurrence. We aimed to keep the algorithm under such premises, and more thought is needed to include problems where we cannot do such integrals. (2) The distance to the exact solution at a given approximation order depends on the generic lineal operator's parameters. It presents a set of minima associated with the ODE's number of solutions. That is an exciting tool to know the possible solutions, especially in Boundary Value Problems. (3) The sequence of regular functions reaches the ODE solution/s exponentially fast. This property permits us to define a formal solution expansion (Ghost Expansion) that we can use as a systematic perturbation scheme in other theories. Let us remark that this may be of great interest in cases where there is no intrinsic perturbation parameter that permits access to non-linear behaviours from initial linear approximation.  (4) The method is scalable for Initial Value Problems with long intervals, and we can improve it by using many known predictor-corrector methods.

We think that the proposal we present in this paper merits being deeply studied in many aspects. However, the authors are not specialists in these matters, and at this point, we think that colleagues with expertise in ODEs and algorithms should be able to develop this proposal. For instance, we believe it is essential to get some rigorous insight into the series convergence or the distance function's behaviour and how many minima they may have. Finally, this method has many exciting applications, for instance, to systems of first-order ODEs, partial differential equations, etc., that we will explore in future works.

\section{Acknowledgments}
This work is part of the Project of I+D+i Ref. PID2020-113681GB-I00, financed by 
MICIN/AEI/10.13039/501100011033 and FEDER “A way to make Europe”.

\end{document}